\documentclass[10pt]{article}
\usepackage{amsmath,amscd}
\usepackage{amssymb,latexsym,amsthm}
\usepackage{color}
\usepackage[dvipsnames]{xcolor}
\usepackage[spanish,english,activeacute]{babel}
\usepackage[pdftex]{hyperref}
\usepackage[utf8]{inputenc}
\usepackage{amssymb}
\usepackage{amsmath,amsthm,amscd}
\usepackage{lscape}
\usepackage{fancyhdr}
\usepackage{amsfonts}
\usepackage{pb-diagram}

\numberwithin{equation}{section}
\newtheorem{theorem}{Theorem}[section]

\newtheorem{definition}[theorem]{Definition}
\newtheorem{proposition}[theorem]{Proposition}
\newtheorem{lemma}[theorem]{Lemma}

\newtheorem{corollary}[theorem]{Corollary}

\theoremstyle{definition}

\newtheorem{example}[theorem]{Example}

\newtheorem{remark}[theorem]{Remark}

\deactivatequoting

\title{\textbf{Noncommutative coding theory \\ and algebraic sets for skew $PBW$ extensions}}
\author{Oswaldo Lezama\\
\texttt{jolezamas@unal.edu.co}\\
Claudia Gallego\\
\texttt{cmgallegoj@unal.edu.co}
\\Seminario de Álgebra Constructiva - SAC$^2$\\ Departamento de Matemáticas\\ Universidad Nacional de
Colombia, Sede Bogot\'a}
\date{}
\begin{document}
\maketitle
\begin{abstract}
\noindent The classical commutative coding theory has been recently extended to noncommutative rings of polynomial type. There are many interesting works in coding theory over single Ore extensions. In this review article we present the most relevant algebraic tools and properties of single Ore extensions used in noncommutative coding theory. The last section represents the novelty of the paper. We will discuss the algebraic sets arising in noncommutative coding theory but for skew $PBW$ extensions. These extensions conform a general class of noncommutative rings of polynomial type and cover several algebras arising in physics and noncommutative algebraic geometry, in particular, they cover the Ore extensions of endomorphism injective type and the polynomials rings over fields.

\bigskip


\bigskip

\end{abstract}

\tableofcontents
%-------------------------------------------------
%-------------------------------------------------

\section{Introduction}

The algebraic tools involved in the study of commutative coding theory have been recently extended to noncommutative skew cyclic codes. The main algebraic object in the classical commutative cyclic coding theory is the $\mathbb{F}$-algebra $\mathbb{F}[x]/\langle x^n-1\rangle$, where $\mathbb{F}$ is a finite field and $\langle x^n-1\rangle$ is the ideal generated by $x^n-1$ (see for example \cite{Gabidulin}, \cite{Huffman}, \cite{MacWilliams}). In the skew cyclic coding theory this object is replaced by $A/\langle x^n-1\rangle$, where $A$ is the Ore extension $\mathbb{F}[x;\sigma]$, with $\sigma$ an automorphism of $\mathbb{F}$, o more general, by $A/Af$, where $A=\mathbb{F}[x;\sigma,\delta]$, with $\delta$ is a $\sigma$-derivation of $\mathbb{F}$ and $f$ a polynomial of $A$ (see \cite{Boucher}). Even more, the field $\mathbb{F}$ can be replaced by an arbitrary ring $R$ (see \cite{Luerssen}).

In this review paper we have in mind two main purposes: To review the basic algebraic tools of the noncommutative coding theory over skew polynomial rings and, using this information, to investigate analogue tools for skew $PBW$ extensions. These extensions conform a general class of noncommutative rings of polynomial type and cover several algebras arising in mathematical physics and noncommutative algebraic geometry, in particular, they cover the Ore extensions of injective type. The focus of this work is algebraic and the paper not should be understood as a contribution to noncommutative coding theory. The novelty of the paper is only in the last section where we extend to skew $PBW$ extensions some results about algebraic sets of Ore extensions arising in coding theory. 

The paper is organized as follows: In the second section we will review the most important algebraic tools and the general theory the of skew cyclic codes, i.e., the coding theory over single Ore extensions. In particular, we will consider the following aspects of the classical coding theory, but in the noncommutative context: linear codes, skew cyclic codes, dual codes, generator matrix, parity check matrix, length and dimension of skew codes, similarity of polynomials, Vandermonde and Wronskian matrices, skew algebraic sets, bound of a polynomial, remainder and evaluation codes, maximum rank distance, minimal Hamming distance, linearized polynomials, among others. The structure and properties of skew cyclic codes depend on $R,\sigma, \delta $ and $f$. The most remarkable case with practical applications is when $R=\mathbb{F}$ is a finite field and $\delta=0$ (see \cite{Boucher}, \cite{Gomez-Torrecillas3}, \cite{Luerssen}, \cite{Gomez-Torrecillas}, \cite{Gomez-Torrecillas2}, \cite{Gomez-Torrecillas4}). Another key case is when $R=\mathbb{F}(t)$, with $\mathbb{F}$ a finite field and $\delta=0$. We will study first the general theory and tools of the skew cyclic codes for the case when $R=\mathbb{F}$ is a finite field and $A:=\mathbb{F}[x;\sigma,\delta]$, and after this, we will study the mentioned remarkable particular cases.

In the third section we will review some basic facts about the skew $PBW$ extensions that will use in the last section. The last section is dedicated to extend to skew $PBW$ extensions some of the results presented in the previous sections,  
more precisely, we will study the algebraic sets, the ideal of points and the relationship between them. Some properties of affine algebraic sets of commutative algebraic geometry will be extended. The results will be illustrated with examples. 
		
\section{Coding theory over single Ore extensions}\label{section2}

Following mainly \cite{Boucher}, \cite{Bueso}, \cite{Luerssen}, \cite{Gomez-Torrecillas}, \cite{Gomez-Torrecillas2}, \cite{Gomez-Torrecillas5}, \cite{Gomez-Torrecillas6}, \cite{Lam} and \cite{Leroy}, we will review in this section the most important algebraic tools and the general theory the of skew cyclic codes, i.e., the coding theory over single Ore extensions. \textbf{All that we will present in this section 
is known and have been published before in the specialized literature}. For simplicity, we will avoid to indicate the reference from which we have taken a definition or a result, only in some places we will point out the precise reference used. For completeness, and in order that the non specialized readers can follow this paper, we will include some proofs omitted in the literature. 

Let $R$ be a ring, $\sigma$ be an automorphism of $R$ and $\delta$ be a \textit{\textbf{$\boldsymbol{\sigma}$-derivation of $\boldsymbol{R}$}}, i.e., a 
function $\delta: R\to R$ such that $\delta(r+s)=\delta(r)+\delta(s)$ and 
$\delta(rs)=\sigma(r)\delta(s)+\delta(r)s$, for all $r,s\in R$. With this, we consider the \textit{\textbf{skew polynomial ring}} $\boldsymbol{R[x;\sigma,\delta]}$, also called \textit{\textbf{Ore extension}} (\cite{Ore}). Recall that $A:=R[x;\sigma,\delta]$ is a noncommutative ring of polynomial type with rule of multiplication
\[xr=\sigma(r)x+\delta(r), \ \text{for every}\ r\in R.
\]
See \cite{McConnell} (or also \cite{Lezama9}) for the algebraic and homological properties of skew polynomial rings. Let $f$ be a monic nonzero polynomial of $A$ of degree 
$n\geq 1$, $f:=x^n+f_{n-1}x^{n-1}+\cdots+f_0$, and let $Af$ be the left ideal of $A$ generated by $f$. We consider the left $A$-module $A/Af$; note that  $A/Af$ is a left $R$-module with the canonical 
product $r\cdot \overline{a}:=\overline{ra}$, for all $r\in R$ and $a\in A$. Since $f$ is monic, we have the following isomorphism of left $R$-modules:
\begin{equation}\label{equation2.1}
\mathfrak{p}_f:R^n \to A/Af,\ \ (r_0,r_1,\dots,r_{n-1}) \mapsto \overline{\sum_{i=0}^{n-1}r_ix^i}.
\end{equation}
Indeed, it is clear that $\mathfrak{p}_f$ is an $R$-homomorphism; let $(r_0,r_1,\dots,r_{n-1})\in \ker(\mathfrak{p}_f)$, then there exists $a\in A$ such that $\sum_{i=0}^{n-1}r_ix^i=af$, assume that $a\neq 0$, let $a:=a_mx^m+\cdots+a_0$, with $m\geq 0$ and $a_m\neq 0$, then
\[\sum_{i=0}^{n-1}r_ix^i=(a_mx^m+\cdots )(x^n+\cdots )=a_mx^{m+n}+\cdots,
\]
a contradiction. So, $a=0$ and hence $r_i=0$ for $0\leq i\leq n-1$. This proves that $\mathfrak{p}_f$ is a monomorphism. Now let $g\in A$, since $f$ is monic there exists $q\in A$ such that $g=qf+h$, where $h\in A$ with $h=0$ or $\deg(h)<n-1$, whence $\overline{g}=\overline{h}$. This proves that $\mathfrak{p}_f$ is an epimorphism. 

\begin{definition}\label{definition2.1}
Let $R$ be a ring and $A:=R[x;\sigma,\delta]$ be as before.
\begin{enumerate}
\item[\rm (i)]An\textbf{ $\boldsymbol{R}$-linear code} $\mathfrak{C}$ is a submodule of the left $R$-module $R^n$. The elements of $R^n$ are called \textbf{words}. 
\item[\rm (ii)]A \textbf{module $\boldsymbol{(\sigma,\delta, f)}$-skew cyclic code} $\mathcal{C}$ is a left $A$-submodule of $A/Af$. If $f=x^n-r$ or $f=x^n-1$, the codes are called \textbf{module $\boldsymbol{(\sigma,\delta, r)}$-skew constacyclic} or \textbf{module $\boldsymbol{(\sigma,\delta)}$-skew cyclic}, respectively. The polynomial $f$ is the \textbf{modulus} of the module skew cyclic code $\mathcal{C}$.
\item[\rm (iii)]If $f$ is a \textbf{two-sided polynomial}, i.e., $Af=fA$, then a module skew cyclic code is called an \textbf{ideal skew cyclic code}.   
\end{enumerate}
\end{definition}

\begin{remark}
(i) Fixing $R,\sigma, \delta $ and $f$, we will simply talk about skew cyclic codes. Observe that if $\mathcal{C}$ is a skew cyclic code, then $\mathfrak{p}_f^{-1}(\mathcal{C})$ is an $R$-linear code. In some papers a skew cyclic code is defined as 
an $R$-linear code $\mathfrak{C}$ such that $\mathfrak{p}_f(\mathfrak{C})$ is a submodule of the left $A$-module $A/Af$ (see \cite{Luerssen}). This is because
\begin{center}
\textit{there exists a bijective correspondence between}
	
\textit{the collection of $R$-linear codes such that its images by $\mathfrak{p}_f$ are submodules of the left $A$-module $A/Af$ and the collection of skew cyclic codes}.
\end{center}

(ii) Consider the particular case when $\delta=0$ and $f=x^n-1$. Then, $\mathcal{C}$ is a skew cyclic code if and only if 
\begin{center}
$(r_0,r_1,\dots,r_{n-1})\in \mathfrak{p}_f^{-1}(\mathcal{C})\Longrightarrow (\sigma(r_{n-1}),\sigma(r_0),\dots,\sigma(r_{n-2}))\in \mathfrak{p}_f^{-1}(\mathcal{C})$.
\end{center}
In fact, observe that $(\sigma(r_{n-1}),\sigma(r_0),\dots,\sigma(r_{n-2}))=\mathfrak{p}_f^{-1}(\overline{x(r_0+r_1x+\cdots+r_{n-1}x^{n-1})})$. The name given to codes as cyclic codes is justified in the particular case when $\mid\sigma\mid=n$.  
\end{remark}

\begin{definition}\label{definition2.3}
For $Z:=(z_1,\dots,z_n), Z':=(z_1',\dots,z_n')\in R^n$, the \textbf{inner product} of $Z$ and $Z'$ in $R^n$ is defined as
\begin{center}
	$Z\cdot Z':=z_1z_1'+\cdots+z_nz_n'$.
\end{center}
Let $\mathfrak{C}$ be an $R$-linear code. The \textbf{dual code} of $\mathfrak{C}$, denoted $\mathfrak{C}^{\perp}$, is the $R$-linear code defined by
\begin{center}
	$\mathfrak{C}^{\perp}:=\{Z\in R^n\mid Z\cdot Z'=0,\ \text{for all}\ Z'\in \mathfrak{C}\}$.
\end{center}
The $R$-linear code $\mathfrak{C}$ is \textbf{self-orthogonal} if $\mathfrak{C}\subseteq \mathfrak{C}^{\perp}$. $\mathfrak{C}$ is \textbf{self-dual} if $\mathfrak{C}= \mathfrak{C}^{\perp}$.
\end{definition}

\begin{definition}
Let $\mathcal{C}$ be a skew cyclic code. The \textbf{dual code} of $\mathcal{C}$, denoted $\mathcal{C}^{\perp}$, is the $R$-linear code $\mathfrak{C}^{\perp}$, where $\mathfrak{C}:=\mathfrak{p}_f^{-1}(\mathcal{C})$.
\end{definition}

\begin{remark}
Observe that if $\mathcal{C}$ be a skew cyclic code, then $\mathcal{C}^{\perp}:=\mathfrak{C}^{\perp}$ is not necessarily a skew cyclic code, i.e., $\mathfrak{p}_f(\mathcal{C}^{\perp})=\mathfrak{p}_f(\mathfrak{C}^{\perp})$ could not be a submodule of $A/Af$. In Subsection \ref{subsection2.3} we will consider a particular situation in which this $\mathcal{C}^{\perp}$ is a skew linear code.  
\end{remark}

\textbf{From now on in this paper} we will assume that $R$ is an $\textbf{\textit{IBN}}$ ring (\textbf{\textit{Invariant Basis Number}}), i.e., if $M$ is a free left $R$-module of finite bases, then all bases of $M$ have the same number of elements (see \cite{Lam2} or also \cite{Lezamadim}). With this in mind, we recall next the notion of generator matrix for $R$-linear codes. In the next subsection we will present this notion for skew cyclic codes over finite fields. 

\begin{definition}
Let $\mathfrak{C}$ be an $R$-linear code of $R^n$. Assume that $\mathfrak{C}$ is $R$- free of finite dimension $s$ and let $X:=\{v_1,\dots,v_s\}$ be a $R$-basis of $\mathfrak{C}$, with $v_i:=(v_{i1},\dots,v_{in})\in \mathfrak{C}\subseteq R^n$, $1\leq i \leq s$. The \textbf{generator matrix of $\mathfrak{C}$ with respect to $X$}, denoted $G_X(\mathfrak{C})$, is defined by
\begin{center}
$G_X(\mathfrak{C}):=\begin{bmatrix}
v_{11} & \cdots & v_{1n}\\
\vdots & \cdots & \vdots\\
v_{s1} & \cdots & v_{sn}
\end{bmatrix}\in M_{s\times n}(R)$.
\end{center} 
\end{definition} 

It is clear that the structure and properties of skew cyclic codes depend on $R,\sigma, \delta $ and $f$. The most remarkable case with practical applications is when $R=\mathbb{F}$ is a finite field and $\delta=0$ (see \cite{Boucher}, \cite{Gomez-Torrecillas3}, \cite{Luerssen}, \cite{Gomez-Torrecillas}, \cite{Gomez-Torrecillas2}, \cite{Gomez-Torrecillas4}). Another key case is when $R=\mathbb{F}(t)$, with $\mathbb{F}$ a finite field and $\delta=0$. We will study first the general theory and tools of the skew cyclic codes for the case when $R=\mathbb{F}$ is a finite field and $A:=\mathbb{F}[x;\sigma,\delta]$, and after this, we will study the mentioned remarkable cases. 

\subsection{Main tools of the skew cyclic coding theory}\label{subsection2.1}

Let $\mathbb{F}$ be a finite field of characteristic $q$, $\sigma$ be an automorphism of $\mathbb{F}$ and $\delta$ be a $\sigma$-derivation of $\mathbb{F}$. We will recall in the present subsection the general theory and the main tools involved in the study of the skew cyclic codes over the Ore extension $A:=\mathbb{F}[x;\sigma,\delta]$. We fix the modulus $f\in A$ of degree $n\geq 1$.

\subsubsection{Basic definitions and elementary facts}

(A) Recall that $|\mathbb{F}|=q^k$, where $k:=\dim_{\mathbb{Z}_q}(\mathbb{F})$ and $\mathbb{Z}_q$ is the prime subfield of $\mathbb{F}$ (see \cite{Waerden}, or also \cite{Lezama5}); up to isomorphism, $\mathbb{F}$ consists of the roots of the polynomial $z^{q^k}-z$; ${\rm Aut}(\mathbb{F})\cong \mathbb{Z}_k$ and every automorphism is a power of the \textit{\textbf{Frobenius automorphism}} defined by
\begin{equation}\label{equ2.2}
\phi: \mathbb{F}\to \mathbb{F},\ \ z\mapsto z^q, \ z\in\mathbb{F}.
\end{equation}
$\delta$ is necessarily a \textit{$\boldsymbol{\sigma}$-\textbf{inner derivation}}, i.e., there exists $w\in \mathbb{F}$ such that 
\begin{equation}\label{equ2.3}
\delta(z)=w(\sigma(z)-z),\ \text{for every} \ z\in \mathbb{F}.
\end{equation}
In fact, recall that the group $\mathbb{F}^*$ is cyclic generated by any $t$-th primitive root $z_0$ of unity, where $t:=q^k-1$. Then, for $\sigma\neq i_{\mathbb{F}}$, $w:=\frac{\delta(z_0)}{\sigma(z_0)-z_0}$ is the claimed element since given $z\in \mathbb{F}^*$ there exists $l\geq 1$ such that $z=z_0^{l}$ and (\ref{equ2.3}) can be proved by induction on $l$. If $\sigma=i_{\mathbb{F}}$, then $\delta=0$: Indeed, $\delta(z_0^q)=qz_0^{q-1}\delta(z_0)=0$ and $z_1:=z_0^q$ also generates $\mathbb{F}^*$ since ${\rm gcd}(q,q^k-1)=1$, whence for every $z\in \mathbb{F}$, $\delta(z)=\delta(z_1^l)=lz_1^{l-1}\delta(z_1)=0$.

(B) Since $A$ is a left (and right) euclidean domain, then $A$ is a left (and right) principal ideal domain, and consequently, a left (and right) noetherian domain (see \cite{McConnell}, or also \cite{Lezama9}). Hence, a skew cyclic code has the form $Ag/Af$, for some $g\in A$ with 
$Af\subseteq Ag$. Thus, there exists $h\in A$ such that $f=hg$, i.e., $g$ is a \textbf{\textit{right divisor}} of $f$, and we write $g\mid_r f$. Therefore,
\begin{center}
	\textit{there exists a bijective correspondence between the skew cyclic codes and the monic right divisors of $f$}.
\end{center}

(C) Observe that the skew cyclic code $Ag/Af$ is a $\mathbb{F}$-vector of dimension $\deg(f)-\deg(g)$. $Ag/Af$ is a subspace of the $\mathbb{F}$-vector space $A/Af$ which has dimension $n=\deg(f)$. There are two trivial skew cyclic codes: $A/Af$ of dimension $n$ and $Af/Af$ of dimension $0$. 

\begin{definition}\label{definition2.6}
Let $\mathcal{C}:=Ag/Af$ be a skew cyclic code.
\begin{enumerate}
\item[\rm (i)]The \textbf{\textit{length}} of $\mathcal{C}$ is $n$.
\item[\rm (ii)]$s:=\deg(f)-\deg(g)=\dim_{\mathbb{F}}(\mathcal{C})$ is the \textbf{\textit{dimension}} of $\mathcal{C}$.
\item[\rm (iii)]The \textbf{generator matrix} of $\mathcal{C}$, denoted $G(\mathcal{C})$, is the matrix of $M_{s\times n}(\mathbb{F})$ defined by the coefficients of $g,xg,\dots,x^{s-1}g$, disposed by rows and completing with zeros.
\end{enumerate}
\end{definition}

\begin{proposition}\label{proposition2.6}
Let $\mathcal{C}:=Ag/Af$ be a skew cyclic code of dimension $s$. Then, $\mathfrak{C}:=\mathfrak{p}_f^{-1}(\mathcal{C})$ is a $\mathbb{F}$-linear code of dimension $s$ and there exists a basis $X$ of $\mathfrak{C}$ such that $G(\mathcal{C})=G_X(\mathfrak{C})$.  
\end{proposition}
\begin{proof}
Since $\mathfrak{p}_f^{-1}$ is a $\mathbb{F}$-isomorphism, then $\mathfrak{C}$ is a $\mathbb{F}$-linear code of dimension $s$, moreover, $Y$ is a $\mathbb{F}$-basis of $\mathcal{C}$ if and only if $\mathfrak{p}_f^{-1}(Y)$ is a $\mathbb{F}$-basis of $\mathfrak{C}$. In the $\mathbb{F}$-vector space $A/Af$, observe that $Y:=\{\overline{g},\overline{xg},\dots, \overline{x^{s-1}g}\}$ is a $\mathbb{F}$-basis of $\mathcal{C}$: In fact, since $A$ is domain, $Y$ has exactly $s$ elements, moreover, if $z_1,\dots,z_{s-1}\in \mathbb{F}$ are such that $z_0\overline{g}+z_1\overline{xg}+\cdots +z_{s-1}\overline{x^{s-1}g}=\overline{0}$, then $(z_0+z_1x+\cdots+z_{s-1}x^{s-1})g=pf$, for some $p\in A$, but since $s+\deg(g)=n$, we conclude that $p=0$, so $z_0=\cdots=z_{s-1}=0$.  

Now let $X:=\mathfrak{p}_f^{-1}(Y)=\{v_0,\dots,v_{s-1}\}$, where
\begin{center}
$v_i:=\mathfrak{p}_f^{-1}(\overline{x^ig})=(v_{i1},\dots,v_{in})$ and $v_{i1},\dots,v_{in}$ are the coefficients of $x^ig$, $0\leq i \leq s-1$.
\end{center}
Now it is clear that $G(\mathcal{C})=G_X(\mathfrak{C})$. 
\end{proof}

\begin{proposition}\label{proposition2.7}
\begin{enumerate}
Let $\mathfrak{C}$ be a $\mathbb{F}$-linear code with $\dim_\mathbb{F}(\mathfrak{C})=s$. 
\item[\rm (i)]$\dim_\mathbb{F}(\mathfrak{C}^{\perp})=n-s$. Moreover, there exists a $\mathbb{F}$-basis $X^{\perp}$ of $\mathfrak{C}^{\perp}$ such that $G_{X^{\perp}}(\mathfrak{C}^{\perp})Z'^T=0$, for every $Z'\in \mathfrak{C}$. In particular,   
$G_{X^{\perp}}(\mathfrak{C}^{\perp})G_X(\mathfrak{C})^T=0$.  
\item[\rm (ii)]$(\mathfrak{C}^{\perp})^{\perp}=\mathfrak{C}$.
\item[\rm (iii)]Let $Z'\in \mathbb{F}^r$. Then, $Z'\in \mathfrak{C}$ if and only if $G_{X^{\perp}}(\mathfrak{C}^{\perp})Z'^T=0$. 
\item[]It is said that $G_{X^{\perp}}(\mathfrak{C}^{\perp})$ is a \textbf{parity check matrix} of $\mathfrak{C}$.
\end{enumerate}
\end{proposition}
\begin{proof}
(i) Since $\mathbb{F}$ is a field, it is clear that every $\mathbb{F}$-linear code has a finite $\mathbb{F}$-basis, thus, $\mathfrak{C}^{\perp}$ has a finite $\mathbb{F}$-basis $X^{\perp}\subseteq \mathbb{F}^n$; of course $X^{\perp}$ satisfies $G_{X^{\perp}}(\mathfrak{C}^{\perp})Z'^T=0$, for every $Z'\in \mathfrak{C}$. In particular,   
$G_{X^{\perp}}(\mathfrak{C}^{\perp})G_X(\mathfrak{C})^T=0$. Only remains to show that $X^{\perp}$ has $n-s$ elements. Let $l$ be the size of $X^{\perp}$, we consider the $\mathbb{F}$-linear operator
\begin{align*}
\mathbb{F}^n & \xrightarrow{T} \mathbb{F}^s\\
Z & \mapsto (Z\cdot Z_1,\dots, Z\cdot Z_s),
\end{align*} 
where $X:=\{Z_1,\dots,Z_s\}$ is a $\mathbb{F}$-basis of $\mathfrak{C}$. Observe that $\ker(T)=\mathfrak{C}^{\perp}$ and ${\rm Im}(T)$ is the $\mathbb{F}$-space generated by the columns of $G_X(\mathfrak{C})$. From this we get that $n=l+s$, so $l=n-s$. 

(ii) It is clear that $\mathfrak{C}\subseteq (\mathfrak{C}^{\perp})^{\perp}$, but from (i), $\dim_{\mathbb{F}}((\mathfrak{C}^{\perp})^{\perp})=n-(n-s)=s$, so $\mathfrak{C}=(\mathfrak{C}^{\perp})^{\perp}$.

(iii) $\Rightarrow)$: Evident.

$\Leftarrow)$: From $G_{X^{\perp}}(\mathfrak{C}^{\perp})Z'^T=0$ we get that $Z'\in (\mathfrak{C}^{\perp})^{\perp}=\mathfrak{C}$.
\end{proof}

(D) Given two polynomials $g_1,g_2\in A$, not both zero, there exists a unique monic polynomial $d\in A$ such that $d\mid_r g_1, d\mid_r g_2$ and if 
$h\in A$ is such that $h\mid_r g_1, h\mid_r g_2$, then $h\mid_r d$. In fact, $d$ is the monic generator of the left ideal $Ag_1+Ag_2$, i.e., $Ad=Ag_1+Ag_2$. The polynomial $d$ is called the \textbf{\textit{greatest common right divisor}} of $g_1$ and $g_2$, denoted by ${\rm gcrd}(g_1, g_2)$. We say that $g_1, g_2$ are \textbf{\textit{relatively right-prime}} if $d=1$. Now we assume that $g_1,g_2\in A$ are nonzero, then there exists a unique monic polynomial $l\in A$ such that $g_1\mid_r l, g_2\mid_r l$, and if $h\in A$ is such that $g_1\mid_r h, g_2\mid_r h$, then $l\mid_r h$. In fact, $l$ is the monic generator of
$Ag_1\cap Ag_2$, i.e., $Al=Ag_1\cap Ag_2$. We denote $l={\rm lclm}(g_1,g_2)$ and $l$ is called the \textbf{\textit{least common left multiple}} of $g_1$ and $g_2$. In \cite{Ore} has been proved that
\begin{equation}\label{equation2.4a}
\deg(g_1)+\deg(g_2)=\deg({\rm gcrd}(g_1, g_2))+\deg({\rm lclm}(g_1,g_2)).
\end{equation}

Similar definitions and statements are true for the left side. 

(E) The following example in \cite{Luerssen} shows that a right divisor is not necessarily a left divisor: Let $\mathbb{F}$ be the field of $4=2^2$ elements, $\mathbb{F}=\{0,1,w,w^2\}$, then $w^3=1$, so $0=(w-1)(w^2+w+1)$, i.e., $w^2=-w-1=w+1$; if $\phi$ is the Frobenius authomorphism, then for every $z\in \mathbb{F}$, $\phi^2(z)=z^4=z$, i.e, $\phi^2=i_{\mathbb{F}}$. Then, in $A:=\mathbb{F}[x;\phi]$ we have
\begin{center}
$(x^2+wx+w)(x+w)=x^3+w^2x+w^2$,
\end{center}
so $x+w$ is a right divisor of $x^3+w^2x+w^2$, but it is not a left divisor: Contrary, let $z_2x^2+z_1x+z_0\in A$ such that $(x+w)(z_2x^2+z_1x+z_0)=x^3+w^2x+w^2$, then 
\begin{center}
$x^3+w^2x+w^2=z_2^2x^3+(wz_2+z_1^2)x^2+(wz_1+z_0^2)x+wz_0$,
\end{center} 
therefore
\begin{center}
$z_2^2=1, wz_2+z_1^2=0, wz_1+z_0^2=w^2, wz_0=w^2$,
\end{center} 
so $z_0=w$ and hence $wz_1=0$, i.e., $z_1=0$, whence $z_2=0$, a contradiction.

(F) Let $g,h\in A$, we say that $g,h$ are \textit{\textbf{right associated}}, denoted $g\sim_r h$, if there exists $u\in A^*$ such that $g=uh$. It is clear that this is an equivalence relation. Observe that
\begin{center}
$g\sim_r h$ if and only if $g\mid_r h$ and $h\mid_r g$. 
\end{center}
In fact, if $g\sim_r h$, then $g=uh$, for some $u\in A^*$, so $h=u^{-1}g$, whence $g\mid_r h$ and $h\mid_r g$. Conversely, assume that $g\mid_r h$ and $h\mid_r g$, then there exist $u,u'\in A$ such that $h=ug$ and $g=u'h$. Then, $g=u'ug$ and $h=uu'h$; if $g=0$, then $h=0$ and $g\sim_r h$; assume that $g\neq 0$, then $h\neq 0$, whence $1=u'u$ and $1=uu'$, i.e., $u\in A^*$. Therefore, $g\sim_r h$.

Left associated polynomials are defined similarly. Observe that these definitions are not left-right symmetric. In the next item we will consider a close related notion which is left-right symmetric.   

(G) Let $g,h\in A$, we say that $g,h$ are left \textit{\textbf{similar}}, denoted $g\approx_l h$, if there exists an isomorphism of left $A$-modules $A/Ag\xrightarrow{\alpha} A/Ah$. 
In a similar way can be defined the right  similarity of polynomials. We will prove next that $g\approx_l h$ if and only if $g\approx_r h$. By symmetry, we will only show that if $g\approx_ l h$, then $g\approx_r h$. Let $A/Ag\xrightarrow{\alpha} A/Ah$ be an isomorphism of left $A$-modules, let $\alpha(\overline{1}):=\widetilde{p}$, with $p\in A$, then $\alpha(\overline{g})=\widetilde{0}=g\alpha(\overline{1})=g\widetilde{p}=\widetilde{gp}$, so there exists $q\in A$ such that $gp=qh$. We define $A/hA\xrightarrow{\beta} A/gA$ by $\beta(\widehat{a}):=\overline{\overline{qa}}$, for $a\in A$. $\beta$ is well-defined since if $\widehat{a}=\widehat{b}$, then $a-b=hc$, for some $c\in A$, so $qa-qb=gpc$, i.e., 
$\overline{\overline{qa}}=\overline{\overline{qb}}$. It is clear that $\beta$ is a homomorphism of right $A$-modules. Thus, we have a function $\Lambda$ defined by $\Lambda(\alpha):=\beta$. Observe that $\Lambda(\alpha_2\circ \alpha_1)=\Lambda(\alpha_2)\circ \Lambda(\alpha_1)$, with $A/Ag\xrightarrow{\alpha_1} A/Ah\xrightarrow{\alpha_2} A/Ak$: Indeed, let $\alpha_1(\overline{1}):=\widetilde{p_1}$, with $gp_1=q_1h$, and $\alpha_2(\widetilde{1}):=\widetilde{\widetilde{p_2}}$, with $hp_2=q_2k$, then $(\alpha_2\circ \alpha_1)(\overline{1})=\widetilde{\widetilde{p_1p_2}}$, with $g(p_1p_2)=(q_1q_2)k$, so $(\beta_2\circ \beta_1)(\widehat{1})=[q_2q_1]$ and $\beta_2(\beta_1(\widehat{1}))=\beta_2(\overline{\overline{q_1}})=\beta_2(\overline{\overline{1}})q_1=[q_2]q_1=[q_2q_1]$. Moreover, $\Lambda(i_{A/Ag})=i_{A/gA}$ and $\Lambda(i_{A/Ah})=i_{A/hA}$. Hence, $\beta=\Lambda(\alpha)$ is an isomorphism and so is $\beta^{-1}$, i.e., $g\approx_r h$.  

Thus, we will simply write  $g\approx h$. It is clear that $\approx$ is an equivalence relation. 

Note that if $g\sim_r h$ then $g\approx h$ (also, if $g\sim_l h$ then $g\approx h$). Indeed, there exists $u\in A^*$ such that $g=uh$, we define 
\begin{equation}\label{equation2.4}
A/Ag\xrightarrow{\alpha} A/Ah,\ \overline{a}\mapsto \widetilde{a},\ \text{for}\ a\in A.
\end{equation}
$\alpha$ is well-defined since if $\overline{a}=\overline{b}$, then $a-b\in yg$ for some $y\in A$, so $a-b=yuh$, i.e., $\tilde{a}=\tilde{b}$. It is clear that $\alpha$ is a homomorphism of left $A$-modules. If $\alpha(\overline{a})=\tilde{0}$, then $a=yh$ for some $y\in A$, so $a=yu^{-1}g$, i.e., $\overline{a}=\overline{0}$. This proves that $\alpha $ is injective. Clearly $\alpha$ is surjective. 

(H) A polynomial $p\in A$ is \textbf{\textit{irreducible}} if it satisfies the following conditions: (i) $p\neq 0$ (ii) $p\notin A^*$ (iii) if $p=hk$ for some $h,k\in A$, then $h\in A^*$ or $k\in A^*$. We know that $A$ is a left and right principal ideal domain, using this, in \cite{Bueso} it is proved that $A$ is an \textit{\textbf{unique factorization domain}} ($UFD$) in the sense of the following definition: Every element $0\neq a\notin A^*$ can be factorized in a finite product $a=p_1\cdots p_t$ of irreducible elements $p_i\in A$, $1\leq i\leq t$, and if $p_1\cdots p_t=q_1\cdots q_s$, where $p_i$, $q_j$ are irreducible, $1\leq i\leq t$, $1\leq j\leq s$, then $s=t$ and there exists a permutation $\pi$ in the symmetric group $S_t$ of all permutations of the set $\{1,\dots, t\}$ such that $q_i\approx p_{\pi(i)}$, $1\leq i\leq t$. 

Related to this we have the following basic facts (see \cite{Bueso}): 

(a) Let $p\in A$, $p$ is \textit{\textbf{right prime}} if $p\notin A^*$ and if $b,c\in A$ are such that $p\mid_r bc$, then $p\mid_r b$ or $p\mid_r c$. 

(b) Consider $\mathbb{C}[x;\sigma]$, with $\sigma(z):=\overline{z}$, $z\in \mathbb{C}$. For any $z\in \mathbb{C}$, $x-z$ is irreducible, but if $z\in \mathbb{C}-\mathbb{R}$, then $x-z$ is not prime: In fact,
\begin{center}
$(x-z)\mid_r (x+\overline{z})(x-z)=x^2-\left| z\right|^2=(x-\left| z\right|)(x+\left| z\right|) $,
\end{center}  
but $(x-z)\nmid_r (x-\left| z\right|)$ and $(x-z)\nmid_r (x+\left| z\right|)$.

(c) Since $A$ is a left principal ideal domain, then for $p\in A$, $p$ is irreducible if and only if $Ap$ is a left maximal ideal of $A$.

(d) Related to (b), since $A$ is a left and right principal ideal domain, every right prime element of $A$ is irreducible. 

\begin{remark}\label{remark2.10}
Some definitions and facts of this subsection are valid without assuming that $\mathbb{F}$ is finite, inclusive, some of them are valid for any domain, i.e., for arbitrary rings without zero divisors (see \cite{Bueso}). 
\end{remark}

\subsubsection{Roots of polynomials and algebraic sets}

Recall that we are assuming that $A:=\mathbb{F}[x;\sigma,\delta]$. Now we will consider right roots, right remainders of polynomials of $A$ and skew algebraic sets. Some related key matrices will be also studied. We will follow \cite{Boucher}, \cite{Luerssen}, \cite{Lam} and \cite{Leroy}.  

Let $g\in A$ and $z\in \mathbb{F}$, then there exist unique $q\in A$ and $r\in \mathbb{F}$ such that
\begin{center}
$g=q(x-z)+r$. 
\end{center}	
	
\begin{definition}
With the previous notation, $r$ is the \textbf{right evaluation} of $g$ in $z$, i.e., $g(z):=r$. We say that $z$ is a \textbf{right root} of $g$ if $g(z)=0$, i.e., if $r=0$.   
\end{definition}
Thus,
\begin{center}
$z$ is a right root of $g$ if and only if $(x-z)\mid_r g$.
\end{center}
The habitual replacing the variable $x$ by $z$ is not a well-defined notion for the evaluation of polynomials, for example, if $\sigma\neq i_\mathbb{F}$ and $u\in \mathbb{F}$ is such that $\sigma(u)\neq u$, then for $g:=ux=x\sigma^{-1}(u)$ in $\mathbb{F}[x;\sigma]$ and any $0\neq z\in \mathbb{F}$ we have $g(z)=uz=z\sigma^{-1}(u)=\sigma^{-1}(u)z$, i.e., $u=\sigma^{-1}(u)$, a contradiction.
 
Example 2.7 in \cite{Luerssen} shows that a polynomial of degree $n$ may have more that $n$ right roots: We consider again $\mathbb{F}[x;\phi]$ as in part (E) of the previous subsection, then
\begin{center}
$x^2+1=(x + 1)(x + 1)=(x + w^2)(x + w)=(x + w)(x + w^2)$,
\end{center} 
i.e., $-1, -w, -w^2$ are right roots of $x^2+1$. On the other hand, if the field of coefficients is infinite, then a polynomial may have infinitely many right roots. In fact, consider $\mathbb{C}[x;\sigma]$, with $\sigma(z):=\overline{z}$, $z\in \mathbb{C}$, then $x^2-1=(x+\overline{z})(x-z)$, for any $z$ on
the unit circle. 

\begin{definition}
	For $i\geq 0$, the $i$-th \textbf{norm} on $\mathbb{F}$ is defined by
	\begin{center}
		$N_0(z):=1$, $N_{i+1}(z):=\sigma(N_i(z))z+ \delta(N_i(z))$, for any  $z\in \mathbb{F}$.
	\end{center}
\end{definition}

The right evaluation of a polynomial can be computed using the norm.

\begin{proposition}[\cite{Leroy}, Lemma 2.4]\label{proposition2.13}
Let $g=g_0+\cdots +g_mx^m\in A$ and $z\in \mathbb{F}$. Then, 
\begin{center}
$g(z)=\sum_{i=0}^m g_iN_i(z)$.
\end{center} 
\end{proposition}
\begin{proof}
Observe that
\begin{center}
$g-\sum_{i=0}^{m}g_iN_i(z)=\sum_{i=0}^{m}g_ix^i-\sum_{i=0}^{m}g_iN_i(z)=\sum_{i=0}^{m}g_i(x^i-N_i(z))$.
\end{center} 
If we prove that for every $i\geq 0$, $x^i-N_i(z)\in A(x-z)$, then there exists $q\in A$ such that $g=q(x-z)+\sum_{i=0}^{m}g_iN_i(z)$, so by the uniqueness of the right evaluation we conclude that $g(z)=\sum_{i=0}^m g_iN_i(z)$.  

We prove the claimed by induction. For $i=0$ we have $x^0-N_0(z)=0\in A(x-z)$. For $i=1$, $x-N_1(z)=x-z\in A(x-z)$. Assume that the assertion is true for $i$, then
\begin{center}
$x^{i+1}-N_{i+1}(z)=x^{i+1}-\sigma(N_i(z))z-\delta(N_i(z))=x^{i+1}+\sigma(N_i(z))(x-z)-\sigma(N_i(z))x-\delta(N_i(z))=\sigma(N_i(z))(x-z)+x(x^i-N_i(z))\in A(x-z)$.
\end{center} 
\end{proof}

\begin{proposition}\label{proposition2.14A}
Let $\mathcal{D}$ be $\sigma$ or $\delta$ and  
\begin{center}
$\mathbb{F}[\mathcal{D};\circ]:=\{\sum_{i=0}^m g_i\mathcal{D}^i\mid g_i\in \mathbb{F}, m\geq 0\}$
\end{center}
where $g_i\mathcal{D}^i:\mathbb{F}\to \mathbb{F}$ is defined by $(g_i\mathcal{D}^i)(z):=g_i\mathcal{D}^i(z)$, $z\in \mathbb{F}$. Then,
\begin{enumerate}
	\item[\rm (i)]$\mathbb{F}[\mathcal{D};\circ]$ is a ring, where the addition is the habitual addition by coefficients and the multiplication is the composition of functions.
	\item[\rm (ii)]The function
	\begin{align*}
		\mathcal{L}: A & \to \mathbb{F}[\mathcal{D};\circ]\\
		g:=\sum_{i=0}^m g_ix^i & \mapsto \mathcal{L}(g):=\mathcal{L}_g:=\sum_{i=0}^m g_i\mathcal{D}^i
	\end{align*}
	is a surjective ring homomorphism, where $A=\mathbb{F}[x;\sigma]$ if $\mathcal{D}=\sigma$ and $\delta=0$, and $A=\mathbb{F}[x;\sigma, \delta]$ if $\mathcal{D}=\delta\neq 0$.  
\end{enumerate} 

\end{proposition}
\begin{proof}
(i) Evident.

(ii) All conditions of a ring homomorphism can be easy checked except, probably, the product. We have to show that
\begin{center}
	$\mathcal{L}(g_ix^ih_jx^j)=\mathcal{L}(g_ix^i)\circ \mathcal{L}(h_jx^j)$, for all $i,j\geq 0$ and $g_i,h_j\in \mathbb{F}$.
\end{center} 
The proof of this identity is by induction on $i$. We consider first that $\mathcal{D}=\sigma$. For $i=0$ the assertion is trivial; let $i=1$, then 
\begin{center}
$g_ixh_jx^j=g_i\sigma(h_j)x^{j+1}\mapsto g_i\sigma(h_j)\mathcal{D}^{j+1}=g_i\sigma(h_j)\sigma^{j+1}$, so for $z\in \mathbb{F}$, $(g_i\sigma(h_j)\sigma^{j+1})(z)=g_i\sigma(h_j)\sigma^{j+1}(z)$; 

$\mathcal{L}(g_ix)\circ \mathcal{L}(h_jx^j)=(g_i\mathcal{D})\circ (h_j\mathcal{D}^j)$, so for $z\in \mathbb{F}$, $(g_i\mathcal{D})\circ (h_j\mathcal{D}^j)(z)=
(g_i\sigma)\circ (h_j\sigma^j)(z)=g_i\sigma(h_j)\sigma^{j+1}(z)$. 
\end{center} 
We assume that the identity is true for $i$, from this we get
\begin{center}
$\mathcal{L}(g_ix^{i+1}h_jx^j)=\mathcal{L}(g_ix^ixh_jx^j)=\mathcal{L}(g_ix^i\sigma(h_j)x^{j+1})=\mathcal{L}(g_ix^i)\circ \mathcal{L}(\sigma(h_j)x^{j+1})=(g_i\mathcal{D}^{i})\circ (\sigma(h_j)\mathcal{D}^{j+1})$, so for $z\in \mathbb{F}$, $(g_i\mathcal{D}^{i})\circ (\sigma(h_j)\mathcal{D}^{j+1})(z)=g_i\sigma^{i+1}(h_j)\sigma^{i+j+1}(z)$;

$\mathcal{L}(g_ix^{i+1})\circ \mathcal{L}(h_jx^j)=(g_i\mathcal{D}^{i+1})\circ (h_j\mathcal{D}^j)$, so for $z\in \mathbb{F}$, $(g_i\mathcal{D}^{i+1})\circ (h_j\mathcal{D}^j)(z)=g_i\sigma^{i+1}(h_j)\sigma^{i+j+1}(z)$.
\end{center}

Now we assume that $\mathcal{D}=\delta$. For $i=0$ the statement is trivial. Let $i=1$, then 
\begin{center}
	$g_ixh_jx^j=g_i\sigma(h_j)x^{j+1}+g_i\delta(h_j)x^j\mapsto g_i\sigma(h_j)\delta^{j+1}+g_i\delta(h_j)\delta^j$, so for $z\in \mathbb{F}$, $(g_i\sigma(h_j)\delta^{j+1}+g_i\delta(h_j)\delta^j)(z)=g_i\sigma(h_j)\delta^{j+1}(z)+g_i\delta(h_j)\delta^j(z)$; 
	
	$\mathcal{L}(g_ix)\circ \mathcal{L}(h_jx^j)=(g_i\delta)\circ (h_j\delta^j)$, so for $z\in \mathbb{F}$, $(g_i\delta)\circ (h_j\delta^j)(z)=g_i\delta(h_j\delta^j(z))
	=g_i\sigma(h_j)\delta^{j+1}(z)+g_i\delta(h_j)\delta^j(z)$. 
\end{center} 
We assume that the identity is true for $i$ and from this we get
\begin{center}
	$\mathcal{L}(g_ix^{i+1}h_jx^j)=\mathcal{L}(g_ix^ixh_jx^j)=\mathcal{L}(g_ix^i\sigma(h_j)x^{j+1}+g_ix^i\delta(h_j)x^j)=\mathcal{L}(g_ix^i)\circ \mathcal{L}(\sigma(h_j)x^{j+1})+\mathcal{L}(g_ix^i)\circ\mathcal{L}(\delta(h_j)x^j)=(g_i\delta^{i})\circ (\sigma(h_j)\delta^{j+1})+(g_i\delta^{i})\circ (\delta(h_j)\delta^{j})$, so for $z\in \mathbb{F}$, $(g_i\delta^{i})\circ (\sigma(h_j)\delta^{j+1})(z)+(g_i\delta^{i})\circ (\delta(h_j)\delta^{j})(z)=
	(g_i\delta^{i})(\sigma(h_j)\delta^{j+1}(z))+(g_i\delta^{i})(\delta(h_j)\delta^{j}(z))$;
	
	$\mathcal{L}(g_ix^{i+1})\circ \mathcal{L}(h_jx^j)=(g_i\delta^{i+1})\circ (h_j\delta^j)=g_i\delta^i\circ(\delta\circ h_j\delta^j)$, so for $z\in \mathbb{F}$, $(g_i\delta^i\circ(\delta\circ h_j\delta^j))(z)=g_i\delta^i[\delta(h_j\delta^j(z))]=g_i\delta^i[\sigma(h_j)\delta^{j+1}(z)+\delta(h_j)\delta^j(z)]=(g_i\delta^{i})(\sigma(h_j)\delta^{j+1}(z))+(g_i\delta^{i})(\delta(h_j)\delta^{j}(z))$.
\end{center}
It is clear that $\mathcal{L}$ is surjective. This completes de proof.
\end{proof}

\begin{remark}
(i) As we observed in (\ref{equ2.3}), $\delta$ is neccesarily a $\sigma$-inner derivation, and it is well-known that $\mathbb{F}[x;\sigma, \delta]\cong \mathbb{F}[y; \sigma]$, with $y:=x-w$ (see \cite{McConnell} or also \cite{Lezama9}). Thus, for $\delta\neq 0$ in the previous proposition, 
we have also a surjective ring homomorphism from $\mathbb{F}[y; \sigma]$ into $\mathbb{F}[\mathcal{D}; \circ]$.

(ii) The \textit{\textbf{operator evaluation}} of $g$ at $z\in \mathbb{F}$ is $ \mathcal{L}_g(z)=\sum_{i=0}^m g_i\mathcal{D}^i(z)$.
\end{remark}

We return again to the general case  $A:=\mathbb{F}[x;\sigma,\delta]$.

\begin{definition}\label{definition2.16}
Let $Z:=(z_1,\dots,z_r)\in \mathbb{F}^r$. Then,
\begin{enumerate}
\item[\rm (i)]The \textbf{Vandermonde matrix} of $Z$ is defined by
\begin{center}
${\rm V}_r(Z):=\begin{bmatrix}1 & 1 & \cdots & 1\\
N_1(z_1) & N_1(z_2) & \cdots & N_1(z_r)\\
\vdots & \vdots & \cdots & \vdots\\
N_{r-1}(z_1) & N_{r-1}(z_2) & \cdots & N_{r-1}(z_r)\end{bmatrix}$.
\end{center}
\item[\rm (ii)]The \textbf{Wronskian matrix} of $Z$ is defined by
\begin{center}
	${\rm Wr}_r(Z):=\begin{bmatrix}z_1 & z_2 & \cdots & z_r\\
	\mathcal{D}(z_1) & \mathcal{D}(z_2) & \cdots & \mathcal{D}(z_r)\\
	\vdots & \vdots & \cdots & \vdots\\
	\mathcal{D}^{r-1}(z_1) & \mathcal{D}^{r-1}(z_2) & \cdots & \mathcal{D}^{r-1}(z_r)\end{bmatrix}$,
\end{center}
where $\mathcal{D}=\sigma$ if $\delta=0$ and $\mathcal{D}=\delta$ if $\delta\neq 0$.
\item[\rm (iii)]For $z\in \mathbb{F}$ and $u\in \mathbb{F}^*$, let $z^u:=\sigma(u)zu^{-1}+\delta(u)u^{-1}$. We say that $z,z'\in \mathbb{F}$ are \textbf{conjugate}, 
denoted $z\equiv z'$, if there exists $u\in \mathbb{F}^*$ such that $z'=z^u$. The \textbf{conjugacy class} of $z$ is
\begin{center}
$\Delta(z):=\{z^u\mid u\in \mathbb{F}^*\}$. 
\end{center}
\item[\rm (iv)]For $z\in \mathbb{F}$,
\begin{center}
$\mathcal{C}(z):=\{u\in \mathbb{F}^*\mid z^u=z\}\cup \{0\}$.
\end{center}
\end{enumerate}
\end{definition}

\begin{proposition}\label{proposition2.9}
$\equiv$ is an equivalence relation in $\mathbb{F}$.
\end{proposition}
\begin{proof}
$z\equiv z$ taking $u=1$; if $z\equiv z'$, then $z'\equiv z$ since $z=(z')^{u^{-1}}$. In fact, 
\begin{center}
$\sigma(u^{-1})z'u+\delta(u^{-1})u=\sigma(u^{-1})[\sigma(u)zu^{-1}+\delta(u)u^{-1}]u+\delta(u^{-1})u=z+\sigma(u^{-1})\delta(u)+\delta(u^{-1})u=z+\delta(u^{-1}u)=z+\delta(1)=z$.
\end{center}
If $z\equiv z'$ with $z'=z^u$ and $z'\equiv z''$ with $z''=(z')^v$, then $z\equiv z''$ since $z''=z^{vu}$. Indeed,
\begin{center}
$z'=\sigma(u)zu^{-1}+\delta(u)u^{-1}$ and $z''=\sigma(v)z'v^{-1}+\delta(v)v^{-1}$,
\end{center}
hence
\begin{center}
$z''=\sigma(v)[\sigma(u)zu^{-1}+\delta(u)u^{-1}]v^{-1}+\delta(v)v^{-1}=\sigma(vu)zu^{-1}v^{-1}+\sigma(v)\delta(u)u^{-1}v^{-1}+\delta(v)v^{-1}
=\sigma(vu)zu^{-1}v^{-1}+\delta(vu)u^{-1}v^{-1}=z^{vu}$.
\end{center}
\end{proof}

\begin{proposition}\label{proposition2.10}
Let $w\in \mathbb{F}$ as in $($\ref{equ2.3}$)$. Then, 
\begin{enumerate}
\item[\rm (i)]For every $z\in \mathbb{F}$ and $u\in \mathbb{F}^*$, 
\begin{center}
$z^u=\sigma(u)(z+w)u^{-1}-w$.
\end{center}
\item[\rm (ii)]$\Delta(-w)=\{-w\}$. For $z\neq -w$, $|\Delta(z)|=\frac{q^k-1}{q^r-1}$, with $r:={\rm gcd}(l,k)$, where $l$ is such that 
$\sigma=\phi^l$. Thus, there are $q^r$ conjugacy classes: The conjugacy class of $-w$ and $q^r-1$ classes with $\frac{q^k-1}{q^r-1}$
elements in each class. In particular, if $\sigma=\phi$, then $r=l=1$ and there are $q$ conjugacy classes.
\item[\rm (iii)]If $\delta=0$, $\sigma=\phi$ and $q=2$, then there are two conjugacy classes: The class of $1$ with $2^k-1$ elements and the class of $0$ which is $\{0\}$. 
\item[\rm (iv)]For $z\in \mathbb{F}$, $\mathcal{C}(z)$ is a subfield of $\mathbb{F}$. If $\delta=0$, then $\mathcal{C}(1)=\{u\in \mathbb{F}\mid \sigma(u)=u\}$. If $\delta\neq 0$, then $\mathcal{C}(0)=\{u\in \mathbb{F}\mid \delta(u)=0\}$.
\item[\rm (v)]For $z,z'\in \mathbb{F}$, $z\equiv z'$, if and only if, there exists $u\in \mathbb{F}^*$ such that $(x-z')u=\sigma(u)(x-z)$, if and only if, $x-z\approx x-z'$.
\item[\rm (vi)]For $z\neq z'\in \mathbb{F}$, ${\rm lclm}(x-z,x-z')=(x-z'^{z'-z})(x-z)=(x-z^{z-z'})(x-z')$.
\item[\rm (vii)]If $g = (x-z_1)\cdots (x-z_m)$, with
$z_i\in \mathbb{F}$, $1\leq i\leq m$, and $g(z)=0$ for some $z\in \mathbb{F}$, then $z\equiv z_i$ for some $1\leq i\leq m$.
\end{enumerate}
\end{proposition}
\begin{proof}
(i) From  (\ref{equ2.3}) we get
\begin{center}
$z^u=\sigma(u)zu^{-1}+\delta(u)u^{-1}=\sigma(u)zu^{-1}+w(\sigma(u)-u)u^{-1}=\sigma(u)(z+w)u^{-1}-w$.
\end{center}

(ii) By (i), for every $u\in \mathbb{F}^*$, $(-w)^u=\sigma(u)(-w+w)u^{-1}-w=-w$, so $\Delta(-w)=\{-w\}$. Let $z\neq -w$, by (i), there exists a bijective correspondence between $\Delta(z)$ and $\{\sigma(u)u^{-1}\mid u\in \mathbb{F}^*\}$, but the cardinality of this last one set is $\frac{q^k-1}{q^r-1}$, with $r$ as in the statement of the proposition. In fact, in $\mathbb{F}^*$ we define the relation $u\asymp v$ if and only if $\sigma(u)u^{-1}=\sigma(v)v^{-1}$. It is clear that $\asymp$ is an equivalence relation, then the cardinality of $\{\sigma(u)u^{-1}\mid u\in \mathbb{F}^*\}$ is the number of different classes. To proof the claimed we want to show first that all classes have the same cardinality, namely, $q^r-1$. Let $[u]$ be the class of $u$, then $[u]=\{v\in \mathbb{F}^*\mid \sigma(u)u^{-1}=\sigma(v)v^{-1}\}=\{v\in \mathbb{F}^*\mid \sigma(uv^{-1})=uv^{-1}\}$. Let
\begin{center}
$(\mathbb{F}^*)^{\sigma}:=\{a\in \mathbb{F}^*\mid \sigma(a)=a\}$;
\end{center}
notice that $h:(\mathbb{F}^*)^{\sigma}\to [u]$, $a\mapsto ua^{-1}$, is a bijective function: Firstly, $h(a)\in [u]$ since $\sigma(ua^{-1})(ua^{-1})^{-1}=\sigma(u)\sigma(a)^{-1}au^{-1}=\sigma(u)a^{-1}au^{-1}=\sigma(u)u^{-1}$; now, if $ua_1^{-1}=ua_2^{-1}$, then $a_1=a_2$, i.e., $h$ is injective; finally, if $v\in [u]$, then we define $a:=uv^{-1}$, and hence, $\sigma(a)=\sigma(uv^{-1})=uv^{-1}=a$, i.e., $a$ is a fixed point of $h$ and $h(a)=ua^{-1}=u(uv^{-1})^{-1}=v$, i.e., $h$ is surjective. We will show that $|(\mathbb{F}^*)^{\sigma}|=q^r-1$, i.e., $[u]=q^r-1$. Since $\sigma=\phi^l$, for some $0\leq l\leq k-1$, we get that $a\in (\mathbb{F}^*)^{\sigma}$ if and only if $a\in \mathbb{F}^*$ and $a^{q^l}=a$. Let $f_l(x):=x^{q^l}-x\in \mathbb{Z}_q[x]$, then all roots of $f_l(x)$ are simple and $F_{q^l}:=\{x\in \mathbb{E}\mid x^{q^l}=x\}$ is a field of size $q^l$, where $\mathbb{E}$ is the field of decomposition of $f_l(x)$. Let $F:=(\mathbb{F}^*)^{\sigma}\cup \{0\}=\mathbb{F}\cap F_{q^l}$, then $|(\mathbb{F}^*)^{\sigma}|=|F|-1=q^r-1$, where $r:={\rm gcd}(l,k)$. Finally, let $t:=|\mathbb{F}^*/\asymp|$, so $|\mathbb{F}^*|=|[u_1]|+\cdots+|[u_t]|=t(q^r-1)$, i.e., $t=\frac{q^k-1}{q^r-1}$, but $t=|\Delta(z)|$. This completes the proof of claimed. The last assertion of (ii) is trivial.

(iii) This follows from (ii) since in this situation $w=0$. Observe that $\Delta(1)=\{\sigma(u)u^{-1}\mid u\in \mathbb{F}^*\}$. 

(iv) Let $z\in \mathbb{F}$ and $0\neq u,v\in \mathcal{C}(z)$. Then, as was observed in the proof of Proposition \ref{proposition2.9}, $(z^{u})^v=z^{uv}$, but $(z^{u})^v=z$ since
$(z^{u})^v=\sigma(v)z^uv^{-1}+\delta(v)v^{-1}=\sigma(v)zv^{-1}+\delta(v)v^{-1}=z^v=z$. Thus, $uv\in \mathcal{C}(z)$. Now, $u^{-1}\in \mathcal{C}(z)$ since from $\sigma(u)zu^{-1}+\delta(u)u^{-1}=z$ we get $z=\sigma(u)^{-1}zu-\sigma(u)^{-1}\delta(u)$, but $-\sigma(u)^{-1}\delta(u)=\delta(u^{-1})u$, so $z=\sigma(u^{-1})z(u^{-1})^{-1}+\delta(u^{-1})(u^{-1})^{-1}$, i.e., $u^{-1}\in \mathcal{C}(z)$. Finally, if $u-v=0$, then $u-v\in \mathcal{C}(z)$; let $u-v\neq 0$, then let $z':=\sigma(u-v)z(u-v)^{-1}+\delta(u-v)(u-v)^{-1}$, so from this we obtain $\sigma(u-v)z+\delta(u-v)=z'(u-v)$, whence $\sigma(u)z+\delta(u)-\sigma(v)z-\delta(v)=z'(u-v)=zu-zv=z(u-v)$, so $z=z'$. Thus, $u-v\in \mathcal{C}(z)$. This completes the proof that $\mathcal{C}(z)$ is a subfield of $\mathbb{F}$. 

Let $\delta=0$ and $u\neq 0$. Then, $u\in \mathcal{C}(1)$ if and only if $\sigma(u)u^{-1}+\delta(u)u^{-1}=1$, i.e., $\sigma(u)=u$.

Let $\delta\neq 0$ and $u\neq 0$. Then, $u\in \mathcal{C}(1)$ if and only if $\delta(u)u^{-1}=0$, i.e., $\delta(u)=0$.

(v) $z\equiv z'$ if and only if there exists $u\in \mathbb{F}^*$ such that $z'=\sigma(u)zu^{-1}+\delta(u)u^{-1}$ if and only if $(x-z')u=(x-(\sigma(u)zu^{-1}+\delta(u)u^{-1}))u=xu-\sigma(u)z-\delta(u)=\sigma(u)x+\delta(u)-\sigma(u)z-\delta(u)=\sigma(u)(x-z)$. 

If $(x-z')u=\sigma(u)(x-z)$, then $(x-z')u\sim_r (x-z)$, whence $(x-z')u\approx (x-z)$, so there exists an isomorphism of left $A$-modules $A/A(x-z')u\rightarrow A/A(x-z)$, and from this we know that there exists an isomorphism of right $A$-modules $A/(x-z')uA\rightarrow A/(x-z)A$, but $A/(x-z')A=A/(x-z')uA$, so we have an isomorphism of right $A$-modules $A/(x-z')A\rightarrow A/(x-z)A$, i.e., $(x-z)\approx (x-z')$. Conversely, assume that  $(x-z)\approx (x-z')$, then we have an isomorphism of left $A$-modules $A/A(x-z')\xrightarrow{\alpha} A/A(x-z)$, and from this $\alpha(\overline{0})=\widetilde{0}=\alpha(\overline{(x-z')})=(x-z')\alpha(\overline{1})=(x-z')\widetilde{g}$, for some $g\in A$, with $\widetilde{g}:=\alpha(\overline{1})\neq \widetilde{0}$ since $\alpha$ is injective. There exist unique $q,r\in A$, with $r\in \mathbb{F}^*$, such that $g=q(x-z)+r$, hence, $\widetilde{g}=\widetilde{r}$. Let $u:=r$, from $\widetilde{0}=(x-z')\widetilde{g}$ we get that $(x-z')u=v(x-z)$, for some $v\in A$. Taking the degree at both sides we get that $v\in \mathbb{F}^*$, moreover, $v=\sigma(u)$. Thus,  $(x-z')u=\sigma(u)(x-z)$, with $u\in \mathbb{F}^*$.    

(vi) We prove first that $(x-z'^{z'-z})(x-z)=(x-z^{z-z'})(x-z')$: Observe that $z'^{z'-z}=az'+b$ and $z^{z-z'}=az+b$, with $a:=\sigma(z-z')(z-z')^{-1}$ and 
$b=\delta(z-z')(z-z')^{-1}$, then $(x-z'^{z'-z})(x-z)=(x-az'-b)(x-z)$ and $(x-z^{z-z'})(x-z')=(x-az-b)(x-z')$, i.e., 
\begin{center}
$(x-z'^{z'-z})(x-z)=x^2-xz-(az'+b)x+(az'+b)z$, $(x-z^{z-z'})(x-z')=x^2-xz'-(az+b)x+(az+b)z'$,
\end{center}
then a direct computation shows that $(x-z'^{z'-z})(x-z)-(x-z^{z-z'})(x-z')=0$.

Let $l:=(x-z'^{z'-z})(x-z)=(x-z^{z-z'})(x-z')$, then $(x-z)\mid_r l$ and $(x-z')\mid_r l$. Now let $h\in A$ such that $(x-z)\mid_r h$ and $(x-z')\mid_r h$, then there are $f_1,f_2\in A$ such that $h=f_1(x-z)=f_2(x-z')$. Set $q_1:=(x-z')(z'-z)^{-1}$ and $q_2:=1+(x-z')(z'-z)^{-1} $. It is not difficult to show that $(f_2-f_1)q_1=f_1$ and $(f_2-f_1)q_2=f_2$, thus $q_1\mid_r f_1$ and $q_2\mid_r f_2$. However $q_1=\sigma((z'-z)^{-1})(x-z'^{z'-z})$ and $q_2=-\sigma((z-z')^{-1})(x-z^{z-z'})$. From this, it follows that $l\mid_r h$.

(vii) From $g(z)=0$ we get $g=q(x-z)$ for some $q\in A$, but since every $(x-z_i)$ is irreducible and $A$ is $UFD$, then there exists $i$ such that $x-z\approx x-z_i$, so from (v), $z\equiv z_i$. 
\end{proof}

Now we pass to study the algebraic sets in the context of the skew polynomial ring $A:=\mathbb{F}[x;\sigma,\delta]$. 

\begin{definition}\label{definition2.11}
Let $g\in A$. The \textbf{right vanishing set} of $g$, also called the \textbf{set of right roots} of $g$, is denoted by $V(g)$, i.e.,  
\begin{equation}
V(g):=\{z\in \mathbb{F}\mid g(z)=0\}.
\end{equation}
A subset $X\subseteq \mathbb{F}$ is \textbf{algebraic} if there exists $0\neq g\in A$ such that $X\subseteq V(g)$. The monic polynomial $0\neq m_X$ 
of smallest degree such that $X\subseteq V(m_X)$ is called the \textbf{minimal polynomial} of $X$. The \textbf{rank} of $X$ is defined by $\boldsymbol{{\rm rank}(X)}:=\deg(m_X)$. A monic polynomial $g\in A$ is a \textbf{$\boldsymbol{W}$-polynomial} $($\textbf{Wedderburn polynomial}\,$)$ if $g=m_X$ for some $X\subseteq \mathbb{F}$.
\end{definition}

Next we will show some results about the sets of right roots, algebraic sets and minimal polynomials. We start with the following elementary facts, some of them well-known, and others, probably new: The Zariski topology for $\mathbb{F}$ and the left ideal of points. These results will generalized in the last section to skew $PBW$ extensions.

\begin{proposition}\label{prop2.12A}
\begin{enumerate}
\item[\rm (i)]Let $X\subseteq \mathbb{F}$ be algebraic. Then, $m_X$ is well-defined. 
\item[\rm (ii)]For any $g,h,t\in A$, $V(h)\subseteq V(g)\Rightarrow V(ht)\subseteq V(gt)$.
\item[\rm (iii)]Let $I$ be a left ideal of $A$ and
\begin{center}
$V(I):=\{z\in \mathbb{F}\mid h(z)=0,\ \text{for every $h\in I$}\}$.
\end{center}
Then, $V(I)=V(g)$, where $I=Ag$.
\item[\rm (iv)]
\begin{enumerate}
	\item[\rm (a)]$V(0)=\mathbb{F}=V(l)$, with $l:={\rm lclm}(x-z\mid z\in \mathbb{F})$. Thus, $\mathbb{F}$ is algebraic.
\item[\rm (b)]$V(A)=V(1)=\emptyset$. Thus, $\emptyset$ is algebraic.
\item[\rm (c)]$V(I)\cup V(J)\subseteq V(I\cap J)$, where $I,J$ are left ideals of $A$.
\item[\rm (d)]$V(\sum_{k\in \mathcal{K}}I_k)=\bigcap_{k\in \mathcal{K}}V(I_k)$. 
\item[\rm (e)]$\mathbb{F}$ has a \textbf{Zariski topology} where the closed sets are the algebraic sets.
\end{enumerate} 
\item[\rm (v)]Let $X\subseteq \mathbb{F}$. Then,
\begin{center}
$I(X):=\{g\in A\mid g(z)=0 \ \text{for every $z\in X$}\}$
\end{center}
is a left ideal of $A$, called the \textbf{left ideal of points} of $X$. Some properties of the left ideal of points are:
\begin{enumerate}
	\item[\rm (a)]$I(\emptyset)=A$.
	\item[\rm (b)]Let $X\neq \emptyset $. Then, $I(X)=Al$, where $l={\rm lclm}(x-z\mid z\in X)$. 
	\item[\rm{(c)}]If $X,Y\subseteq \mathbb{F}$, $X\subseteq Y\Rightarrow I(Y)\subseteq I(X)$.
	\item[\rm{(d)}]If $g\in A$, then $Ag\subseteq I(V(g))$. Thus, if $I$ is a left ideal of $A$, then $I\subseteq I(V(I))$. 
	\item[\rm{(e)}]$X\subseteq V(I(X))$. Thus, every subset $X$ is algebraic and the above topology is the discrete topology.
	\item[\rm{(f)}]If $g\in A$, $V(I(V(g)))=V(g)$. Thus, if $X= V(g)$, then 
	$V(I(X))=X$.
	\item[\rm{(g)}] $I\left(  V\left(  I\left( X\right)  \right)  \right)  =I\left(
	X\right)$. 
	\item[\rm{(h)}]$I(\bigcup_{k\in \mathcal{K}}X_k)=\bigcap_{k\in \mathcal{K}}I(X_k)$.
	\end{enumerate}
\end{enumerate}
\end{proposition}
\begin{proof}
(i) If $X=\emptyset$ then $\emptyset=V(1)$, i.e., $\emptyset$ is algebraic and $m_\emptyset=1$. 

Assume that $X\neq \emptyset$  and let $0\neq g,h\in A$ be monic polynomials of smallest degree such that $X\subseteq V(g)$ and $X\subseteq V(h)$, we have to show that $g=h$. For this we observe first that if $p,q\in A$ and $z\in \mathbb{F}$, then
\begin{center}
$(p+q)(z)=p(z)+q(z)$.
\end{center}
In fact, there exist unique $q_1,q_2\in A$ and $r_1,r_2\in \mathbb{F}$ such that $p=q_1(x-z)+r_1$ and $q=q_2(x-z)+r_2$, so $p+q=(q_1+q_2)(x-z)+ (r_1+r_2)$, so the claimed follows by unicity. 

For $g,h$ we have $g=qh+r$ for some unique $q,r\in A$ with $\deg(r)<\deg(h)$ or $r=0$. Let $z\in X$, since  $h=c(x-z)$ for some $c\in A$, then $qh=qc(x-z)$, i.e., $(qh)(z)=0$. 
From this, $g(z)=0=(qh+r)(z)=qh(z)+r(z)=r(z)$, so necessarily $r=0$ because of the condition on the degree of $h$, whence, $g=qh$, but from the degree condition we get that $q\in \mathbb{F}^*$, but since $g$ and $h$ are monic, then $g=h$.

(ii) There exist unique $c,r\in A$ such that $g=ch+r$; let $z\in V(h)$, then $h=q(x-z)$ for some $q\in A$, whence $g=cq(x-z)+r$, but $z\in V(g)$, so $r=0$. Thus, $g=ch$, and hence, $gt=cht$. From this, if $z\in V(ht)$, then $ht=a(x-z)$, for some $a\in A$, whence $gt=ca(x-z)$, i.e., $z\in V(gt)$.   

(iii) Since in $A$ every left ideal is principal, then there exists $g\in A$ such that $I=Ag$. From this it is clear that $V(I)\subseteq V(g)$. Let $qg\in I$ and $z\in V(g)$, then $g(z)=0$, and as we saw in (i), this implies that $(qg)(z)=0$, i.e., $z\in V(I)$. This proves that $V(g)\subseteq V(I)$.

(iv) (a) It is clear that $V(0)=\mathbb{F}$. For the second equality, $x-z\mid_r l$, for every $z\in \mathbb{F}$, so $\mathbb{F}\subseteq V(l)$, i.e., $\mathbb{F}=V(l)$.

 (b) Evident. 
 
 (c) Since $I\cap J\subseteq I,J$, then $V(I)\cup V(J)\subseteq V(I\cap J)$.

(d) Since $I_k\subseteq \sum_{k\in \mathcal{K}}I_k$ for every $k\in \mathcal{K}$, then $V(\sum_{k\in \mathcal{K}}I_k)\subseteq \bigcap_{k\in \mathcal{K}}V(I_k)$. Let $z\in \bigcap_{k\in \mathcal{K}}V(I_k)$ and let $g\in \sum_{k\in \mathcal{K}}I_k$, then $g=g_{k_1}+\cdots+g_{k_t}$, with $g_{k_j}\in I_{k_j}$, $1\leq j\leq t$, then, as we saw in (i), $g(z)=g_{k_1}(z)+\cdots+g_{k_t}(z)=0$, whence $z\in V(\sum_{k\in \mathcal{K}}I_k)$. Thus, $\bigcap_{k\in \mathcal{K}}V(I_k)\subseteq V(\sum_{k\in \mathcal{K}}I_k)$.

(e) This follows from (a)-(d).

(v) From the proof of (i) we get that $I(X)$ is a left ideal of $A$.

(a) This is trivial.

(b) For every $z\in X$, $x-z\mid_r l$, so $l(z)=0$, whence $(al)(z)=0$ for every $a\in A$. Hence, $Al\subseteq I(X)$. Let $h\in A$ such that $I(X)=Ah$, then $h(z)=0$ for every $z\in X$, so $x-z\mid_r h$ for every $z\in X$, hence $l\mid_r h$, whence $h=pl$ for some $p\in A$. This means that $Ah\subseteq Al$. Thus, $Ah=Al$.

(c) This is trivial. 

(d) Since $g(z)=0$ for every $z\in V(g)$, then $g\in I(V(g))$, whence $Ag\subseteq I(V(g))$. The second assertion in (d) follows from (iii). 

(e) For $X=\emptyset$ the assertion follows from (a) and (iv)-(b). Let $X\neq \emptyset$. If $z\in X$, then for every $g\in I(X)$, $g(z)=0$, and this means that $z\in V(I(X))$. Therefore, $X\subseteq V(I(X))$. 

Thus, if $X=\emptyset $, then $\emptyset =V(A)=V(1)$ is algebraic; if $X\neq \emptyset $, then from (b), $I(X)=Al$ and $X\subseteq V(l)$, so $X$ is algebraic.

(f) From (e), $V(g)\subseteq V(I(V(g)))$. Let $z\in V(I(V(g)))$, from (d), $g\in I(V(g))$, so $g(z)=0$, i.e., $z\in V(g)$. Therefore, $V(I(V(g)))\subseteq V(g)$. 

(g) From (d), $I(X)\subseteq I(V(I(X)))$. From (e), $X\subseteq V(I(X))$, so from (c), $I( V(I(X))\subseteq I(X)$. 

(h) Since $X_k\subseteq \bigcup_{k\in \mathcal{K}}X_k$ for every $k\in \mathcal{K}$, then $I(\bigcup_{k\in \mathcal{K}}X_k)\subseteq I(X_k)$, so $I(\bigcup_{k\in \mathcal{K}}X_k)\subseteq \bigcap_{k\in \mathcal{K}}I(X_k)$. Let $g\in \bigcap_{k\in \mathcal{K}}I(X_k)$ and let $z\in \bigcup_{k\in \mathcal{K}}X_k$, then there exists $k\in \mathcal{K}$ such that $z\in X_k$, then $g(z)=0$, whence $g\in I(\bigcup_{k\in \mathcal{K}}X_k)$. This completes the proof.
\end{proof}

\begin{proposition}\label{proposition2.14}
Let $X\subseteq \mathbb{F}$. Then, 
\begin{enumerate}
	\item[\rm (i)]${\rm rank}(X)\leq |X|$.
	\item[\rm (ii)]If $X=\{z_1,\dots,z_r\}$, then 
	\begin{center}
		$m_X={\rm lclm}(x-z_1, \dots, x-z_r)$.
	\end{center}
	\item[\rm (iii)]Let $Y\subseteq \mathbb{F}$. Then, $m_{X\cup Y}= {\rm lcml}(m_X,m_Y)$ and ${\rm rank}(X\cup Y)\leq {\rm rank}(X)+ {\rm rank}(Y)$. 
	\item[\rm (iv)]Let $g\in A$ of degree $\geq 1$. Then $g$ is a $W$-polynomial if and only if $g={\rm lclm}(x-z_1,\dots,x-z_r)$ for some distinct elements $z_1,\dots ,z_r\in \mathbb{F}$.
	\item[\rm (v)]If $X=\{z_1,\dots,z_r\}$, then ${\rm rank}(X)={\rm rank}({\rm V}_r(z_1,\dots,z_r))$.
	\end{enumerate}
\end{proposition}
\begin{proof}
	(i) We know that $X$ is algebraic. If $X=\emptyset$, then $m_{\emptyset}=1$ and hence ${\rm rank}(X)=0=|X|$. Let $X\neq \emptyset$. The idea is to find a monic polynomial $g\in A$ of degree $\leq m:=|X|$ such that 
	$g(z)=0$ for every $z\in X$. As in \cite{Lam}, the proof is by induction on $m$. For $m=1$, let $X:=\{z\}$, then the statement is trivial taking $g:=x-z$. Let $X:=\{z_1,\dots,z_{m-1},z_m\}$. By induction, there exists a monic polynomial $g'\in A$ of degree $\leq m-1$ such that $g'(z_i)=0$ for every $1\leq i\leq m-1$. Let $g:={\rm lclm}(x-z_m,g')$, then $g(z)=0$ for every $z\in X$ since $x-z_m\mid_r g$ and $g'\mid_r g$; moreover, $g$ is monic, and from (\ref{equation2.4a}), $\deg(g)\leq m$.
	
	(ii) Let $l:={\rm lclm}(x-z_1,\dots,x-z_r)$. It is clear that $X\subseteq V(l)$; if $h\in A$ is monic such that $X\subseteq V(h)$, then $x-z_i\mid_r h$ for every $1\leq i\leq r$, so $l\mid_r h$, whence $\deg(l)\leq \deg(h)$. This proves that $l$ is the monic polynomial of smallest degree such that $X\subseteq V(l)$, i.e., $l=m_X$.
	
	(iii) $X\subseteq V(m_X)$ and $Y\subseteq V(m_Y)$, then 
	\begin{center} 
	$X\cup Y\subseteq V(m_X)\cup V(m_Y)=V(Am_X)\cup V(Am_Y)\subseteq V(Am_X\cap Am_Y)=V({\rm lclm}(m_X,m_Y))$.
	\end{center}
Thus, $X\cup Y\subseteq V(l)$, with $l:={\rm lclm}(m_X,m_Y)$; let $h\in A$ monic such that $X\cup Y\subseteq V(h)$. We have to show that $\deg(h)\geq \deg(l)$. We have $X\subseteq V(h)$ and $Y\subseteq V(h)$, then $m_X\mid_r h$ and $m_Y\mid_r h$. In fact, $h=qm_X+r$ for some $q,r\in A$, with $r=0$ or $\deg(r)<\deg (m_X)$; let $z\in X$, then $h(z)=0=(qm_X)(z)+r(z)=0+r(z)=r(z)$, and this implies that $r=0$ (contrary, there exists a monic polynomial $r'$ of degree $<\deg(m_X)$ such that $X\subseteq V(r')$, a contradiction). Thus, $h=qm_X$, so  $m_X\mid_r h$. Similarly,  $m_Y\mid_r h$. Therefore, $l\mid_r h$, whence $\deg(h)\geq \deg(l)$.

For the second assertion, ${\rm rank}(X\cup Y)=\deg(m_{X\cup Y})=\deg({\rm lcml}(m_X,m_Y))$, but from (\ref{equation2.4a}) we know that $\deg({\rm lcml}(m_X,m_Y)\leq \deg(m_X)+\deg(m_Y)={\rm rank}(X)+{\rm rank}(Y)$.

(iv) $\Rightarrow)$ This part follows from (ii).

$\Leftarrow)$ The idea is to show that $g=m_X$, with $X=\{z_1,\dots,z_r\}$. Since for every $1\leq i\leq r$, $x-z_i\mid_r g$, then $X\subseteq V(g)$; let $h\in A$ monic such that 
$X\subseteq V(h)$, then for every $1\leq i\leq r$, $x-z_i\mid_r h$, hence $g\mid_r h$, so $\deg(g)\leq \deg(h)$. This proves that $g=m_X$. 

(v) We will follow the ideas of \cite{Lam}). Let $V:= {\rm V}_r(z_1,\dots,z_r)$. For proving that ${\rm rank}(X)={\rm rank}({\rm V}_r(z_1,\dots,z_r))$, we will start by noting that ${\rm rank}(X)=r$ if and only if ${\rm rank}({\rm V})=r$. Suppose that ${\rm rank}(X)=r$ and the rows of ${\rm V}$ are not linearly independent, then there exist $c_0, c_1,\ldots c_{r-1}\in \mathbb{F}$, at least one of which is not zero, such that $c_0 {\rm V}_{(1)}+\cdots+c_{r-1}{\rm V}_{(r)}= \textbf{0}$. If $g:= \sum_{i=0}^{r-1}c_ix^{i}$ then $g\neq 0$ and,  using Proposition \ref{proposition2.13}, it follows that $g(z_j)=0$ for each $1\leq j\leq r$. Thus ${\rm rank}(X)<r$, a contradiction. Reciprocally, if ${\rm rank}({\rm V})=r$, then ${\rm V}_{(1)},\ldots, {\rm V}_{(r)}$ are linearly independent. Therefore,  $c_0 {\rm V}_{(1)}+\cdots+c_{r-1}{\rm V}_{(r)}= \textbf{0}$ is true only in the case that $c_0=\cdots=c_{r-1}=0$. As a consequence, given  $g=\sum_{i=1}^{n}c_i x^i$ with $n\leq r-1$, if $g(z_j)=0$ for each  $1\leq j\leq r$, necessarily $g=0$. Whence, ${\rm rank}(X)=r$. Now, let $t:= {\rm rank}({\rm V})$, with $t\leq r-1$. Without lost of generality, it is possible to assume that the first $t$ columns of ${\rm V}$ conform a basis for its column space. Taking $X':=\{z_1,\ldots, z_t\}$ and $g':=m_{X'}$, we have ${\rm rank}(X')\leq t$. We may suppose that $g'=\sum_{i=1}^{r-1} b_ix^i$, adding zero coefficients if it is necessary. Thus, $g'(z_j)=0$ for all $1\leq j\leq t$. But $V^{(t+1)},\ldots, V^{(r)}$ can be written as a linear combination of $V^{(1)},\ldots, V^{(t)}$. Hence  $\begin{pmatrix}b_0 &\cdots& b_{r-1} \end{pmatrix}\cdot {\rm V}^{(j)}=0$, for $t+1\leq j \leq s$. Thus, $g'(z_j)=0$ for all $1\leq j\leq t$ and, therefore, ${\rm rank}(X)={\rm deg}(m_{X})\leq t$. For to show that $t\leq {\rm rank}(X)={\rm deg}(m_{X})$, let $s:={\rm deg}(m_{X})$. Then, $x^i=q_im_X+r_i$, with $r_i=0$ or ${\rm deg}(r_i)\leq {\rm deg}(m_x)$, for $s+1\leq i\leq r$. If $r_{i}=\sum_{k=0}^{s-1}d_{k}^{i}x^i$, evaluating in $z_j$, we obtain that  $N_{i}(z_j)=r_i(x_j)=\sum_{k=0}^{s-1}d_{k}^{i}N_{k}(z_j)$ for $1\leq j\leq r$, that is, ${\rm V}_{(i)}= d_0^{(i)}{\rm V}_{(1)}+\cdots+d_{s-1}^{(i)}{\rm V}_{(s)}$, for each $s+1\leq i\leq r$. Therefore, the row space of ${\rm V}$ is generated by its first $s$ rows. In consequence, $t={\rm rank}({\rm V})\leq s={\rm rank}(X)$. 
\end{proof}

\begin{corollary}
Let $X\subseteq \mathbb{F}$. Then, $I(X)$ is generated by a $W$-polynomial. 
\end{corollary}
\begin{proof}
If $X=\emptyset $, then $I(\emptyset)=A=\langle 1\}$ and $1=m_{\emptyset}$. For $X\neq \emptyset $, the assertion follows from (v)-(b) of Proposition \ref{prop2.12A} and (ii) of Proposition \ref{proposition2.14}. 
\end{proof}

\subsubsection{The bound of a polynomial} 

Following \cite{Gomez-Torrecillas6}, now we will define the bound of a given polynomial of $A$ . Some elementary properties will be considered. The results can be applied in particular to the modulus $f$ (see Definition \ref{definition2.1}), however, in this section $f$ will represent an arbitrary polynomial of $A$. Recall that $A=\mathbb{F}[x;\sigma,\delta]$.

\begin{proposition}\label{proposition2.15}
Let $f\in A$. Then,
\begin{enumerate}
\item[\rm (i)]Let $I$ be the largest two-sided ideal of $A$ contained in $Af$. Then, there exists a two-sided polynomial $f^*\in A$ such that $I=Af^*=f^*A$. $f^*$ is unique up to a non-zero element of $\mathbb{F}$. Moreover, the largest two-sided ideal of $A$ contained in $fA$ coincides with $I$. 
\item[\rm (ii)]Let $f^*\neq 0$. Then, $f^*$ is a two-sided multiple of $f$ of least degree.
\item[\rm (iii)]$Af^*={\rm Ann}_A(A/Af)$.
\item[\rm (iv)]The lattice of left $A$-submodules of $A/Af$ coincides with the lattice of left $A/Af^*$-submodules of $A/Af$.
\end{enumerate}		 
\end{proposition}
\begin{proof}
(i) Since $A$ is a left and right principal ideal domain there exist polynomials $f^*,f'\in A$ such that $I=Af^*=f'A$. If $f^*=0$, then clearly $Af^*=f^*A$. Suppose that $f^*\neq 0$, then $f'\neq 0$ and from $Af^*=f'A$ we get that $\deg(f^*)=\deg(f')$, hence $f'=zf^*$, with $z\in \mathbb{F}^*$, Then, $Af^*=zf^*A$, so $f^*A=z^{-1}Af^*=Af^*$.

Uniqueness of $f^*$: For $f^*=0$ the statement is true.  Let $f^*\neq 0$ and $g\in A$ such that $Ag=I=gA$. Then $Af^*=Ag$, $g\neq 0$ and $g=zf^*$, with $z\in \mathbb{F}^*$. By symmetry, $g=f^*z'$ with $z'\in \mathbb{F}^*$.

Finally, let $J$ be the largest two-sided ideal of $A$ contained in $fA$. By symmetry, there exists $g^*\in A$ such that $J=g^*A=Ag^*$. Let $h\in A$ be a polynomial such that $Ah=Af+Ag^*$. Then, $h=af+bg^{*}$, for some $a,b\in A$. Since $Ag^{*}\subseteq fA$, it follows that $g^*=fa_{1}$, for a non-zero $a_1\in A$. Thus, 
\begin{center}
   $ha_1=afa_1+bg^*a_1=afa_1+ba_{1}'g^*=afa_1+ba_{1}'fa_1=(a+ba_1')fa_1$. 
\end{center}
Therefore, $h=(a+ba_1')f$, i.e., $h\in Af$ and $Ah\subseteq Af$. However, $Ag^*\subseteq Ah$, so that $J=Ag^{*}\subseteq Af$. The maximality of $I$ implies that $J\subseteq I$. Similarly, it is shown that $I\subseteq J$. From the latter it follows that
 $g^*=zf^*$, with $z\in \mathbb{F}^*$.

(ii) Let $h\in A$ be a two-sided multiple of $f$. Then $h=af=fb$, for some $a,b\in A$. By (i), $I=f^*A=Af^*$ is the largest two-sided ideal of $A$ contained in $Af$ and $fA$. Thus, $hA=Ah \subseteq I$ and, therefore, $h=f^* p=qf^{*}$ for certain $p,q\in A$. This implies that $\deg(f^*)\leq \deg(h)$.

(iii) Since $f^*\in Af$, $Af^*\subseteq {\rm Ann}_A(A/Af)$. If $h\in {\rm Ann}_A(A/Af)$, then $h\overline{1}=\overline{0}$, i.e., $h\in Af$. Thus, the two-sided ideal ${\rm Ann}_A(A/Af)$ is contained in $Af$, whence ${\rm Ann}_A(A/Af)\subseteq Af^*$. 

(iv) This follows from the fact that $A/Af$ is a left $A/Af^*$-module since $Af^*(A/Af)=0$.           
\end{proof}
\begin{definition}\label{definition2.24}
	Let $f\in A$. The polynomial $f^*$ is called the\textbf{ bound polynomial of $\boldsymbol{f}$}. Moreover,  $f$ is \textbf{bounded} if $f^*\neq 0$.   
\end{definition}

The part (iv) of the previous proposition explains the importance of the bound polynomial in coding theory. Some other properties involving $f^*$ are presented next.

\begin{proposition}\label{proposition2.17}
Let $f\in A$. 
\begin{enumerate}
	\item[\rm (i)]$f$ is a two-sided polynomial if and only if $f^*=f$. 
	\item[\rm (ii)]If $f$ is bounded, then $\dim_\mathbb{F}(A/Af^*)<\infty$ and $A/Af^*$ is an artinian ring $($left and right$)$. 
	\item[\rm (iii)]If $f$ is irreducible and bounded, then $Af^*$ is maximal and $A/Af^*$ is a simple artinian ring. 
\end{enumerate} 
\end{proposition}
\begin{proof}
(i) $\Rightarrow)$ Since $f^*\in Af$, $\deg(f^*)\geq \deg(f)$, but according to Proposition \ref{proposition2.15} (ii), $\deg(f^*)\leq \deg(f)$. From this we get that $f^*=zf$, with $z\in \mathbb{F}^*$. From \ref{proposition2.15} (i), $f^*=f$.

$\Leftarrow)$ This follows from \ref{proposition2.15} (i).

(ii) Let $\deg(f^*)=n$, since $A$ is a left euclidean domain, $\{1,x,\dots,x^{n-1}\}$ is a $\mathbb{F}$-basis of $A/Af^*$. Any left ideal of $A/Af^*$ is a vector subspace of $A/Af^*$, so $A/Af^*$ is left artinian. Since $\sigma$ is bijective and $Af^*=f^*A$, then considering the coefficients of $A$ on the right side we obtain that $A/Af^*$ is right artinian. 

(iii) $Af^*$ is maximal: $Af^*\neq A$, contrary $Af=A=fA$, so $f$ is invertible, a contradiction. Since $f$ is irreducible, then $Af$ is a maximal left ideal of $A$. In fact, let $Ag$ such that $Af\subseteq Ag$, then $f=pg$, for some $p\in A$, so $p\in A^*$ or $g\in A^*$. In the first case $Ag=Af$ and in the second case $Ag=A$. 

In order to proof of the maximality of $Af^*$, we need the following two preliminary facts.  

(a) Let $I\neq 0$ be a two-sided ideal of $A$. Then, there exists $a^*\in A-\{0\}$ such that $I=a^{*}A=Aa^{*}$: Since $A$ is a domain of left and right principal ideals, then $I=aA=Aa'$, for some $a,a'\in A-\{0\}$. Therefore, there exist $u,v\in A$ such that $a=ua'$ and $a'=av$, so $a=uav$ and $a'=ua'v$. Since $ua\in I=aA$, then $ua=au'$, for some $u'\in A$. Hence, $a=au'v$, so $u'v=1$, i.e., $v\in \mathbb{F}^{*}$. Thus, $aA=Aa'=Aav$, so $aAv^{-1}=Aa$, i.e., $aA=Aa$. Similarly, $u\in \mathbb{F}^{*}$ and $a'A=Aa'$, so we can take $a^*:=a$ or $a^*:=a'$.

(b) Let $I_1, I_2\neq 0$ be two-sided ideals of $A$. If $I_1\subseteq I_2$, then there exists a two-sided ideal $I_3$ of $A$ such that $I_1=I_2I_3$, with $I_1\subseteq I_3$: By (a), $I_1=a^{*}A=Aa^{*}$ and $I_2=b^{*}A=Ab^{*}$. Then, $a^*=b^{*}c$, for some $c\in A$. Given $f\in A$, there exist $f',\overline{f}\in A$ such that $fa^{*}=a^{*}f'$ and $fb^{*}=b^{*}\overline{f}$. Therefore, $b^*\overline{f}c=fb^{*}c=fa^{*}=a^{*}f'= b^{*}cf'$, whence $\overline{f}c=cf'$. Since $f$ runs through $A$, then $f'$ and $\overline{f} $ also run through $A$, hence, $c^{*}:=c$ defines a two-sided ideal $I_3:=c^{*}A=Ac^{*}$ that satisfies $I_1=I_2I_3$, with $I_1\subseteq I_3$.

Now we can complete the proof of the maximality of $Af^*$. Let $L$ be a two-sided ideal of $A$ such that $Af^*\subseteq L$. We have to show that either $L=Af^*$ or $L=A$. By (b), taking $I_1:=Af^*$ and $I_2:=L$, there exists a two-sided ideal $K$ such that $Af^*=LK$, with $Af^*\subseteq K$. If $K\nsubseteq Af$, then $A=Af+K$, so $L=LAf+LK=LAf+Af^*\subseteq Af+Af^*\subseteq Af$, whence $L\subseteq Af^*$, i.e., $L=Af^*$. If $K\subseteq Af$, then $K\subseteq Af^*$, so $K=Af^*$, whence $Af^*=LAf^*$. This implies that $L=A$ since $A$ is a domain and $f^*\neq 0$.  

$A/Af^*$ is a simple artinian ring: This follows from the just proved and (ii). 
\end{proof}

In \cite{Gomez-Torrecillas6} has been computed the bound $f^*$ of a given bounded polynomial $f\in A$ assuming that $A$ is finitely generated as module over its center $Z(A)$.   

\begin{proposition}[\cite{Gomez-Torrecillas6}, Proposition 2.4]\label{proposition2.18}
Assume that $A$ is finitely generated as module over its center by $a_1,\dots, a_r\in A$. Let $0\neq f\in A$. Then, 
$f$ is bounded and 
\begin{center}
$f^*={\rm lclm}(f, f_{a_1},\dots,f_{a_r})$, where $f_{a_i}$ is such that $Af_{a_i}={\rm Ann}_A(\overline{a_i})$, with $\overline{a_i}:=a_i+Af$, $1\leq i\leq r$.
\end{center}
\end{proposition}
\begin{proof}
$f$ is bounded: First observe that $Af\cap Z(A)$ is a non zero two-sided ideal of $A$ contained in $Af$ (see \cite{Gomez-Torrecillas6}, p. 273), so $Af\cap Z(A)\subseteq Af^*$, whence $f^*\neq 0$. 

Now we will show that $f^*={\rm lclm}(f, f_{a_1},\dots,f_{a_r})$. Initially we will prove that $f^*={\rm lclm}(f_{a_1},\dots,f_{a_r})$: We have
\begin{center}
$A{\rm lclm}(f_{a_1},\dots,f_{a_r})=Af_{a_1}\cap\cdots\cap Af_{a_r}={\rm Ann}_A(\overline{a_1})\cap \cdots\cap {\rm Ann}_A(\overline{a_r})$ and $Af^*={\rm Ann}_A(A/Af)$,
\end{center}
hence $Af^*\subseteq A{\rm lclm}(f_{a_1},\dots,f_{a_r})$. On the other hand, let $p\in A{\rm lclm}(f_{a_1},\dots,f_{a_r})$ and $q\in A$, then there exist $q_1,\dots,q_r\in Z(A)$ such that $q=q_1a_1+\cdots+q_ra_r$. Therefore,
\begin{center}
$p\overline{q}=p\overline{q_1a_1}+\cdots+p\overline{q_ra_r}=q_1p\overline{a_1}+\cdots+q_rp\overline{a_r}=\overline{0}$,
\end{center} 
so $p\in {\rm Ann}_A(A/Af)=Af^*$. Hence, $Af^*=A{\rm lclm}(f_{a_1},\dots,f_{a_r})$, i.e., $f^*={\rm lclm}(f_{a_1},\dots,f_{a_r})$.

Now, since $f|_l f^*$, then $f^*={\rm lclm}(f,f_{a_1},\dots,f_{a_r})$.
\end{proof}

\subsubsection{$\mathbb{F}$-linear evaluation codes}

There are two key types of $\mathbb{F}$-linear codes that we present next.  

\begin{definition}
Let $r\geq 1$ and $k\in \{1,\dots,r\}$.
\begin{enumerate}
\item[\rm (i)]Let $Z:=(z_1,\dots,z_r)\in \mathbb{F}^r$ such that ${\rm rank}({\rm V}_r(Z))\geq k$. The \textbf{remainder
evaluation code of length $\boldsymbol{r}$} and\textbf{ support $\boldsymbol{Z}$} is defined as
\begin{center} 
$\mathfrak{C}_k(Z):=\{(g(z_1),\dots,g(z_r)\in \mathbb{F}^r)\mid g\in A,\deg(g)\leq k-1\}$.
\end{center}
\item[\rm (ii)]Let $Z:=(z_1,\dots,z_r)\in \mathbb{F}^r$ such that ${\rm rank}({\rm Wr}_r(Z))\geq k$. The \textbf{operator
evaluation code of length $\boldsymbol{r}$} and\textbf{ support $\boldsymbol{Z}$} is defined as
\begin{center} 
	$\mathfrak{C}_{k,\mathcal{L}}(Z):=\{(\mathcal{L}_g(z_1),\dots,\mathcal{L}_g(z_r)\in \mathbb{F}^r)\mid g\in A,\deg(g)\leq k-1\}$.
\end{center} 
\end{enumerate} 
\end{definition}

\begin{proposition}\label{proposition2.27}
Let $r\geq 1$, $Z:=(z_1,\dots,z_r)\in \mathbb{F}^r$ and $k\in \{1,\dots,r\}$. Then, 
\begin{enumerate}
\item[\rm (i)]$\mathfrak{C}_k(Z)$ is a $\mathbb{F}$-linear code of dimension $k$.
\item[\rm (ii)]$\mathfrak{C}_{k,\mathcal{L}}(Z)$ is a $\mathbb{F}$-linear code of dimension $k$.  
\end{enumerate} 	
\end{proposition}
\begin{proof}
(i) Since for $g,g'\in A$, $z,z'\in \mathbb{F}$, we have $(g+g')(z)=g(z)+g'(z)$ and $(z'g)(z)=z'g(z)$, then $\mathfrak{C}_k(Z)$ is a $\mathbb{F}$-subspace of $\mathbb{F}^r$, i.e., $\mathfrak{C}_k(Z)$ is a $\mathbb{F}$-linear code. Next we will show that $\dim_\mathbb{F}(\mathfrak{C}_k(Z))=k$. Observe first that $\mathfrak{C}_k(Z)$ is generated as $\mathbb{F}$-space by $k$ vectors:
\begin{center}
$(1(z_1),\dots,1(z_r)),(x(z_1),\dots,x(z_r))\dots, (x^{k-1}(z_1),\dots,x^{k-1}(z_r))$,
\end{center}
thus, we have the matrix
\begin{center}
$M:=\begin{bmatrix}
1(z_1) & \cdots & 1(z_r)\\
x(z_1) & \cdots & x(z_r)\\
\vdots & \cdots & \vdots\\
x^{k-1}(z_1) & \cdots & x^{k-1}(z_r)
\end{bmatrix}=
\begin{bmatrix}
1 & \cdots & 1\\
z_1 & \cdots & z_r\\
\vdots & \cdots & \vdots\\
x^{k-1}(z_1) & \cdots & x^{k-1}(z_r)
\end{bmatrix}$.
\end{center}
Observe that the rows of $M$ are the first $k$ rows of ${\rm V}_r(Z)$, so ${\rm rank}(M)={\rm min}\{l,k\}$, where $l:={\rm rank}({\rm V}_r(Z))$, hence $l\geq k$. Thus, ${\rm rank}(M)=k$. This means that $\dim_\mathbb{F}(\mathfrak{C}_k(Z))=k$.

(ii) For $g,g'\in A$, $z,z'\in \mathbb{F}$, we have $\mathcal{L}_g(z)+\mathcal{L}_{g'}(z)=\mathcal{L}_{g+g'}(z)$ and $(z'\mathcal{L}_g)(z)=\mathcal{L}_{z'g}(z)$, so $\mathfrak{C}_{k,\mathcal{L}}(Z)$ is a $\mathbb{F}$-subspace of $\mathbb{F}^r$, i.e., $\mathfrak{C}_{k,\mathcal{L}}(Z)$ is a $\mathbb{F}$-linear code. Moreover,  $\mathfrak{C}_{k,\mathcal{L}}(Z)$ is generated as $\mathbb{F}$-space by $k$ vectors:
\begin{center}
	$(z_1,\dots,z_r),(\mathcal{D}(z_1),\dots,(\mathcal{D}(z_r))\dots, ((\mathcal{D}^{k-1}(z_1),\dots,(\mathcal{D}^{k-1}(z_r))$,
\end{center} 
then we have the matrix
\begin{center}
	$N:=\begin{bmatrix}
	z_1 & \cdots & z_r\\
	\mathcal{D}(z_1) & \cdots & \mathcal{D}(z_r)\\
	\vdots & \cdots & \vdots\\
	\mathcal{D}^{k-1}(z_1) & \cdots & \mathcal{D}^{k-1}(z_r)
	\end{bmatrix}$.
\end{center}
Observe that the rows of $N$ are the first $k$ rows of ${\rm Wr}_r(Z)$, so ${\rm rank}(N)={\rm min}\{s,k\}$, where $s:={\rm rank}({\rm Wr}_r(Z))$, hence $s\geq k$. Thus, ${\rm rank}(N)=k$. This means that $\dim_\mathbb{F}(\mathfrak{C}_{k,\mathcal{L}}(Z))=k$.
\end{proof}

\subsubsection{Distance of $\mathbb{F}$-linear codes}

In this subsection we recall two classical notions of distance for $\mathbb{F}$-linear codes. We will follow \cite{Boucher} and \cite{Gabidulin}.

\begin{definition}
Let $r\geq 1$ and $\mathbb{F}^{\sigma}:=\{z\in \mathbb{F}\mid \sigma(z)=z\}$. 
\begin{enumerate}
\item[\rm (i)]If $Z:=(z_1,\dots,z_r)\in \mathbb{F}^r$, the \textbf{rank} of $Z$, denoted ${\rm rank}(Z)$, is the dimension of the $\mathbb{F}^{\sigma}$-vector subspace of $\mathbb{F}$ spanned by $z_1,\dots,z_r$. 
\item[\rm (ii)]For $Z=(z_1,\dots,z_r),Z'=:=(z_1',\dots,z_r')\in \mathbb{F}^r$ the \textbf{rank distance} is defined by
\begin{center}
$d_{{\rm rank}}(Z,Z'):={\rm rank}(Z-Z')$
\end{center}
and the \textbf{Hamming distance} is defined by
\begin{center}
$d_H(Z,Z'):=\mid \{i\in \{1,\dots,r\}\mid z_i\neq z_i'\}\mid$.
\end{center} 
\item[\rm (iii)]If $\mathfrak{C}$ is a $\mathbb{F}$-linear code of $\mathbb{F}^r$, the \textbf{minimal rank distance} of $\mathfrak{C}$ is defined by
\begin{center}
$d_{{\rm rank}}(\mathfrak{C}):={\rm min}\{d_{{\rm rank}}(Z,Z')\mid Z,Z'\in \mathfrak{C}, Z\neq Z'\}$
\end{center}
and the \textbf{minimal Hamming distance} of $\mathfrak{C}$ is defined by 
\begin{center}
	$d_H(\mathfrak{C}):={\rm min}\{d_H(Z,Z')\mid Z,Z'\in \mathfrak{C}, Z\neq Z'\}$.
\end{center}  
\item[\rm (iv)]Let $\mathfrak{C}$ be a $\mathbb{F}$-linear code of $\mathbb{F}^r$. $\mathfrak{C}$ is a $\textbf{MRD}$ $($\textbf{Maximum Rank Distance}$)$ code if
\begin{center}
$d_{{\rm rank}}(\mathfrak{C})=r-\dim_\mathbb{F}(\mathfrak{C})+1$.	
\end{center}

$\mathfrak{C}$ is a $\textbf{MDS}$ $($\textbf{Maximum Distance Separable}$)$ code if
\begin{center}
$d_H(\mathfrak{C})=r-\dim_\mathbb{F}(\mathfrak{C})+1$. 
\end{center}
\end{enumerate}
\end{definition}

Some properties of the defined distances are presented in the next proposition. 

\begin{proposition}
With the notation of the previous definition,
\begin{enumerate}
\item[\rm (i)]$d_{{\rm rank}}(Z,Z')\leq d_H(Z,Z')$, for all $Z,Z'\in \mathbb{F}^r$.
\item[\rm (ii)]$d_{{\rm rank}}(Z,Z')\geq 0$; $d_{{\rm rank}}(Z,Z')=0$ if and only if $Z=Z'$; $d_{{\rm rank}}(Z,Z')=d_{{\rm rank}}(Z',Z)$; $d_{{\rm rank}}(Z,Z')\leq d_{{\rm rank}}(Z,Z'')+d_{{\rm rank}}(Z'',Z')$, for every $Z''\in \mathbb{F}^r$. Thus, $d_{{\rm rank}}$ is a metric over $\mathbb{F}^r$.
\item[\rm (iii)]$d_H(Z,Z')\geq 0$; $d_H(Z,Z')=0$ if and only if $Z=Z'$; $d_H(Z,Z')=d_H(Z',Z)$; $d_H(Z,Z')\leq d_H(Z,Z'')+d_H(Z'',Z')$, for every $Z''\in \mathbb{F}^r$. Thus, $d_H$ is a metric over $\mathbb{F}^r$.
\item[\rm (iv)]Let $\mathfrak{C}$ be a $\mathbb{F}$-linear code of $\mathbb{F}^r$ with $\dim_\mathbb{F}(\mathfrak{C})=k$. Then,
\begin{center}$d_H(\mathfrak{C})\leq r-k+1$ and $d_{{\rm rank}}(\mathfrak{C})\leq r-k+1$.
\end{center}
\item[\rm (v)]Let $\mathfrak{C}$ be a $\mathbb{F}$-linear code of $\mathbb{F}^r$ with $\dim_\mathbb{F}(\mathfrak{C})=k$. $\mathfrak{C}$ is $MRD$ if and only if for every matrix $Y\in M_{(r-k)\times r}(\mathbb{F}^{\sigma})$ of rank $r-k$, the rank of $YG_{X^{\perp}}(\mathfrak{C}^{\perp})^T$ over $\mathbb{F}$ is $r-k$. Moreover, $\mathfrak{C}$ is $MRD$ if and only if $\mathfrak{C}^{\perp}$ is $MRD$. 
\item[\rm (vi)]Let $\mathfrak{C}$ be a $\mathbb{F}$-linear code of $\mathbb{F}^r$ with $\dim_\mathbb{F}(\mathfrak{C})=k$. $\mathfrak{C}$ is $MDS$ if and only if any $r-k$ columns of $G_{X^{\perp}}(\mathfrak{C}^{\perp})$ are linearly independent. Moreover, $\mathfrak{C}$ is $MDS$ if and only if $\mathfrak{C}^{\perp}$ is $MDS$.
\item[\rm (vii)]Let $Z:=(z_1,\dots,z_r)\in \mathbb{F}^r$ such that ${\rm rank}({\rm V}_r(Z))=r$. Then, for every $k\in \{1,\dots,r\}$, $\mathfrak{C}_k(Z)$ is $MDS$.
\item[\rm (viii)]Let $Z:=(z_1,\dots,z_r)\in \mathbb{F}^r$ such that $z_1,\dots,z_r$ are linearly independent over $\mathbb{F}^{\sigma}$. Then, for every $k\in \{1,\dots,r\}$,	$\mathfrak{C}_{k,\mathcal{L}}(Z)$ is $MRD$.
\end{enumerate}
\end{proposition}
\begin{proof}
(i)-(iii) are evident.

(iv) Since $\dim_\mathbb{F}(\mathfrak{C})=k$, there exists at least one word $Z\neq 0$ in $\mathfrak{C}$ such that the number of its non zero entries is $\leq r-k+1$. In fact, consider a generator matrix $G_X(\mathfrak{C})$ of $\mathfrak{C}$,
\begin{center}
	$G_X(\mathfrak{C})=\begin{bmatrix}
	v_{11} & \cdots & v_{1r}\\
	\vdots & \cdots & \vdots\\
	v_{k1} & \cdots & v_{kr}
	\end{bmatrix}$.
\end{center}
In the first row of $G_X(\mathfrak{C})$ there is at least one non zero entry, say, $v_{1j_1}$. By elementary row operations we can assume that
\begin{center}
	$G_X(\mathfrak{C})=\begin{bmatrix}
	v_{11} & \cdots & 1 & \cdots & v_{1r}\\
	\vdots & \cdots & 0 & \vdots\\
	v_{k1} & \cdots & 0 & \cdots & v_{kr}
	\end{bmatrix}$,
\end{center} 
where the $j_1$-th column of $G_X(\mathfrak{C}$ is $(1,0,\dots,0)^T$. We can repeat this procedure for any other row; observe that $j_1, \dots, j_k$ are different. Thus, considering any row of $G_X(\mathfrak{C})$, for example, the first one, we conclude that the number of non zero entries in this first row is $\leq r-k+1$. This proves the claimed. Therefore, there exists $0\neq Z$ in $\mathfrak{C}$ such that $d_H(Z,0)\leq r-k+1$, whence  $d_H(\mathfrak{C})\leq r-k+1$.

Now, from (i), $d_{{\rm rank}}(Z,Z')\leq d_H(Z,Z')$, for all $Z,Z'\in \mathfrak{C}$, so $d_{{\rm rank}}(\mathfrak{C})\leq d_{H}(\mathfrak{C})\leq r-k+1$. 

(v) See \cite{Gabidulin}, Theorem 3.

(vi) $\Rightarrow)$: By the hypothesis, $d_H(\mathfrak{C})=r-k+1$, then every nonzero word of $\mathfrak{C}$ has at least $r-k+1$ non zero entries. Assume that $G_{X^{\perp}}(\mathfrak{C}^{\perp})$ has $r-k$ columns linearly dependent, then there exists $0\neq Z\in \mathbb{F}^r$ such that $G_{X^{\perp}}(\mathfrak{C}^{\perp})Z^T=0$, where the number of non zero entries of $Z$ is $\leq r-k$. But since $G_{X^{\perp}}(\mathfrak{C}^{\perp})$ is a check parity matrix of $\mathfrak{C}$, then $Z\in \mathfrak{C}$ (Proposition \ref{proposition2.7}), a contradiction. 

$\Leftarrow)$: Assume that there exist $Z,Z'\in \mathfrak{C}$ such that $d_H(Z,Z')<r-k+1$. We have $G_{X^{\perp}}(\mathfrak{C}^{\perp})Z=0=G_{X^{\perp}}(\mathfrak{C}^{\perp})Z'$, so $G_{X^{\perp}}(\mathfrak{C}^{\perp})(Z-Z')=0$. But the number of non zero entries of $Z-Z'$ is $<r-k+1$, then $G_{X^{\perp}}(\mathfrak{C}^{\perp})$ has $r-k$ columns linearly dependent.

Since $(\mathfrak{C}^{\perp})^{\perp}=\mathfrak{C}$, for the second statement it is enough to prove that if $\mathfrak{\mathfrak{C}}$ is $MDS$, then $\mathfrak{C}^{\perp}$ is $MDS$. We have to show that any $k$ columns of $G_{X}(\mathfrak{C})$ are linearly independent. Contrary, assume that $G_{X}(\mathfrak{C})$ has $k$ columns linearly dependent. Consider the submatrix $M$ of size $k\times k$ formed by these columns. Then ${\rm rank}(M)<k$, whence the rows of $M$ are linearly dependent, i.e., there exists a linear combination of the rows of $M$ that is zero, and hence, for the same linear combination on the rows of $G_{X}(\mathfrak{C})$, we have that at least $k$ entries in this combination are zero. This implies that we have a word in $\mathfrak{C}$ where the number of non zero entries is $\leq r-k$. This is a contradiction since $\mathfrak{C}$ is $MDS$. 

(vii) From (iv) we know that for every $k\in \{1,\dots,r\}$, $d_H(\mathfrak{C}_k(Z))\leq r-k+1$. Suppose that there exists $k\in \{1,\dots,r\}$ such that $d_H(\mathfrak{C}_k(Z))<r-k+1$, then there exist $Z',Z''\in \mathfrak{C}_k(Z)$ such that $d_H(Z',Z'')<r-k+1$, so $Z'-Z''\in \mathfrak{C}_k(Z)$ is such that at least $k$ coordinates of $Z'-Z''$ are null, i.e., there exists $g\in A$ such that $\deg(g)\leq k-1$ and $g(z_{i_1})=\cdots=g(z_{i_k})=0$. Then, $(x-z_{i_j})\mid_r g$ for every $1\leq j\leq k$, whence ${\rm lclm}(x-z_{i_1},\dots,x-z_{i_k})\mid_r g$. Since ${\rm rank}({\rm V}_r(Z))=r$, then ${\rm rank}({\rm V}_r(z_{i_1},\dots,z_{i_k}))=k$, but from Proposition \ref{proposition2.14}, for the algebraic set $X:=\{z_{i_1},\dots,z_{i_k}\}$, we have
\begin{center}
 ${\rm rank}({\rm V}_r(z_{i_1},\dots,z_{i_k}))={\rm rank}(X)=\deg(m_X)=\deg({\rm lclm}(x-z_{i_1},\dots,x-z_{i_k}))$,
\end{center}
so $\deg(g)\geq k$, a contradiction. This completes the proof of (vii). 
  
(viii) By Proposition \ref{proposition2.7}, $(\mathfrak{C}_{k,\mathcal{L}}(Z)^{\perp})^{\perp}=\mathfrak{C}_{k,\mathcal{L}}(Z)$, hence, from (v), we will show that for every $k\in \{1,\dots,r\}$, $\mathfrak{C}_{k,\mathcal{L}}(Z)^{\perp}$ is $MRD$. Thus, we will prove that $d_{{\rm rank}}(\mathfrak{C}_{k,\mathcal{L}}(Z)^{\perp})=r-\dim_\mathbb{F}(\mathfrak{C}_{k,\mathcal{L}}(Z)^{\perp})+1=r-(r-k)+1=k+1$, for every $k\in \{1,\dots,r\}$. The idea is to show that for every $k\in \{1,\dots,r\}$, there is not a word in $\mathfrak{C}_{k,\mathcal{L}}(Z)^{\perp}$ of rank $<k+1$ over $\mathbb{F}^{\sigma}$. Contrary, assume that there is a $k\in \{1,\dots,r\}$  and a word $Z':=(z_1',\dots,z_r')\in \mathfrak{C}_{k,\mathcal{L}}(Z)^{\perp}$ of rank $t<k+1$, i.e., $\dim_{\mathbb{F}^{\sigma}}\langle z_1',\dots,z_r'\rangle =t<k+1$; let $\{u_1,\dots,u_t\}$ be a $\mathbb{F}^{\sigma}$-basis of $\langle z_1',\dots,z_r'\rangle$, then there exists a matrix $M\in M_{t\times r}(\mathbb{F}^{\sigma})$ such that $Z'=(u_1,\dots,u_t)M$, and from Proposition \ref{proposition2.7}, $G_{X}(\mathfrak{C}_{k,\mathcal{L}}(Z))Z'^T=0$, i.e., $G_{X}(\mathfrak{C}_{k,\mathcal{L}}(Z))M^T(u_1,\dots,u_t)^T=0$, but according to the proof of Proposition \ref{proposition2.27}, $G_{X}(\mathfrak{C}_{k,\mathcal{L}}(Z))={\rm Wr}_{k,r}(Z)$ (the rectangular Wronskian matrix, defined as in Definition \ref{definition2.16}, but with only $k$ rows; $k\leq r$), but since $t\leq k$, ${\rm Wr}_{t,r}(Z)M^T(u_1,\dots,u_t)^T=0$. Let $Z'':=(z_1'',\dots,z_t''):=ZM^T$; since $\mathcal{D}$ is $\mathbb{F}$-linear (see Definition \ref{definition2.16}), we obtain that ${\rm Wr}_{t,r}(Z)M^T={\rm Wr}_{t}(Z'')$, hence ${\rm Wr}_{t}(Z'')(u_1,\dots,u_t)^T=0$. Now, as $z_1,\dots,z_r$ are linearly independent over $\mathbb{F}^{\sigma}$ and the rank of $M$ over $\mathbb{F}^{\sigma}$ is $t$ (the fact that the rank of $M$ is $t$ can be explained in the following way: there exists a matrix $N\in M_{r\times t}(\mathbb{F}^{\sigma})$ such that $U=Z'N$, where $U:=(u_1,\dots,u_t)$, so $MN=I_t$, whence, $t={\rm rank}(MN)\leq {\rm min}\{{\rm rank}(M),{\rm rank}(N)\}$, but ${\rm rank}(M)\leq t$, ${\rm rank}(N)\leq t$), then $z_1'',\dots,z_t''$ are linearly independent over $\mathbb{F}^{\sigma}$, hence $\det({\rm Wr}_{t}(Z''))\neq 0$, a contradiction. 
  
\end{proof}

\subsection{The case $A:=\mathbb{F}[x;\sigma]$}

Let $\mathbb{F}[x;\sigma,\delta]$ be as in the previous subsections. As we observed at the beginning of subsection \ref{subsection2.1}, since $\mathbb{F}$ is a finite field then $\delta$ necessarily is a $\sigma$-inner derivation, but in such case it is well-known that $\mathbb{F}[x;\sigma,\delta]$ is isomorphic to $\mathbb{F}[y;\sigma]$, where $y:=x-w$, with $w$ as in (\ref{equ2.3}) (see \cite{McConnell} or also \cite{Lezama9}, Chapter 1). Thus, over finite fields, the Ore extension $\mathbb{F}[x;\sigma,\delta]$ can be reduced to one of automorphism type, i.e, where the $\sigma$-derivation is trivial. In this subsection we will study the case $A:=\mathbb{F}[x;\sigma]$ and we fix the modulus $f\in A$ of degree $n\geq 1$. Of course all results of the previous subsection can be applied in this situation, but some additional properties will be added. We will mainly follow \cite{Luerssen}, and also, \cite{Gomez-Torrecillas} and \cite{Gomez-Torrecillas4}.

\subsubsection{Two-sided polynomials}

The next proposition establishes conditions under which a given polynomial of $A$ is two-sided. For this, recall that $Z(A)=\mathbb{F}^{\sigma}[x^s]$, where $s$ is the order of $\sigma$ (see \cite{Lezama-sigmaPBW} or also \cite{Lezama9}).

\begin{proposition}\label{proposition2.30}
$g\in A$ is a two-sided polynomial if and only if $g=cx^th$, for some $c\in \mathbb{F}$, $t\geq 0$ and $h\in Z(A)$. In particular, for $a\in \mathbb{F}^*$, $x^n-a$ is a two-sided polynomial if and only if $x^n-a\in Z(A)$ if and only if $\sigma(a)=a$ and $s\mid n$, where $a\in \mathbb{F}$ and $s$ is the order of $\sigma$.   
\end{proposition}
\begin{proof}
$\Rightarrow)$: If $g=0$, then we can take $c=0$, $t=0$, $h=1$. Let $g:=g_0+\cdots +g_mx^m\neq 0$, with $g_m\neq 0$. We can assume that $g$ is monic. In fact, $g=g_m(g_m^{-1}g_0+\cdots+x^m)=g_mg'$, with $g':=g_m^{-1}g_0+\cdots+x^m$. We have $Ag=gA$, so, since $Ag_m=g_mA$, we get $Ag_mg'=g_mg'A=g_mAg'$, but since $g_m\neq 0$, then $g'A=Ag'$. Therefore, if we prove the claimed for monic polynomials, then $g'=c'x^th$, with $c'\in \mathbb{F}$, $t\geq 0$ and $h\in Z(A)$, so $g=cx^th$, with $c=g_mc'$.  

We will prove the claimed by induction on $m$. If $m=0$, then $c=g=1\in \mathbb{F}$, $t=0$ and $h=1$. If $m=1$, then $g=g_0+x$; if $g_0=0$, then we get the claimed with $c=1$, $t=1$ and $h=1$. Assume that $g_0\neq 0$; if $\sigma=i_\mathbb{F}$, then $A$ is commutative and $g\in Z(A)$, thus $c=1,t=0, h=g$. Hence, assume that $g_0\neq 0$ and $\sigma\neq i_\mathbb{F}$. Then, there exists $z\in \mathbb{F}$ such that $\sigma(z)\neq z$; consider $gz\in gA=Ag$, then $(g_0+x)z=z'(g_0+x)$, with $z'\in \mathbb{F}$, from this we get that $g_0z+\sigma(z)x=z'g_0+z'x$, so $z=z'$ and $\sigma(z)=z'=z$, a contradiction. This completes the proof of case $m=1$.

Assume the claimed proved for non zero monic two-sided polynomials of degree $\leq m-1$. We will consider two possible cases. 

\textit{Case 1}. $g_0\neq 0$. If $\sigma=i_\mathbb{F}$, then $A$ is commutative and $g\in Z(A)$, thus $c=1,t=0, h=g$. Now assume that $\sigma\neq i_\mathbb{F}$. Let $G$ be the set of non zero coefficients of $g$. Arise two cases.

\textit{Case 1.1}. There exists $g_l\in G$ such that $\sigma^l\neq i_\mathbb{F}$. Then, there exists $u\in \mathbb{F}$ such that $\sigma^l(u)\neq u$, so $gu=u'g$, for some $u'\in \mathbb{F}$, from this we get that $u=u'$ and $\sigma^l(u)=u'=u$, a contradiction. 

\textit{Case 1.2}. For every $g_l\in G$, $\sigma^l=i_\mathbb{F}$. Let $s$ be the order of $\sigma$, then $s\mid l$. In particular, $m=ts$, with $t\in \mathbb{N}$. Thus, $ g=g_0+g_sx^s+g_{2s}x^{2s}+\cdots+g_{st}x^{ts}$, where $G\subseteq \{g_0,g_s,g_{2s},\dots,g_{ts}\}$, with $g_{ts}=1$. We want to show that $\sigma(g_l)=g_l$ for every $g_l\in G$, if so, then $g\in Z(A)$, and we finish the case 1. Assume that there exists $g_l\in G$ such that $\sigma(g_l)\neq g_l$, with $l=rs$, for some $0\leq r\leq t-1$. Consider $gx\in gA$, then $gx=(p_0+x)g$, for some $p_0\in \mathbb{F}$ (recall that $g$ is monic). From this we get that $\sigma(g_{rs})=g_{rs}$, i.e., $\sigma(g_l)=g_l$, a contradiction. Thus, $g\in Z(A)$, $c=1$ and $t=0$. 

\textit{Case 2}. $g_0=0$. Let $l$ be minimum such that $g_l\neq 0$. Then $g=x^{l}g'$, with $l\geq 1$ and $0\neq g'\in A$ monic. Observe that $x^l$ and $g'$ are two-sided: In fact, $x^lA=Ax^l$ and $Ag=Ax^lg'=x^lAg'=x^lg'A=gA$, but since $A$ is a domain, then $Ag'=g'A$. By induction, there exist $c'\in \mathbb{F}$, $l'\geq 0$ and $h\in Z(A)$ such that $g'=c'x^{l'}h$, hence $g=x^{l}c'x^{l'}h=\sigma^l(c')x^{l+l'}h$, so we obtain the claimed with $c:=\sigma^l(c')$ and $t:=l+l'$.    

$\Leftarrow)$: If $c=0$, then clearly $Ag=gA$. Let $c\neq 0$ and $p=p_0+p_1x+\cdots+p_lx^l\in A$, then
\begin{center} 
$gp=cx^thp=cx^tph=c(\sigma^t(p_0)x^t+\sigma^t(p_1)x^{t+1}+\cdots+\sigma^t(p_l)x^{t+l})h=[\sigma^t(p_0)cx^t+\sigma^t(p_1)cxc^{-1}cx^t+\cdots+\sigma^t(p_l)cx^lc^{-1}cx^t]h=[\sigma^t(p_0)+\sigma^t(p_1)cxc^{-1}+\cdots+\sigma^t(p_l)cx^lc^{-1}]cx^th=[\sigma^t(p_0)+\sigma^t(p_1)cxc^{-1}+\cdots+\sigma^t(p_l)cx^lc^{-1}]g$,  
\end{center}
so $gA\subseteq Ag$. Now,
\begin{center} 
	$pg=(p_0+p_1x+\cdots+p_lx^l)cx^th=h(p_0+p_1x+\cdots+p_lx^l)cx^t=cx^th[\sigma^{-t}(p_0)+\sigma^{-t}(p_1c^{-1}\sigma(c))x+\cdots+\sigma^{-t}(p_lc^{-1}\sigma^l(c))x^l]=g[\sigma^{-t}(p_0)+\sigma^{-t}(p_1c^{-1}\sigma(c))x+\cdots+\sigma^{-t}(p_lc^{-1}\sigma^l(c))x^l]$,  
\end{center}
so $Ag\subseteq gA$.

Now we consider the particular case when $g$ is the modulus $f$, with $f:=x^n-a$, $a\in \mathbb{F}^*$. If $x^n-a$ is a two-sided polynomial, we will show that $x^n-a$ is central. For $b\in \mathbb{F}$, $(x^n-a)b=p(x^n-a)$, with $p\in A$, then necessarily $p\in \mathbb{F}$ and $ab=pa$, so $p=b$, i.e., $x^n-a$ commutes with any element of $\mathbb{F}$; now, $(x^n-a)x=q(x^n-a)$ for some $q\in A$, then $\deg(q)=1$, so $q=q_0+q_1x$, and from this necessarily $q_1=1$ and $q_0=0$, i.e., $x^n-a$ commutes with $x$. This implies that $x^n-a$ is central. Conversely, if $x^n-a\in Z(A)$, then obviously $x^n-a$ is a two-sided polynomial. The last assertion of the proposition is trivial.
\end{proof}

\subsubsection{Similarity of polynomials}

\begin{definition}
Let $g=g_0+\cdots +x^m\in A$ be a monic polynomial, $m\geq 1$. The \textbf{companion matrix} of $g$ is defined by 
\begin{center}
$C_g:=\begin{bmatrix}
0 & 1 & 0 & \cdots & 0\\
0 & 0 & 1 & \cdots & 0\\
\vdots & \vdots & & \ddots & \vdots\\
0 & 0 & 0 & \cdots & 1\\
-g_0 & -g_1 & -g_2 & \cdots & -g_{m-1}
\end{bmatrix}\in M_{m\times m}(\mathbb{F})$.
\end{center} 
\end{definition}

\begin{proposition}[\cite{Lam-Leroy}, Theorem 4.9]
Let $g,h\in A$ be monic polynomials of degree $m\geq 1$. Then, $g\approx_l h$ if and only if there exists a matrix $B\in GL_m(\mathbb{F})$ such that $C_g=\sigma(B)C_hB^{-1}$.   
\end{proposition}
\begin{proof}
$\Rightarrow)$: Let $\alpha: A/Ag\to A/Ah$ be an isomorphism of left $A$-modules. Then $\alpha$ is an $\mathbb{F}$-isomorphism. Consider the function $S_g:A/Ag\to A/Ag$, $\overline{a}\mapsto \overline{xa}$, for $a\in A$. Observe that $S_g$ is a well-defined $\boldsymbol{\sigma}$-\textit{\textbf{semi-linear homomorphism}}, i.e., $S_g(z\overline{a})=\sigma(z)S_g(\overline{a})$, for $z\in \mathbb{F}$. Despite of $S_g$ is not $\mathbb{F}$-linear, notice that $C_g$ is the matrix of $S_g$ in the canonical $\mathbb{F}$-basis $X_g:=\{\overline{1}, \overline{x}, \dots,\overline{x^{m-1}}\}$ of $A/Ag$ (we dispose the scalars by rows, see \cite{Lam-Leroy}, p.10). Similarly, $C_h$ is the matrix of $S_h:A/Ah\to A/Ah$ in the canonical $\mathbb{F}$-basis $X_h:=\{\overline{1}, \overline{x}, \dots,\overline{x^{m-1}}\}$ of $A/Ah$. Let $B$ be the matrix of $\alpha$ in the canonical $\mathbb{F}$-bases $X_g$ and $X_h$, then $B\in GL_m(\mathbb{F})$. We want to show that $C_gB=\sigma(B)C_h$. Since $\alpha$ is $A$-linear we have $\alpha\circ S_g=S_h\circ \alpha$. In fact, let $a\in A$, then $\alpha\circ S_g(\overline{a})=\alpha(\overline{xa})=\alpha(x\overline{a})=x\alpha(\overline{a})$ and $S_h\circ \alpha(\overline{a})=S_h(\overline{b})=\overline{xb}=x\overline{b}=x\alpha(\overline{a})$. Let $\overline{a}=a_0\overline{1}+a_1\overline{x}+\cdots+a_{m-1}\overline{x^{m-1}}\in A/Ag$, with $a_i\in \mathbb{F}$, $0\leq i\leq m-1$, then $\alpha(S_g(\overline{a}))=S_h(\alpha(\overline{a})) $ and we can express this equality by the coordinates in the canonical bases, i.e., in a matrix form: 
\begin{center} 
$(\sigma(a_0),\dots,\sigma(a_{m-1}))C_gB=(\sigma(a_0),\dots,\sigma(a_{m-1}))\sigma(B)C_h$, for every $(a_0,\dots,a_{m-1})$. 
\end{center} 
As $\sigma$ is bijective, then $C_gB=\sigma(B)C_h$.  

$\Leftarrow)$: Now we assume that there exists a matrix $B\in GL_m(\mathbb{F})$ such that $C_g=\sigma(B)C_hB^{-1}$. Let $\alpha: A/Ag\to A/Ah$ be the $\mathbb{F}$-linear isomorphism corresponding to $B$ in the canonical $\mathbb{F}$-bases as before. We want to show that $\alpha$ is $A$-linear. From $C_gB=\sigma(B)C_h$ we get that $\alpha \circ S_g=S_h\circ \alpha$, but this means that $\alpha$ is $A$-linear, so $g\approx_l h$.   
\end{proof}

\subsubsection{Skew polynomials and linearized polynomials}

Recall that $|\mathbb{F}|=q^k$, where $k:=\dim_{\mathbb{Z}_q}(\mathbb{F})$ and $\mathbb{Z}_q$ is the prime subfield of $\mathbb{F}$. In this subsection we will assume that $\sigma=\phi$ is the Frobenius automorphism. Related to Proposition \ref{proposition2.14A} we have the following subset of the commutative polynomial ring $\mathbb{F}[y]$. 

\begin{definition}
Let $\mathfrak{L}:=\{\sum_{i=0}^m g_iy^{q^i}\mid g_i\in \mathbb{F}, m\geq 0\}$. The elements of $\mathfrak{L}$ are called $\boldsymbol{q}$-\textbf{linearized polynomials}.  
\end{definition}

Observe that if $g\in A$, the function $\mathbb{F} \to \mathbb{F}$ defined by $z \mapsto g(z)$ is not, in general, $\mathbb{Z}_q$-linear. But if $g\in \mathfrak{L}$, then the function, also denoted by $g$, and defined by 
\begin{align*}
\mathbb{F} & \xrightarrow{g} \mathbb{F}\\
z & \mapsto g(z) 
\end{align*}
is $\mathbb{Z}_q$-linear. In fact, this function is additive since
\begin{center}
$(z+z')^{q^i}=z^{q^i}+z'^{q^i}$, and if $u\in \mathbb{Z}_q$ and $z\in \mathbb{F}$, then $(uz)^{q^i}=u^{q^i}z^{q^i}=uz^{q^i}$.  
\end{center}
This justifies the name of the elements of $\mathfrak{L}$. If $X:=\{z_0,\dots,z_{k-1}\}$ is a $\mathbb{Z}_q$-basis of $\mathbb{F}$, then $M_g$ will denote the matrix of $g$ in the basis $X$.

\begin{proposition}
Let $\mathfrak{L}$ be the set of $q$-linearized polynomials. Then,
\begin{enumerate}
\item[\rm (i)] $\mathfrak{L}$ is a ring, where the addition is the usual addition in $\mathbb{F}[y]$ and the product $\circ$ is the \textbf{composition of polynomials}, i.e., $zy^{q^i}\circ z'y^{q^j}:=zz'^{q^i}y^{q^{i+j}}$, with $z,z'\in \mathbb{F}$.
\item[\rm (ii)]The function
\begin{align*}
\Lambda:A & \to \mathfrak{L}\\
\sum_{i=0}^m g_ix^i & \mapsto \sum_{i=0}^m g_iy^{q^i}
\end{align*}
is a ring isomorphism.
\item[\rm (iii)]Let $X:=\{z_0,\dots,z_{k-1}\}$ be a $\mathbb{Z}_q$-basis of $\mathbb{F}$. The \textbf{Moore matrix} of $X$ is defined by
\begin{center}
${\rm M}(X):=\begin{bmatrix}
z_0 & z_1 & \cdots & z_{k-2} & z_{k-1}\\
z_{0}^q & z_1^q & \cdots & z_{k-2}^q & z_{k-1}^q\\
\vdots & \vdots & & \vdots & \vdots\\
z_0^{q^{k-1}} & z_1^{q^{k-1}} & \cdots & z_{k-2}^{q^{k-1}} & z_{k-1}^{q^{k-1}}
\end{bmatrix}$.
\end{center}
Then, ${\rm M}(X)\in GL_k(\mathbb{F})$.
\item[\rm (iv)]For $g:=\sum_{i=0}^{k-1} g_iy^{q^i}\in \mathfrak{L}$, the \textbf{Dickson matrix} of $g$ $($also known as $\boldsymbol{q}$-\textbf{circulant matrix}$)$ is defined by
\begin{center}
$D_g:=\begin{bmatrix}
g_0 & g_1 & \cdots & g_{k-2} & g_{k-1}\\
g_{k-1}^q & g_0^q & \cdots & g_{k-3}^q & g_{k-2}^q\\
\vdots & \vdots & & \vdots & \vdots\\
g_1^{q^{k-1}} & g_2^{q^{k-1}} & \cdots & g_{k-1}^{q^{k-1}} & g_0^{q^{k-1}}
\end{bmatrix}$.
\end{center}
Then, $D_g={\rm M}(X)M_g{\rm M}(X)^{-1}$. 
\item[\rm (v)]The following ring isomorphisms hold: 
\begin{center}
	$A/\langle x^k-1\rangle\cong \mathfrak{L}/\langle y^{q^k}-y\rangle\cong M_{k\times k}(\mathbb{Z}_q)$.
\end{center}
\end{enumerate}     
\end{proposition}
\begin{proof}
(i) The composition of polynomials is associative; the unity element is $1y^{q^0}=y$; the distributive law holds. Thus, $\mathfrak{L}$ is a ring.

(ii) It is clear that $\Lambda$ is additive; $\Lambda(g_ix^ig_jx^j)=\Lambda(g_i\sigma^i(g_j)x^{i+j})=\Lambda(g_ig_j^{q^i}x^{i+j})=g_ig_j^{q^i}y^{q^{i+j}}$ and $\Lambda(g_ix^i)\Lambda(g_jy^{q^j})=g_iy^{q^i}\circ g_jy^{q^j}=g_ig_j^{q^i}y^{q^{i+j}}$, so $\Lambda$ is multiplicative; $\Lambda(1)=y$. Thus, $\Lambda$ is a ring homomorphism. Since $\mathbb{F}[y]$ is $\mathbb{F}$-free and $\mathfrak{L}\subset \mathbb{F}[y]$, then $\Lambda$ is injective; it is clear that $\Lambda$ is surjective.

(iii) Observe that $M(X)={\rm Wr}_k(Z)$, with $Z:=(z_0,\dots,z_{k-1})$ and $\mathcal{D}:=\sigma=\phi$ (see Definition \ref{definition2.16}). Since $z_0,\dots,z_{k-1}$ are linearly independent over $\mathbb{Z}_q$, Corollary 4.14 of \cite{Leroy} says that ${\rm Wr}_k(Z)$ is invertible.  

(iv) We will follow the proof of Lemma 4.1 in \cite{Wu}. First observe that
\begin{center}
$
\begin{bmatrix}
g\\
g^q\\
\vdots\\
g^{q^{k-1}}
\end{bmatrix}=\begin{bmatrix}
g_0 & g_1 & \cdots & g_{k-2} & g_{k-1}\\
g_{k-1}^q & g_0^q & \cdots & g_{k-3}^q & g_{k-2}^q\\
\vdots & \vdots & & \vdots & \vdots\\
g_1^{q^{k-1}} & g_2^{q^{k-1}} & \cdots & g_{k-1}^{q^{k-1}} & g_0^{q^{k-1}}
\end{bmatrix}\begin{bmatrix}
y\\
y^q\\
\vdots\\
y^{q^{k-1}}
\end{bmatrix}=D_g\begin{bmatrix}
y\\
y^q\\
\vdots\\
y^{q^{k-1}}
\end{bmatrix}
$, 
\end{center}
so $[g(z_j)^{q^i}]=D_g[z_j^{q^i}]=D_gM(X)$. On the other hand, $[g(z_j)^{q^i}]=M(X)M_g$. Indeed, let $M_g:=[m_{ij}]$, then $g(z_i)=m_{i0}z_0+\cdots+m_{ik-1}z_{k-1}$, for $0\leq i\leq k-1$ (recall that we dispose the scalars by rows), but since $m_{ij}\in \mathbb{Z}_q$, then $g(z_i)^{q^j}=m_{i0}z_0^{q^j}+\cdots+m_{ik-1}z_{k-1}^{q^j}=z_0^{q^j}m_{i0}+\cdots+z_{k-1}^{q^j}m_{ik-1}$. Hence, $g(z_j)^{q^i}=z_0^{q^i}m_{0j}+\cdots+z_{k-1}^{q^i}m_{k-1j}$ (see \cite{Wu}, p. 85), i.e., $[g(z_j)^{q^i}]=M(X)M_g$. Therefore, $D_gM(X)=M(X)M_g$, i.e., $D_g={\rm M}(X)M_g{\rm M}(X)^{-1}$.

(v) Since the order of the Frobenius automorphism is $k$, then from Proposition \ref{proposition2.30}, $x^k-1\in Z(A)$ is a two-sided polynomial, whence, $\langle x^k-1\rangle$ is the two-sided ideal $A(x^k-1)=(x^k-1)A$. As $\Lambda(x^k-1)=y^{q^k}-y$ and $\Lambda$ is an isomorphism, then $A/\langle x^k-1\rangle\cong \mathfrak{L}/\langle y^{q^k}-y\rangle$.

For the second isomorphism of (v), we recall that $ M_{k\times k}(\mathbb{Z}_q)$ and ${\rm End}_{\mathbb{Z}_q}(\mathbb{F})$ are isomorphic rings, so we will prove that ${\rm End}_{\mathbb{Z}_q}(\mathbb{F})$ and $\mathfrak{L}/\langle y^{q^k}-y\rangle$ are isomorphic rings. For this, note first that ${\rm End}_{\mathbb{Z}_q}(\mathbb{F})$ is an $\mathbb{F}$-space with natural product 
$(z\cdot \theta)(u):=z\cdot \theta(u)$, for $\theta\in {\rm End}_{\mathbb{Z}_q}(\mathbb{F})$ and $z,u\in \mathbb{F}$. Moreover, observe that $\{\sigma^i\mid 0\leq i\leq k-1\}$ is an $\mathbb{F}$-basis of ${\rm End}_{\mathbb{Z}_q}(\mathbb{F})$. In fact, notice first that $\sigma\in {\rm End}_{\mathbb{Z}_q}(\mathbb{F})$ since $\sigma(z'z)=\sigma(z')\sigma(z)=z'\sigma(z)$, with $z'\in \mathbb{Z}_q$ and $z\in \mathbb{F}$; $\dim_{\mathbb{F}}({\rm End}_{\mathbb{Z}_q}(\mathbb{F}))=k$ since $\dim_{\mathbb{Z}_q}({\rm End}_{\mathbb{Z}_q}(\mathbb{F}))=k^2$ and $\dim_{\mathbb{Z}_q}(\mathbb{F})=k$; finally, $\{\sigma^i\mid 0\leq i\leq k-1\}$ is $\mathbb{F}$-linearly independent since if not so, there exist $g_0,\dots,g_{k-1}\in \mathbb{F}$, which are not all zero, such that $\sum_{i=0}^{k-1}g_i\sigma^i =0$, so for every $z\in \mathbb{F}$, $\sum_{i=0}^{k-1} g_iz^{q^i}=0$, i.e., the non zero polynomial $\sum_{i=0}^{k-1} g_iy^{q^i}\in \mathbb{F}[y]$ has $q^k$ different roots, a contradiction. 

Therefore, we define
\begin{align*}
	{\rm End}_{\mathbb{Z}_q}(\mathbb{F}) & \xrightarrow{\Theta} \mathfrak{L}/\langle y^{q^k}-y\rangle\\
	\sum_{i=0}^{k-1}g_i\sigma^i & \mapsto (\sum_{i=0}^{k-1} g_iy^{q^i})+\langle y^{q^k}-y\rangle. 
\end{align*}
Clearly, $\Theta$ is additive, $\Theta(i_\mathbb{F})=y+\langle y^{q^k}-y\rangle$, $\Theta$ is bijective. Finally, $\Theta$ is multiplicative since for $z\in \mathbb{F}$,  
\begin{center}
	$(g_i\sigma^i\circ g_j\sigma^j)(z)=g_i\sigma^i(g_jz^{q^j})=g_ig_j^{q^j}z^{q^{i+j}}$, i.e., $g_i\sigma^i\circ g_j\sigma^j=g_ig_j^{q^j}\sigma^{i+j}$. 
\end{center}
\end{proof}

\begin{remark}
Observe that $A$, $\mathfrak{L}$ and $M_{k\times k}(\mathbb{Z}_q)$ are $\mathbb{Z}_q$-algebras and the isomorphisms in the previous proposition are isomorphisms of $\mathbb{Z}_q$-algebras.
\end{remark}

\subsection{The case $A:=\mathbb{F}(t)[x;\sigma]$ and $f=x^n-1$}\label{subsection2.3}

In this subsection we consider skew cyclic codes in the particular case when $A:=\mathbb{F}(t)[x;\sigma]$, where $\mathbb{F}$ is a finite field. This special situation of noncommutative coding theory was studied in \cite{Gomez-Torrecillas}. Observe that in this situation the finite field $\mathbb{F}$ was replaced by the field of fractions $\mathbb{F}(t)$ which is not finite. Moreover, in \cite{Gomez-Torrecillas} it is additionally assumed that if the order of $\sigma$ is $n\geq 1$, then the modulus $f$ is taken as $f=x^n-1$. This special restriction is not assumed in the general theory that we studied in the previous subsections. In this subsection we will review some key topics investigated in \cite{Gomez-Torrecillas}. 

\subsubsection{Algebraic structure}

According to Definition \ref{definition2.1}, we will have $R$-linear codes as submodules of the left $R$-module $R^n$, where $R:=\mathbb{F}(t)$, i.e., a $\mathbb{F}(t)$-linear code $\mathfrak{C}$ is a $\mathbb{F}(t)$-subspace of the $\mathbb{F}(t)$-vector space $\mathbb{F}(t)^n$; a skew cyclic code $\mathcal{C}$ is a left $A$-submodule of $A/Af$, but since $A$ is a left (and right) euclidean domain, and hence a left (and right) principal ideal domain, $\mathcal{C}$ has the form $Ag/Af$, for some $g\in A$, with 
$Af\subseteq Ag$. Recall that $A/Af$ is a $\mathbb{F}(t)$-vector space with canonical $\mathbb{F}(t)$-basis $X:=\{\overline{1},\overline{x},\dots, \overline{x^n-1}\}$; the length of $\mathcal{C}$ is $n$ and the dimension of $\mathcal{C}$ (as $\mathbb{F}(t)$-vector subspace of $A/Af$) is $n-\deg(g)$ (see Definition \ref{definition2.6}). Finally, $\mathbb{F}(t)^\sigma:=\{z\in \mathbb{F}(t)\mid \sigma(z)=z\}$ and $Z(A)=\mathbb{F}(t)^\sigma[x^n]$, hence, the modulus $f$ is a two-sided polynomial, i.e., $Af=fA$ is a two-sided ideal of $A$ and $\mathcal{A}:=A/Af$ is a ring. Actually, $\mathcal{A}$ is a $\mathbb{F}(t)^\sigma$-algebra.

\begin{proposition}[\cite{Gomez-Torrecillas}, Theorem 1]
Let $\mathbb{F}$, $\sigma$, $A$ and $f$ be as before. Then, 
\begin{enumerate}
\item[\rm (i)]${\rm Aut}(\mathbb{F}(t))\cong {\rm PGL}_2(\mathbb{F})$ $($the projective linear group$)$.
\item[\rm (ii)]For every $\theta\in {\rm Aut}(\mathbb{F}(t))$, $\mid \theta\mid<\infty$. 
\item[\rm (iii)]$\mathcal{A}\cong M_{n\times n}(\mathbb{F}(t)^\sigma)$ $($isomorphism of $\mathbb{F}(t)^\sigma$-algebras
$)$.
\item[\rm (iv)]For every $0\leq k\leq n$, there exists a skew cyclic code $\mathcal{C}$ of dimension $k$.
\end{enumerate}
\end{proposition}
\begin{proof}
(i) Recall first that ${\rm PGL}_2(\mathbb{F}):={\rm GL}_2(\mathbb{F})/Z({\rm GL}_2(\mathbb{F}))$, where ${\rm GL}_2(\mathbb{F})$ is the full linear group of all $2\times 2$ invertible matrices over $\mathbb{F}$; moreover, for $\theta\in {\rm Aut}(\mathbb{F}(t))$, $\theta(z)=z$, for every $z\in \mathbb{F}$. So we define
\begin{center}
$\Omega:{\rm GL}_2(\mathbb{F})\to {\rm Aut}(\mathbb{F}(t))$, 

$G:=\begin{bmatrix}a & b\\ c & d \end{bmatrix}\mapsto \theta_G$,

$\theta_G: \mathbb{F}(t)\to \mathbb{F}(t)$

$t\mapsto  \frac{at+b}{ct+d}$.
\end{center} 
It is well-known that $\Omega$ is a surjective function (see \cite{Lezama5}, Corollary 2.6.14). Actually, $\Omega$ is a surjective group homomorphism, where the composition of automorphisms in ${\rm Aut}(\mathbb{F}(t))$ works from the left to the right; moreover, $G\in \ker(\Omega)$ if and only if $\frac{at+b}{ct+d}=t$, i.e., $b=c=0$ and $a=d$. Thus, $\ker(\Omega)=Z({\rm GL}_2(\mathbb{F}))$. 

(ii) This follows from (i) since ${\rm PGL}_2(\mathbb{F})$ is a finite group.

(iii) Since $f$ is linear over $Z(A)$, then $f$ is irreducible over $Z(A)$. This implies that $Af=fA$ is a maximal ideal of $A$. In fact, $Z(A)f$ is a maximal ideal of $Z(A)$ and $Z(A)/Z(A)f$ is a field; consider the the canonical injective ring homomorphism $\theta: Z(A)/Z(A)f\to A/Af$, $\widehat{z}\to \overline{z}$, with $z\in Z(A)$; let $J$ be a two-sided ideal of $A$ such that $Af\subseteq J$, then $J=Ag$, for some $g\in A$; the idea is to show that $J=A$ or $J=Af$. We have $\theta^{-1}(Ag/Af)=Z(A)/Z(A)f$ or $\theta^{-1}(Ag/Af)=\{\widehat{0}\}$. In the first case, $\widehat{1}\in \theta^{-1}(Ag/Af)$, so $\overline{1}=\overline{ag}$, for some $a\in A$, whence $ag-1=qf=qpg$, for some $q,p\in A$, thus $(a-qp)g=1$, so $J=Ag=A$. In the second case, $\theta^{-1}(\{\overline{g}\})=\{\widehat{0}\}$, so $\theta(\widehat{0})=\overline{0}=\overline{g}$, so $g\in Af$, i.e., $J=Ag=Af$.  

Therefore, $\mathcal{A}$ is a simple left (and right) artinian ring. Indeed, $\mathcal{A}$ is simple since $Af$ is maximal. Let $I_1\supseteq I_2\supseteq I_3\cdots$ be a chain of left ideals of $\mathcal{A}$, then $I_i=Ag_i/Af$, with $g_i\in A$, for $i\geq 1$, then $Ag_1\supseteq Ag_2\supseteq Ag_3\cdots$, but $\deg(g_i)\leq \deg(f)=n$, for every $i\geq 1$, so the chain is finite.   

By the Artin-Wedderburn theorem, $\mathcal{A}\cong M_{r\times r}(\mathbb{D})$, where $\mathbb{D}={\rm End_{\mathcal{A}}}(M)$ is a division ring, $M$ is a simple left $\mathcal{A}$-module and  $r:=\dim_\mathbb{D}(M)$. We can take $M:=A/A(x-1)$. In fact, we have the isomorphism of left $\mathcal{A}$-modules $M\cong (A/Af)/(A(x-1)/Af)$ and $A(x-1)/Af$ is a maximal $\mathcal{A}$-submodule of $A/Af$. 

$\mathbb{D}\cong \mathbb{F}(t)^{\sigma}$: In fact, given $\phi\in \mathbb{D}$, let $\overline{\overline{g_{\phi}}}:=\phi(\overline{\overline{1}})\in M$, with $g_{\phi}\in A$, and let $g_{\phi}=p_{\phi}(x-1)+\lambda_{\phi}$, with $p_{\phi}\in A$ and $\lambda_{\phi}\in \mathbb{F}(t)$, we will prove that $\lambda_{\phi}\in \mathbb{F}(t)^{\sigma}$. We define
\begin{center}
$\Lambda: \mathbb{D}\to \mathbb{F}(t)^{\sigma}$, $\phi \mapsto \lambda_{\phi}$. 
\end{center}
$\lambda_{\phi}\in \mathbb{F}(t)^{\sigma}$: $\phi(\overline{\overline{x-1}})=\overline{\overline{0}}=\overline{x-1}\phi(\overline{\overline{1}})$, hence $\overline{\overline{(x-1)p_{\phi}(x-1)+(x-1)\lambda_{\phi}}}=\overline{\overline{0}}$, i.e., $\overline{\overline{(x-1)\lambda_{\phi}}}=\overline{\overline{0}}$, so $\sigma(\lambda_{\phi})x-\lambda_{\phi}=q(x-1)$, for some $q\in A$, whence $q\in \mathbb{F}$ and $\sigma(\lambda_{\phi})=q=\lambda_{\phi}$.

It is clear that $\Lambda$ is additive and $\Lambda(i_M)=1$. $\Lambda$ is multiplicative: 
\begin{center}
$(\phi_1\circ \phi_2)(\overline{\overline{1}})=\phi_1(\overline{\overline{g_{\phi_2}}})=\phi_1(\overline{g_{\phi_2}}\overline{\overline{1}})=\overline{g_{\phi_2}}\phi_1(\overline{\overline{1}})=
\overline{\overline{g_{\phi_2}g_{\phi_1}}}$,

$g_{\phi_1}=p_{\phi_1}(x-1)+\lambda_{\phi_1}$, $g_{\phi_2}=p_{\phi_2}(x-1)+\lambda_{\phi_2}$,

$g_{\phi_1}g_{\phi_2}=p(x-1)+\lambda_{\phi_1}\lambda_{\phi_2}$, for some $p\in A$ (here we use that $\lambda_{\phi_2}\in \mathbb{F}(t)^{\sigma}$).
\end{center}

Since $\mathbb{D}$ is a division ring, $\Lambda$ is injective. $\Lambda$ is surjective: Let $\lambda\in \mathbb{F}(t)^\sigma$, we define $\phi(\overline{\overline{g}}):=\overline{\overline{g\lambda}}$, with $\overline{\overline{g}}\in M$; $\phi$ is well-defined since $\lambda\in \mathbb{F}(t)^\sigma$; it is clear that  $\phi\in \mathbb{D}$ and $\Lambda(\phi)=\lambda$.

Finally, $r=n$: It is well-known (see \cite{Lang}, Theorem 1.8, Chapter VI) that $\dim_{\mathbb{F}(t)^\sigma}(\mathbb{F}(t))=\mid \sigma\mid=n$, but since $\dim_{\mathbb{F}(t)}(\mathcal{A})=n$, then $\dim_{\mathbb{F}(t)^\sigma}(\mathcal{A})=n$, so $n^2=r^2$, i.e., $r=n$.      

(iv) Observe that in the particular case of skew cyclic codes that we study in this subsection, 
\begin{center}
	\textit{there exists a bijective correspondence between}
	
	\textit{the collection of left ideals of $\mathcal{A}$ and the collection of skew cyclic codes}.
\end{center}

So, the claimed follows from the ring isomorphism $\mathcal{A}\cong  M_{n\times n}(\mathbb{F}(t)^\sigma)$. In fact, for every $0\leq k\leq n$, in $M_{n\times n}(\mathbb{F}(t)^\sigma)$ there exists a left ideal $I_k$ of dimension $nk$ as $\mathbb{F}(t)^\sigma$-vector space, namely, the matrices with the last $n-k$ columns are null. Thus, $I_k\cong (\mathbb{F}(t)^\sigma)^{nk}\cong \mathbb{F}(t)^k$, i.e., $I_k$ is a left ideal of $M_{n\times n}(\mathbb{F}(t)^\sigma)$ with $\mathbb{F}(t)$-dimension equals $k$. 
\end{proof}

Since the left $A$-submodules of  $\mathcal{A}$ coincide with its left ideals, we have the following consequence of (iii) of the previous proposition.

\begin{corollary}
If $\mathcal{C}$ is a skew cyclic code, then there exists and idempotent $\overline{e}\in \mathcal{A}$ such that $\mathcal{C}=\mathcal{A}\overline{e}$, and conversely. 
\end{corollary}
\begin{proof}
According to (iv) of the previous proposition, $\mathcal{A}$ is a semisimple ring, and since $\mathcal{C}$ is a left ideal of $\mathcal{A}$, then $\mathcal{C}$ is a direct summand of $\mathcal{A}$, so we get the claimed.
\end{proof}

\subsubsection{Dual codes}

Let $\mathcal{C}$ be a skew cyclic code. In this section we will prove that $\mathcal{C}^{\perp}$ is a skew cyclic code. For this we introduce the following notation.

\begin{definition}
For $\overline{g}\in \mathcal{A}$, let $.\overline{g}$ be the homomorphism of left $\mathcal{A}$-modules defined by $.\overline{g}(\overline{h}):=\overline{h}\overline{g}=\overline{hg}$, for $\overline{h}\in \mathcal{A}$. The matrix of $.\overline{g}$ in the canonical basis $X$ is denoted by $M(\overline{g})$. 
\end{definition}

\begin{proposition}[\cite{Gomez-Torrecillas}, Section 3]
Let $g:=g_0+g_1x+\cdots+g_{n-1}x^{n-1}\in A$. Then, 
\begin{enumerate}
\item[\rm (i)]
$M(\overline{g})=\begin{bmatrix}g_0 & g_1 & \cdots & g_{n-1}\\
\sigma(g_{n-1}) & \sigma(g_0) & \cdots & \sigma(g_{n-2})\\
\vdots & \vdots & \vdots & \vdots\\
\sigma^{n-1}(g_1) & \sigma^{n-1}(g_2) & \cdots & \sigma^{n-1}(g_0)
\end{bmatrix}\in M_{n\times n}(\mathbb{F}(t))$.
\item[\rm (ii)]For every $Z\in \mathbb{F}(t)^n$, $\mathfrak{p}_f(ZM(\overline{g}))=\mathfrak{p}_f(Z)\overline{g}$.
\item[\rm (iii)]The function
\begin{center}
$M:\mathcal{A}\to M_{n\times n}(\mathbb{F}(t))$, $\overline{g}\mapsto M(\overline{g})$
\end{center}
is an injective homomorphism of $\mathbb{F}(t)^\sigma$-algebras.
\item[\rm (iv)]${\rm Im}(.M(\overline{g}))=\mathfrak{p}_f^{-1}(\mathcal{A}\overline{g})$.
\item[\rm (v)]The function
\begin{center}
$\Theta: \mathcal{A}\to \mathcal{A}$, $\Theta(\overline{g}):=\overline{\sigma^{n}(g_0)x^n+\sigma^{n-1}(g_1)x^{n-1}+\cdots+\sigma(g_{n-1})x}$
\end{center}
is a ring anti-isomorphism such that $\Theta^2=i_{\mathcal{A}}$.
\item[\rm (vi)]$M(\Theta(\overline{g}))=M(\overline{g})^T$.
\item[\rm (vii)]If $\mathcal{C}$ is a skew cyclic code, then $\mathcal{C}^{\perp}$ is a skew cyclic code. More exactly, if $\mathcal{C}=\mathcal{A}\overline{e}$, where $\overline{e}$ is an idempotent of $\mathcal{A}$, then $\mathfrak{p}_f(\mathcal{C}^{\perp})=\mathcal{A}(\Theta(\overline{1}-\overline{e}))$.
\end{enumerate}
\end{proposition} 
\begin{proof}
(i) The result is trivial computing $.\overline{g}(\overline{x^i})$, for $0\leq i\leq n-1$. For example, $.\overline{g}(\overline{1})=\overline{g}=\overline{g_0+g_1x+\cdots+g_{n-1}x^{n-1}}=g_0\overline{1}+g_1\overline{x}+\cdots+g_{n-1}\overline{x^{n-1}}$ and $.\overline{g}(\overline{x^{n-1}})=\overline{x^{n-1}(g_0+g_1x+\cdots+g_{n-1}x^{n-1})}=\sigma^{n-1}(g_0)\overline{x^{n-1}}+\sigma^{n-1}(g_1)\overline{1}+\sigma^{n-1}(g_2)\overline{x}+\cdots+\sigma^{n-1}(g_n)\overline{x^{n-2}}$.

(ii) Since $\mathfrak{p}_f$ is $\mathbb{F}(t)$-linear, it is enough to show the claimed for $Z=e_i$, $1\leq i\leq n$, where $e_i$ is the $i$-th canonical vector of $\mathbb{F}(t)^n$. $\mathfrak{p}_f(e_iM(\overline{g}))=\mathfrak{p}_f(M(\overline{g})_i)$, where $(M(\overline{g})_i$ is the $i$-th row of $M(\overline{g})$, but $\mathfrak{p}_f(e_i)\overline{g}=\overline{x^i}\overline{g}=\overline{x^ig}=\mathfrak{p}_f(M(\overline{g})_i)$.   

(iii) It is clear that $M$ is additive, injective, $M(\overline{1})=I_n$ (the identical matrix), $M$ is $\mathbb{F}(t)^\sigma$-linear. Finally, $M$ is multiplicative since, from (ii), for every $1\leq i\leq n$,
\begin{center}
$\mathfrak{p}_f(e_iM(\overline{g})M(\overline{g'}))=\mathfrak{p}_f(e_iM(\overline{g}))\overline{g'}=\mathfrak{p}_f(e_i)\overline{g}\,\overline{g'}=\mathfrak{p}_f(e_i)\overline{gg'}$,

$\mathfrak{p}_f(e_iM(\overline{gg'}))=\mathfrak{p}_f(e_i)\overline{gg'}$,
\end{center}
but $\mathfrak{p}_f$ is injective, so $e_iM(\overline{g})M(\overline{g'})=e_iM(\overline{gg'})$, hence $M(\overline{g})M(\overline{g'})=M(\overline{gg'})=M(\overline{g}\,\overline{g'})$.     

(iv) ${\rm Im}(.M(\overline{g}))$ is the $\mathbb{F}(t)$-space generated by the rows of $M(\overline{g})$, i.e., the $\mathbb{F}(t)$-space generated by the vectors $e_1M(\overline{g}),\dots,e_nM(\overline{g})$, but $e_i=\mathfrak{p}_f^{-1}(\overline{x^i})$, hence, from (ii), ${\rm Im}(.M(\overline{g}))$ is the $\mathbb{F}(t)$-space generated by $\{e_iM(\overline{g})\}_{i=1}^n=\{\mathfrak{p}_f^{-1}(\overline{x^i}\overline{g})\}_{i=1}^n$, i.e., ${\rm Im}(.M(\overline{g}))=\mathfrak{p}_f^{-1}(\mathcal{A}\overline{g})$. 

(v) Is is clear that $\Theta$ is additive, and since $\sigma^n=i_{\mathbb{F}(t)}$, then $\Theta(\overline{1})=\overline{1}$. $\Theta$ is anti-multiplicative:  
\begin{center}
$\Theta(\overline{zx^i}\,\overline{z'x^j})=\Theta(\overline{z\sigma^i(z')x^{i+j}})=\overline{\sigma^{n-(i+j)}(z)\sigma^{n-j}(z')x^{n-(i+j)}}$,

$\Theta(\overline{z'x^j})\Theta(\overline{zx^i})=\overline{\sigma^{n-j}(z')x^{n-j}}\,\overline{\sigma^{n-i}(z)x^{n-i}}=\overline{\sigma^{n-j}(z')x^{n-j}\sigma^{n-i}(z)x^{n-i}}=\overline{\sigma^{n-j}(z')\sigma^{n-(i+j)}(z)x^{n-(i+j)}}$. 
\end{center}
Now we prove that $\Theta^2=i_{\mathcal{A}}$:
\begin{center}
	$\Theta^2(\overline{zx^i})=\Theta(\overline{\sigma^{n-i}(z)x^{n-i}})=\overline{\sigma^{n-(n-i)}(\sigma^{n-i}(z))x^{n-(n-i)}}=\overline{zx^i}$, for $z\in \mathbb{F}(t)$ and $0\leq i\leq n-1$.
\end{center} 

(vi) This follows from (v) and (i).

(vii) By (ii), we have the following commutative diagram of $\mathbb{F}(t)$-vector spaces: 
\begin{equation*}
\begin{CD}
\mathbb{F}(t)^n @>{.M(\overline{e})}>> \mathbb{F}(t)^n @>{.M(\overline{1}-\overline{e})}>> \mathbb{F}(t)^n\\
@VV{\mathfrak{p}_f}V @VV{\mathfrak{p}_f}V @VV{\mathfrak{p}_f}V\\
\mathcal{A} @>{.\overline{e}}>> \mathcal{A} @>{.(\overline{1}-\overline{e})}>> \mathcal{A}
\end{CD}
\end{equation*}
Observe that the second row is exact: $\overline{g}\in {\rm ker}(.(\overline{1}-\overline{e})$ if and only if $\overline{g}(\overline{1}-\overline{e})=\overline{0}$ if and only if $\overline{g}=\overline{g}\,\overline{e}$ if and only if $\overline{g}\in {\rm Im}(.\overline{e})$. Since $\mathfrak{p}_f$ is an isomorphism, then the first row is exact. Hence, if $\mathfrak{C}:=\mathfrak{p}_f^{-1}(\mathcal{C})$, then $\mathcal{C}=\mathfrak{p}_f(\mathfrak{C})=\mathcal{A}\overline{e}={\rm Im}(.\overline{e})$, whence $\mathfrak{C}=\mathfrak{p}_f^{-1}({\rm Im}(.\overline{e}))={\rm Im}(.M(\overline{e}))=\ker(.M(\overline{1}-\overline{e}))$. Moreover, from (v), $\overline{g}\in \mathcal{A}\Theta(\overline{1}-\overline{e})$ if and only if $\overline{g}=\overline{g}\Theta(\overline{1}-\overline{e})$ if and only if $\overline{g}\Theta(\overline{e})=\overline{0}$, so from the exactness of a similar commutative diagram we get that ${\rm Im}(.M(\Theta(\overline{1}-\overline{e})))=\ker(.M(\Theta(\overline{e})))$. But from Definition \ref{definition2.3}, $\mathfrak{C}^{\perp}=\ker(.M(\overline{e})^T)$, from (vi) we conclude that $\mathfrak{C}^{\perp}=\ker(.M(\Theta(\overline{e})))={\rm Im}(.M(\Theta(\overline{1}-\overline{e})))$. Hence, from (iv),
\begin{center}
 $\mathfrak{C}^{\perp}=\mathfrak{p}_f^{-1}(\mathcal{A}\Theta(\overline{1}-\overline{e}))$, i.e., $\mathfrak{p}_f(\mathcal{C}^{\perp})=\mathfrak{p}_f(\mathfrak{C}^{\perp})=\mathcal{A}\Theta(\overline{1}-\overline{e})$.
\end{center} 
\end{proof}

\subsubsection{Computation of generating idempotents}

\begin{proposition}[\cite{Gomez-Torrecillas}, Proposition 9]
Let $\mathcal{C}=\mathcal{A}\overline{e}$ be a skew cyclic code, with $\overline{e}\in \mathcal{A}$ be an idempotent. Then, for $g:={\rm grcd}(e,x^{n}-1)$, $\overline{g}$ is a minimal generator of $\mathcal{C}$. 
\end{proposition}
\begin{proof}
Since $g\mid_r e$, then $e=pg$, for some $p\in A$, so $\overline{e}\in \mathcal{A}\overline{g}$; moreover, $g=ue+v(x^{n}-1)$, for some $u,v\in A$, whence $\overline{g}=\overline{u}\overline{e}=\overline{u}\,\overline{e}$, so $\overline{g}\in \mathcal{A}\overline{e}$. Hence, $\mathcal{C}=\mathcal{A}\overline{e}=\mathcal{A}\overline{g}$. Let $h\in A$ such that $\mathcal{C}=\mathcal{A}\overline{h}$, then $\deg(g)\leq \deg(h)$. Indeed, $\overline{h}=\overline{p}\,\overline{g}$, for some $p\in A$, hence $h-pg=q(x^n-1)$, for some $q\in A$, but $x^n-1=tg$, for some $t\in A$, so $h=pg+qtg=(p+gt)g$, thus $\deg(h)\geq \deg(g)$. 
\end{proof}

\begin{proposition}[\cite{Gomez-Torrecillas}, Proposition 10]
Let $g,h\in A$ such that ${\rm lclm}(g,h)=x^{n}-1$ and $\deg(g)+\deg(h)=n$. Then,
\begin{enumerate}
\item[\rm (i)]${\rm grcd}(g,h)=1$.
\item[\rm (ii)]Let $u,v\in A$ such that $1=ug+vh$. Then, $e:=ug$ is such that $\overline{e}$ is a generating idempotent of $\mathcal{A}\overline{g}$.
\end{enumerate} 
\end{proposition}
\begin{proof}
(i) We have $\deg({\rm lclm}(g,h))=n=\deg(g)+\deg(h)-\deg({\rm grcd}(g,h))$ (see (\ref{equation2.4a}) and Remark \ref{remark2.10}), whence $\deg({\rm grcd}(g,h))=0$, so ${\rm grcd}(g,h)=1$.

(ii) From $1=ug+vh$ we get $\overline{1}=\overline{ug}+\overline{vh}=\overline{u}\,\overline{g}+\overline{v}\overline{h}$, i.e., $\mathcal{A}=\mathcal{A}\overline{g}+\mathcal{A}\overline{h}$. Actually, $\mathcal{A}=\mathcal{A}\overline{g}\oplus\mathcal{A}\overline{h}$. In fact, $g\mid_r(x^n-1)$ and $h\mid_r(x^n-1)$, so $Af\subseteq Ag$ and $Af\subseteq Ah$, hence   $\dim_{\mathbb{F}(t)}(\mathcal{A}\overline{g})+\dim_{\mathbb{F}(t)}(\mathcal{A}\overline{h})=\dim_{\mathbb{F}(t)}(Ag/Af)+\dim_{\mathbb{F}(t)}(Ah/Af)=n-\deg(g)+n-\deg(h)=n=\dim_{\mathbb{F}(t)}(\mathcal{A})$, so $\dim_{\mathbb{F}(t)}(\mathcal{A}\overline{g}\cap \mathcal{A}\overline{h})=0$, i.e., $\mathcal{A}\overline{g}\cap \mathcal{A}\overline{h}=0$. 
Therefore, $\overline{e}=\overline{ug}$ is a generator idempotent of the left ideal $\mathcal{A}\overline{g}$. 	
\end{proof}

\section{Skew $PBW$ extensions}

In this section we recall some basic facts about the class of noncommutative rings of polynomial type known as skew $PBW$ extensions.  

\begin{definition}[\cite{LezamaGallego},\cite{Lezama-sigmaPBW}]\label{gpbwextension}
	Let $R$ and $A$ be rings. We say that $A$ is a \textit{\textbf{skew $PBW$
			extension of $R$}} $($also called a $\sigma-PBW$ extension of
	$R$$)$ if the following conditions hold:
	\begin{enumerate}
		\item[\rm (i)]$R\subseteq A$.
		\item[\rm (ii)]There exist finitely many elements $x_1,\dots ,x_n\in A$ such $A$ is an $R$-free left module with basis
		\begin{center}
			${\rm Mon}(A):= \{x^{\alpha}=x_1^{\alpha_1}\cdots
			x_n^{\alpha_n}\mid \alpha=(\alpha_1,\dots ,\alpha_n)\in
			\mathbb{N}^n\}$, with $\mathbb{N}:=\{0,1,2,\dots\}$.
		\end{center}
		In this case we say that $A$ is a \textbf{ring of left polynomial type} over $R$ with respect to
		$\{x_1,\dots,x_n\}$. The set ${\rm Mon}(A)$ is called the set of \textbf{standard monomials} of
		$A$.
		\item[\rm (iii)]For every $1\leq i\leq n$ and $r\in R-\{0\}$ there exists $c_{i,r}\in R-\{0\}$ such that
		\begin{equation}\label{sigmadefinicion1}
		x_ir-c_{i,r}x_i\in R.
		\end{equation}
		\item[\rm (iv)]For every $1\leq i,j\leq n$ there exists $c_{i,j}\in R-\{0\}$ such that
		\begin{equation}\label{sigmadefinicion2}
		x_jx_i-c_{i,j}x_ix_j\in R+Rx_1+\cdots +Rx_n.
		\end{equation}
		Under these conditions we will write $A:=\sigma(R)\langle
		x_1,\dots ,x_n\rangle$.
	\end{enumerate}
\end{definition}
Associated to a skew $PBW$ extension $A=\sigma(R)\langle x_1,\dots
,x_n\rangle$ there are $n$ injective endomorphisms
$\sigma_1,\dots,\sigma_n$ of $R$ and $\sigma_i$-derivations, as
the following proposition shows.

\begin{proposition}[\cite{LezamaGallego}, Proposition 3]\label{sigmadefinition}
	Let $A$ be a skew $PBW$ extension of $R$. Then, for every $1\leq
	i\leq n$, there exist an injective ring endomorphism
	$\sigma_i:R\rightarrow R$ and a $\sigma_i$-derivation
	$\delta_i:R\rightarrow R$ such that
	\begin{center}
		$x_ir=\sigma_i(r)x_i+\delta_i(r)$,
	\end{center}
	for each $r\in R$.
\end{proposition}

Two remarkable particular cases of skew $PBW$ extensions are recalled next. 

\begin{definition}[\cite{Lezama-sigmaPBW}, Chapter 1]\label{sigmapbwderivationtype}
	Let $A$ be a skew $PBW$ extension.
	\begin{enumerate}
		\item[\rm (a)]
		$A$ is \textbf{quasi-commutative}\index{quasi-commutative} if conditions {\rm(}iii{\rm)} and {\rm(}iv{\rm)} in Definition
		\ref{gpbwextension} are replaced by
		\begin{enumerate}
			\item[\rm ($iii'$)]For every $1\leq i\leq n$ and $r\in R-\{0\}$ there exists a $c_{i,r}\in R-\{0\}$ such that
			\begin{equation}
			x_ir=c_{i,r}x_i.
			\end{equation}
			\item[\rm ($iv'$)]For every $1\leq i,j\leq n$ there exists $c_{i,j}\in R-\{0\}$ such that
			\begin{equation}
			x_jx_i=c_{i,j}x_ix_j.
			\end{equation}
		\end{enumerate}
		\item[\rm (b)]$A$ is \textbf{bijective}\index{bijective} if $\sigma_i$ is bijective for
		every $1\leq i\leq n$ and $c_{i,j}$ is invertible for any $1\leq i,j\leq n$.
	\end{enumerate}
\end{definition}

If $A=\sigma(R)\langle x_1,\dots,x_n\rangle$ is a skew $PBW$
extension of the ring $R$, then, as was observed in Proposition
\ref{sigmadefinition}, $A$ induces injective endomorphisms
$\sigma_k:R\to R$ and $\sigma_k$-derivations $\delta_k:R\to R$,
$1\leq k\leq n$. Moreover, from the Definition
\ref{gpbwextension}, there exists a unique finite set of constants
$c_{ij}, d_{ij}, a_{ij}^{(k)}\in R$, $c_{ij}\neq 0$, such that
\begin{equation}\label{equation1.2.1}
x_jx_i=c_{ij}x_ix_j+a_{ij}^{(1)}x_1+\cdots+a_{ij}^{(n)}x_n+d_{ij},
\ \text{for every}\  1\leq i<j\leq n.
\end{equation}
If $A$ is quasi-commutative, then $\delta_k=0$ for every $1\leq k\leq n$ and $p_{\alpha,r},p_{\alpha, \beta}=0$ in Proposition \ref{coefficientes}.

Many important algebras and rings coming from mathematical physics and non-commutative algebraic geometry
are particular examples of skew $PBW$ extensions: \textbf{Habitual ring of
	polynomials in several variables}, Weyl algebras, enveloping
algebras of finite dimensional Lie algebras, algebra of
$q$-differential operators, many important types of Ore algebras, in particular, the single Ore extensions of Section \ref{section2}, 
algebras of diffusion type, additive and multiplicative analogues
of the Weyl algebra, dispin algebra $\mathcal{U}(osp(1,2))$,
quantum algebra $\mathcal{U}'(so(3,K))$, Woronowicz algebra
$\mathcal{W}_{\nu}(\mathfrak{sl}(2,K))$, Manin algebra
$\mathcal{O}_q(M_2(K))$, coordinate algebra of the quantum group
$SL_q(2)$, $q$-Heisenberg algebra \textbf{H}$_n(q)$, Hayashi
algebra $W_q(J)$, differential operators on a quantum space
$D_{\textbf{q}}(S_{\textbf{q}})$, Witten's deformation of
$\mathcal{U}(\mathfrak{sl}(2,K))$, multiparameter Weyl algebra
$A_n^{Q,\Gamma}(K)$, quantum symplectic space
$\mathcal{O}_q(\mathfrak{sp}(K^{2n}))$, some quadratic algebras in
3 variables, some 3-dimensional skew polynomial algebras,
particular types of Sklyanin algebras, homogenized enveloping
algebra $\mathcal{A}(\mathcal{G})$, Sridharan enveloping algebra
of 3-dimensional Lie algebra $\mathcal{G}$, among many others. For
a precise definition of any of these rings and algebras see 
\cite{lezamareyes1} and \cite{Lezama-sigmaPBW}. The skew $PBW$ has been intensively studied in the last years (see \cite{Lezama-sigmaPBW}).

Next we will fix some notation and a monomial order in $A$ (see \cite{Lezama-sigmaPBW}, Chapter 1). 

\begin{definition}\label{1.1.6}
	Let $A$ be a skew $PBW$ extension of $R$ with endomorphisms $\sigma_i$ as in Proposition
	\ref{sigmadefinition}, $1\leq i\leq n$.
	\begin{enumerate}
		\item[\rm (i)]For $\alpha=(\alpha_1,\dots,\alpha_n)\in \mathbb{N}^n$,
		$\boldsymbol{\sigma^{\alpha}}:=\sigma_1^{\alpha_1}\cdots \sigma_n^{\alpha_n}$,
		$\boldsymbol{|\alpha|}:=\alpha_1+\cdots+\alpha_n$. If $\beta=(\beta_1,\dots,\beta_n)\in \mathbb{N}^n$, then
		$\boldsymbol{\alpha+\beta}:=(\alpha_1+\beta_1,\dots,\alpha_n+\beta_n)$.
		\item[\rm (ii)]For $X=x^{\alpha}\in \mathrm{Mon}(A)$,
		$\boldsymbol{\exp(X)}:=\alpha$ and $\boldsymbol{\deg(X)}:=|\alpha|$.
		\item[\rm (iii)]Let $0\neq f\in A$. If $t(f)$ is the finite
		set of terms that conform $f$, i.e., if $f=c_1X_1+\cdots +c_tX_t$, with $X_i\in \mathrm{Mon}(A)$ and $c_i\in
		R-\{0\}$, then $\boldsymbol{t(f)}:=\{c_1X_1,\dots,c_tX_t\}$.
		\item[\rm (iv)]Let $f$ be as in {\rm(iii)}, then $\boldsymbol{\deg(f)}:=\max\{\deg(X_i)\}_{i=1}^t.$
	\end{enumerate}
\end{definition}

In $\mathrm{Mon}(A)$ we define
\begin{center}
	$x^{\alpha}\succeq x^{\beta}\Longleftrightarrow
	\begin{cases}
	x^{\alpha}=x^{\beta}\\
	\text{or} & \\
	x^{\alpha}\neq x^{\beta}\, \text{but} \, |\alpha|> |\beta| & \\
	\text{or} & \\
	x^{\alpha}\neq x^{\beta},|\alpha|=|\beta|\, \text{but $\exists$ $i$ with} &
	\alpha_1=\beta_1,\dots,\alpha_{i-1}=\beta_{i-1},\alpha_i>\beta_i.
	\end{cases}$
\end{center}
It is clear that this is a total order on $\mathrm{Mon}(A)$, called \textit{\textbf{deglex}} order. If
$x^{\alpha}\succeq x^{\beta}$ but $x^{\alpha}\neq x^{\beta}$, we write $x^{\alpha}\succ x^{\beta}$.
Each element $f\in A-\{0\}$ can be represented in a unique way as $f=c_1x^{\alpha_1}+\cdots
+c_tx^{\alpha_t}$, with $c_i\in R-\{0\}$, $1\leq i\leq t$, and $x^{\alpha_1}\succ \cdots \succ
x^{\alpha_t}$. We say that $x^{\alpha_1}$ is the \textit{\textbf{leading monomial}} of $f$ and we write
$lm(f):=x^{\alpha_1}$; $c_1$ is the \textit{\textbf{leading coefficient}} of $f$, $lc(f):=c_1$, and
$c_1x^{\alpha_1}$ is the \textit{\textbf{leading term}} of $f$ denoted by $lt(f):=c_1x^{\alpha_1}$. We say that $f$ is \textit{\textbf{monic}} if $lc(f):=1$. If $f=0$,
we define $lm(0):=0$, $lc(0):=0$, $lt(0):=0$, and we set $X\succ 0$ for any $X\in \mathrm{Mon}(A)$. We observe that
\begin{center}
	$x^{\alpha}\succ x^{\beta}\Rightarrow lm(x^{\gamma}x^{\alpha}x^{\lambda})\succ
	lm(x^{\gamma}x^{\beta}x^{\lambda})$, for every $x^{\gamma},x^{\lambda}\in \mathrm{Mon}(A)$.
\end{center}

The next proposition complements Definition \ref{gpbwextension}.

\begin{proposition}[\cite{LezamaGallego},\cite{Lezama-sigmaPBW}]\label{coefficientes}
	Let $A$ be a ring of a left polynomial type over $R$ w.r.t.\ $\{x_1,\dots,x_n\}$. $A$ is a skew
	$PBW$ extension of $R$ if and only if the following conditions hold:
	\begin{enumerate}
		\item[\rm (a)]For every $x^{\alpha}\in \mathrm{Mon}(A)$ and every $0\neq
		r\in R$ there exist unique elements $r_{\alpha}:=\sigma^{\alpha}(r)\in R-\{0\}$ and $p_{\alpha
			,r}\in A$ such that
		\begin{equation}\label{611}
		x^{\alpha}r=r_{\alpha}x^{\alpha}+p_{\alpha , r},
		\end{equation}
		where $p_{\alpha ,r}=0$ or $\deg(p_{\alpha ,r})<|\alpha|$ if $p_{\alpha , r}\neq 0$. Moreover, if
		$r$ is left invertible, then $r_\alpha$ is left invertible.
		
		\item[\rm (b)]For every $x^{\alpha},x^{\beta}\in \mathrm{Mon}(A)$ there
		exist unique elements $c_{\alpha,\beta}\in R$ and $p_{\alpha,\beta}\in A$ such that
		\begin{equation}\label{612}
		x^{\alpha}x^{\beta}=c_{\alpha,\beta}x^{\alpha+\beta}+p_{\alpha,\beta},
		\end{equation}
		where $c_{\alpha,\beta}$ is left invertible, $p_{\alpha,\beta}=0$ or
		$\deg(p_{\alpha,\beta})<|\alpha+\beta|$ if $p_{\alpha,\beta}\neq 0$.
	\end{enumerate}
\end{proposition}

We conclude this subsection recalling some of the main ingredients of the Gröbner theory of skew $PBW$ extensions, namely, the Division Algorithm and the notion of Gröbner basis of a left ideal of $A$. For all details see \cite{Lezama-sigmaPBW}, Chapter 13. For the condition (ii) in Definition \ref{reductionsigmapbw} below, some natural computational conditions on $R$ will be assumed.
\begin{definition}\label{LGSring}
	A ring $R$ is \textbf{left Gr\"obner soluble}\index{left Grobner soluble@left Gr\"obner soluble} {\rm(}$LGS${\rm)}\index{LGS ring@$LGS$ ring} if the
	following conditions hold:
	\begin{enumerate}
		\item[\rm (i)]$R$ is left noetherian.
		\item[\rm (ii)]Given $a,r_1,\dots,r_m\in R$ there exists an
		algorithm which decides whether $a$ is in the left ideal
		$Rr_1+\cdots+Rr_m$, and if so, finds $b_1,\dots,b_m\in R$ such that
		$a=b_1r_1+\cdots+b_mr_m$.
		\item[\rm (iii)]Given $r_1,\dots,r_m\in R$ there exists an
		algorithm which finds a finite set of generators of the left
		$R$-module
		\begin{center}
			$\mathrm{Syz}_R[r_1\ \cdots \ r_m]:=\{(b_1,\dots,b_m)\in
			R^m\mid b_1r_1+\cdots+b_mr_m=0\}$.
		\end{center}
	\end{enumerate}
\end{definition}

\begin{definition}
	Let $x^{\alpha},x^{\beta}\in \mathrm{Mon}(A)$. We say that $x^{\alpha}$
	\textbf{divides}\index{divides} $x^{\beta}$, denoted by $x^{\alpha}\mid x^{\beta}$, if there exists a unique $x^{\theta}\in \mathrm{Mon}(A)$ such that
	$x^{\beta}=lm(x^{\theta}x^{\alpha})=x^{\theta+\alpha}$ and hence
	$\beta=\theta+\alpha$.
\end{definition}

\begin{definition}\label{reductionsigmapbw}
	Let $F$ be a finite set of nonzero elements of $A$, and let
	$f,h\in A$. We say that $f$ \textbf{reduces to $h$ by $F$ in one step},\index{reduces in one step}
	denoted $f\xrightarrow{\,\, F\,\, } h$, if there exist elements
	$f_1,\dots,f_t\in F$ and $r_1,\dots,r_t\in R$ such that
	\begin{enumerate}
		\item[\rm (i)]$lm(f_i)\mid lm(f)$, $1\leq i\leq t$, i.e., there exists an
		$x^{\alpha_i}\in \mathrm{Mon}(A)$ such that
		$lm(f)=lm(x^{\alpha_i}lm(f_i))$, i.e.,
		$\alpha_i+\exp(lm(f_i))=\exp(lm(f))$.
		\item[\rm
		(ii)]$lc(f)=r_1\sigma^{\alpha_1}(lc(f_1))c_{\alpha_1,f_1}+\cdots+r_t\sigma^{\alpha_t}(lc(f_t))c_{\alpha_t,f_t}$,
		where $c_{\alpha_i,f_i}$ are defined as in Theorem
		\ref{coefficientes}, i.e.,
		$c_{\alpha_i,f_i}:=c_{\alpha_i,\exp(lm(f_i))}$.
		\item[\rm (iii)]$h=f-\sum_{i=1}^tr_ix^{\alpha_i}f_i$.
	\end{enumerate}
	We say that $f$ \textbf{reduces}\index{reduces} to $h$ by $F$, denoted $f\xrightarrow{\,\,
		F\,\, }_{+}h$, if there exist $h_1,\dots ,h_{t-1}\in A$ such that
	\begin{center}
		$\begin{CD} f @>{F}>> h_1 @>{F}>> h_2 @>{F}>>\cdots @>{F}>>h_{t-1}
		@>{F}>>h.
		\end{CD}$
	\end{center}
	$f$ is \textbf{reduced}\index{reduced} {\rm(}also called \textbf{minimal}{\rm)}\index{minimal} w.r.t.\ $F$ if $f =
	0$ or there is no one step reduction  of $f$ by $F$, i.e., one of
	the conditions $(i)$ or $(ii)$
	fails. Otherwise, we will say that $f$ is \textbf{reducible}\index{reducible} w.r.t.\ $F$. If
	$f\xrightarrow{\,\, F\,\, }_{+} h$ and $h$ is reduced w.r.t.\ $F$,
	then we say that $h$ is a \textbf{remainder}\index{remainder} for $f$ w.r.t.\ $F$.
\end{definition}
By definition we will assume that $0\xrightarrow {F} 0$.

\begin{proposition}[Division algorithm]\label{algdivforPBW}
	Let $F=\{f_1,\dots ,f_t\}$ be a finite set of nonzero polynomials
	of $A$ and $f\in A$, then there exist
	polynomials $q_1,\dots ,q_t,h\in A$, with $h$ reduced w.r.t. $F$,
	such that $f\xrightarrow{\,\, F\,\, }_{+} h$ and
	\[
	f=q_1f_1+\cdots +q_tf_t+h,
	\]
	with
	\[
	lm(f)=\max\{lm(lm(q_1)lm(f_1)),\dots
	,lm(lm(q_t)lm(f_t)),lm(h)\}.
	\]
\end{proposition}
\begin{definition}\label{definition3.10}
	Let $I\neq 0$ be a left ideal of $A$ and
	let $G$ be a nonempty finite subset of nonzero polynomials of
	$I$. $G$ is a \textbf{Gröbner basis} \index{Grobner basis@Gr\"obner basis} for $I$ if each element
	$0\neq f\in I$ is reducible w.r.t.\ $G$.
\end{definition}
\begin{proposition}\label{153}
	Let $I\neq 0$ be a left ideal of $A$. Then,
	\begin{enumerate}
		\item[\rm(i)]If $G$ is a Gröbner basis for $I$, then $I=\langle G\}$ $($the left ideal of $A$ generated by $G$$)$.
		\item[\rm(ii)]Let $G$ be a Gröbner basis for $I$. If $f\in I$ and
		$f\xrightarrow{\,\, G\,\, }_{+} h$, with $h$ reduced, then $h=0$.
		\item[\rm(iii)]Let $G=\{g_1,\dots,g_t\}$ be a set of nonzero polynomials of $I$ with $lc(g_i)\in R^{*}$ for each $1\leq i\leq t$. Then,
		$G$ is a Gröbner basis of $I$ if and only if given $0\neq r\in I$
		there exists an $i$ such that $lm(g_i)$ divides $lm(r)$.
	\end{enumerate}
\end{proposition}

\begin{remark}\label{reamrk3.12}
	(i) We remark that the Gröbner theory of skew $PBW$ extensions and some of its important applications in homological algebra have been implemented in Maple in
	\cite{Fajardo2} and \cite{Fajardo3} (see also \cite{Lezama-sigmaPBW}). This implementation is based
	on the library \textbf{\texttt{SPBWE.lib}} specialized for working with bijective skew $PBW$
	extensions. The library has utilities to calculate Gröbner bases, and
	it includes some functions that compute the module of syzygies,
	free resolutions and left inverses of matrices, among other things. For the implementation was assumed that
	$A=\sigma(R)\langle x_1,\dots,x_n\rangle$ is a bijective skew $PBW$ extension of an $LGS$
	ring $R$ and $\mathrm{Mon}(A)$ is endowed with some monomial order $\succeq$.	
	
	(ii) \textbf{From now on in this paper we will assume that $\boldsymbol{A:=\sigma(R)\langle
			x_1,\dots ,x_n\rangle}$ is a bijective skew $\textbf{\textit{PBW}}$ extension of $\boldsymbol{R}$, where $\boldsymbol{R}$ is a left noetherian domain}. This implies that $A$ is a left noetherian domain (see \cite{Lezama-sigmaPBW}, Chapter 1). In the examples where we use the library \texttt{SPBWE.lib} we have assumed additionally that $R$ is $LGS$. This implies that $A$ is $LGS$ (see \cite{Lezama-sigmaPBW}, Chapter 15).   
\end{remark}

\section{Algebraic sets and ideals of points for skew $PBW$ extensions}

This last section represents the novelty of the present work and is dedicated to extend to skew $PBW$ extensions some results of the previous sections, 
more precisely, we will study the algebraic sets, the ideal of points and the relationship between them. Some properties of affine algebraic sets of commutative algebraic geometry  (see \cite{Fulton}, Chapter 1) will be extended in this section, as well as Proposition \ref{prop2.12A}. We will assume on $A$ the conditions in (ii) of Remark \ref{reamrk3.12}. 

As was pointed out in the introduction, the focus of the section is algebraic and not should be understood as a contribution to noncommutative coding theory.

\subsection{Roots of polynomials}

For $n\geq 1$, let $R^n$ be the left $R$-module of vectors over $R$ of $n$ components. Let $f\in A$ and $Z:=(z_1,\dots,z_n)\in R^n$. By Proposition \ref{algdivforPBW},there exist
polynomials $q_1,\dots ,q_t,h\in A$, with remainder $h$ reduced w.r.t.\ $F:=\{x_1-z_1,\dots,x_n-z_n\}$,
such that $f\xrightarrow{\,\, F\,\, }_{+} h$ and
\[
f=q_1(x_1-z_1)+\cdots +q_n(x_n-z_n)+h.
\]
In general, $h$ is not unique, and even worse, it could not belong to $R$, as the next example shows.

\begin{example}
Consider the Witten algebra (see \cite{Lezama-sigmaPBW}, Chapter 2) $A:=\sigma(\mathbb{Q})\langle x,y,z\rangle$ defined by
\begin{center}
	$zx = xz-x$, $zy = yz+2y$, $yx = 2xy$.
\end{center}
For $Z:=(1,-2-3)\in \mathbb{Q}^3$ and $f:=x^2y+xz+yz\in A$, with \textbf{\texttt{SPBWE.lib}} the Algorithm Division produces
\begin{center}
	$f=(\frac{1}{2}xy+\frac{1}{4}y)(x-1)+\frac{1}{4}(y+2)+0(z+3)+xz+yz-\frac{1}{2}$, 
\end{center}
i.e., $q_1=\frac{1}{2}xy+\frac{1}{4}y$, $q_2:=\frac{1}{4}$, $q_3=0$ and $h=xz+yz-\frac{1}{2}$. 

Even for quasi-commutative skew $PBW$ extensions the situation is similar. In fact, consider a $3$-multiparametric quantum space (see \cite{Lezama-sigmaPBW}, Chapter 4) $A:=\sigma(\mathbb{C})\langle x,y,z\rangle$ defined by
\begin{center}
	$yx = 2ixy$, $zx=3ixz$, $zy=-iyz$.
\end{center}
For $Z:=(i,2i,3i)\in \mathbb{C}^3$ and $f:=x^2y+yz^2+xz\in A$, with \textbf{\texttt{SPBWE.lib}} we found that
\begin{center}
	$f=(\frac{1}{2}ixy-\frac{1}{4}iy)(x-i)+\frac{1}{4}(y-2i)+0(z-3i)+yz^2+xz+\frac{1}{2}i$, 
\end{center}
i.e., $q_1=\frac{1}{2}ixy-\frac{1}{4}iy$, $q_2:=\frac{1}{4}$, $q_3=0$ and $h=yz^2+xz+\frac{1}{2}i$.  
\end{example}  

Thus, the evaluation of a polynomial $f\in A$ in a given $Z\in R^n$ as the remainder in the Division Algorithm is not a good idea. However, the following notion does not depend on the Division Algorithm. 
\begin{definition}
	Let $n\geq 1$, $f\in A$ and $Z:=(z_1,\dots,z_n)\in R^n$. $Z$ is a \textbf{root} of $f$ if and only if $f$ is in the two-sided ideal generated by $x_1-z_1,\dots, x_n-z_n$. This condition is denoted by $f(Z)=0$.
\end{definition}
Thus,
\begin{center}
$f(Z)=0$ if and only if $f\in \langle Z\rangle$,
\end{center}
where the two-sided ideal generated by $x_1-z_1,\dots, x_n-z_n$ is simply denoted by $\langle Z\rangle$, i.e., 
\begin{equation}
\langle Z\rangle:=\langle x_1-z_1, \dots, x_n-z_n\rangle.
\end{equation}	
\begin{definition}\label{definition4.14}
	Let $f\in A$. The \textbf{\textit{vanishing set}} of $f$, also called the \textbf{\textit{set of roots}} of $f$, is denoted by $V(f)$, and defined by  
\begin{equation}
V(f):=\{Z\in R^n\mid f(Z)=0\}.
\end{equation}
If $S\subseteq A$, then 
\begin{equation}
V(S):=\{Z\in R^n\mid f(Z)=0,\ \text{for every $f\in S$}\}.
\end{equation}
A subset $X\subseteq R^n$ is \textbf{algebraic} if either $X=R^n$ or there exists $g\neq 0\in A$ such that $X\subseteq V(g)$.
\end{definition} 

\subsection{Algebraic sets and ideals of points}

Some classical properties of affine algebraic sets of commutative algebraic geometry (see \cite{Fulton}, Chapter 1) will be extended in this subsection as well as Proposition \ref{prop2.12A}. 

\begin{theorem}\label{theorem4.4}
	\begin{enumerate}
		\item[\rm (i)]Let $f,g,h\in A$ and $Z:=(z_1,\dots,z_n)\in R^n$. 
		\begin{enumerate}
		\item[\rm (a)]If $f(Z)=0=g(Z)$, then $(f+g)(Z)=0$.
		\item[\rm (b)]$V(f)\subseteq V(gfh)$.
		\end{enumerate}
		\item[\rm (ii)]Let $I:=Ag$ be a left principal ideal of $A$. Then, $V(I)=V(g)$. The same is true for right and two-sided principal ideals of $A$.
		\item[\rm (iii)] 
			\begin{enumerate}
				\item[\rm (a)]$V(0)=R^n$. 
				\item[\rm (b)]$\emptyset$ is algebraic.
				\item[\rm (c)]If $S\subseteq T\subseteq A$, then $V(T)\subseteq V(S)$.
				\item[\rm (d)]If $S\subseteq A$, then $V(S)= V(AS)=V(SA)=V(ASA)$.
				\item[\rm (e)]$V(I)\cup V(J)\subseteq V(I\cap J)$, where $I,J$ are left $($right, two-sided$)$ ideals of $A$.
				\item[\rm (f)]$V(\sum_{k\in \mathcal{K}}I_k)=\bigcap_{k\in \mathcal{K}}V(I_k)$, where $I_k$ is a left left $($right, two-sided$)$ ideal of $A$.
				\item[\rm{(h)}]Let $Z:=(z_1,\dots,z_n)\in R^n$. Then, $\{Z\}\subseteq V(\langle Z\rangle)$. 
				\end{enumerate} 
		\item[\rm (iv)]Let $X\subseteq R^n$. Then,
		\begin{center}
		$I(X):=\{g\in A\mid g(Z)=0, \ \text{for every $Z\in X$}\}$
		\end{center}
		is a two-sided ideal of $A$, called the \textbf{ideal of points} of $X$. Some properties of $I(X)$ are:
		\begin{enumerate}
			\item[\rm (a)]$I(\emptyset)=A$. 
			\item[\rm{(b)}]For $X,Y\subseteq R^n$, $X\subseteq Y\Rightarrow I(Y)\subseteq I(X)$.
			\item[\rm{(c)}]If $I$ is a left $($right, two-sided$)$ ideal of $A$, then $I\subseteq I(V(I))$. 
			\item[\rm{(d)}]$X\subseteq V(I(X))$. 
			\item[\rm{(e)}]If $g\in A$, $V(I(V(g)))=V(g)$. Thus, if $X= V(g)$, then  $V(I(X))=X$.
			\item[\rm{(f)}] $I\left(  V\left(  I\left( X\right)  \right)  \right)  =I\left(
			X\right)$. 
			\item[\rm{(g)}]$I(\bigcup_{k\in \mathcal{K}}X_k)=\bigcap_{k\in \mathcal{K}}I(X_k)$.
			\item[\rm{(h)}]Let $Z:=(z_1,\dots,z_n)\in R^n$. Then, $I(\{Z\})=\langle Z\rangle$. 
			\end{enumerate}
	\end{enumerate}
\end{theorem}
\begin{proof}
	(i) (a) We have $f,g\in \langle Z\rangle$, so $f+g\in \langle Z\rangle$, i.e., $(f+g)(Z)=0$.
	
(b)	Let $Z\in V(f)$, then $f\in \langle Z\rangle$, then $gfh\in \langle Z\rangle$, i.e., $Z\in V(gfh)$.

	(ii) It is clear that $V(I)\subseteq V(g)$. From (i)-(b) we get that $V(g)\subseteq V(I)$.
	
	(iii) (a) It is clear that $V(0)=R^n$. 
	
	(b) Observe that for any $Z:=(z_1,\dots,z_n)\in R^n$, $Z\in V(x_1-z_1+\cdots+x_n-z_n)$. Thus, we have the nonzero polynomial $g:=x_1-z_1+\cdots+x_n-z_n$ and $\emptyset\subseteq V(g)$.
	
	(c) Evident.
	
	(d) Since $S\subseteq AS$, then $V(AS)\subseteq V(S)$; let $Z\in V(S)$ and $g\in AS$, then $g=p_1s_1+\cdots+p_ts_t$, with $p_i\in A$ and $s_i\in S$, $1\leq i\le t$. Since $p_is_i\in \langle Z\rangle$, then $g\in \langle Z\rangle$, so $V(S)\subseteq V(AS)$. Similarly, $V(S)=V(SA)=V(ASA)$.
	
	(e) Since $I\cap J\subseteq I,J$, then $V(I)\cup V(J)\subseteq V(I\cap J)$.
	
	(f) Since $I_k\subseteq \sum_{k\in \mathcal{K}}I_k$ for every $k\in \mathcal{K}$, then $V(\sum_{k\in \mathcal{K}}I_k)\subseteq \bigcap_{k\in \mathcal{K}}V(I_k)$. Let $Z\in \bigcap_{k\in \mathcal{K}}V(I_k)$ and let $g\in \sum_{k\in \mathcal{K}}I_k$, then $g=g_{k_1}+\cdots+g_{k_t}$, with $g_{k_j}\in I_{k_j}$, $1\leq j\leq t$, then from (i)-(a) , $g(Z)=0$, whence $Z\in V(\sum_{k\in \mathcal{K}}I_k)$. Thus, $\bigcap_{k\in \mathcal{K}}V(I_k)\subseteq V(\sum_{k\in \mathcal{K}}I_k)$.
	
	(h) Evident. 	
	
	(iv) (a)-(c) are evident from the definitions.
	
	(d) For $X=\emptyset$ the assertion follows from (a) since $\emptyset\subseteq V(A)$. Let $X\neq \emptyset$. If $Z\in X$, then for every $g\in I(X)$, $g(Z)=0$, and this means that $Z\in V(I(X))$. Therefore, $X\subseteq V(I(X))$. 
			
	(e) From (d), $V(g)\subseteq V(I(V(g)))$. Let $Z\in V(I(V(g)))$, since $g\in I(V(g))$, then $g(Z)=0$, i.e., $Z\in V(g)$. Therefore, $V(I(V(g)))\subseteq V(g)$. 
	
	(f) From (c), $I(X)\subseteq I(V(I(X)))$. From (d), $X\subseteq V(I(X))$, so from (b), $I( V(I(X))\subseteq I(X)$. 
	
	(g) Since $X_k\subseteq \bigcup_{k\in \mathcal{K}}X_k$ for every $k\in \mathcal{K}$, then $I(\bigcup_{k\in \mathcal{K}}X_k)\subseteq I(X_k)$, so $I(\bigcup_{k\in \mathcal{K}}X_k)\subseteq \bigcap_{k\in \mathcal{K}}I(X_k)$. Let $g\in \bigcap_{k\in \mathcal{K}}I(X_k)$ and let $Z\in \bigcup_{k\in \mathcal{K}}X_k$, then there exists $k\in \mathcal{K}$ such that $Z\in X_k$, then $g(Z)=0$, whence $g\in I(\bigcup_{k\in \mathcal{K}}X_k)$. 
	
	(h) Evident.
	\end{proof}

Definition \ref{definition4.14} and the previous theorem induces the following consequences.

\begin{corollary}\label{corollary4.5}
	\begin{enumerate}
		\item[\rm (i)]$R^n$ has a \textbf{Zariski topology} where the closed sets are the algebraic sets. 
		\item[\rm (ii)]If $X\subseteq R^n$ is finite, then $X$ is algebraic, and hence, closed.
			\end{enumerate}
\end{corollary}
\begin{proof}
	(i) By Definition \ref{definition4.14}, $R^n$ is algebraic. From Theorem \ref{theorem4.4} we know that $X=\emptyset $ is algebraic. Moreover, let $X\subseteq V(f)$ and $Y\subseteq V(g)$ be algebraic, with $0\neq f\in A$ and $0\neq g\in A$, then, since $A$ is a left noetherian domain, $A$ is a left Ore domain, i.e, $Af\cap Ag\neq 0$, whence 
	\begin{center} 
		$X\cup Y\subseteq V(f)\cup V(g)\subseteq V(l)$,
	\end{center}
	where $0\neq l\in Af\cap Ag$. Finally, let $\{X_k\}_{k\in \mathcal{K}}$ be a family of algebraic sets of $R^n$, then for every $k\in \mathcal{K}$ there exists $0\neq g_k\in A$ such that $X_k\subseteq V(g_k)$, hence
	\[
	\bigcap_{k\in \mathcal{K}} X_k\subseteq \bigcap_{k\in \mathcal{K}}V(g_k)=\bigcap_{k\in \mathcal{K}}V(Ag_k)=V(\sum_{k\in \mathcal{K}}Ag_k)\subseteq V(Ag_k)=V(g_k),\ \text{for any $k$}.
	\]
	
	(ii) We know that $X=\emptyset $ is algebraic. Let $\emptyset \neq X:=\{Z_1,\dots,Z_r\}$, then $I(X)\neq 0$. In fact, let
	\begin{center} 
		$f_i:=(x_1-z_{i1})+\cdots+(x_n-z_{in})$, with $Z_i:=(z_{i1},\dots,z_{in})$, $1\leq i\leq r$.
	\end{center}
	Let $0\neq l\in Af_1\cap \cdots \cap Af_r$. Observe that $l\in I(X)$: Indeed, for every $i$, $l=p_if_i$, for some $p_i\in A$, so $l(Z_i)=p_if_i(Z_i)=0$. This shows that $I(X)\neq 0$. 	
	Thus, $X\subseteq V(I(X))\subseteq V(l)$ is algebraic.
	
	\end{proof}

\begin{definition}
	Let $f\in A-R$. $V(f)$ is called the \textbf{skew hypersurface} defined by $f$. In particular,
	\begin{enumerate}
	\item[\rm (i)]$V(f)$ is a \textbf{skew plane curve} if $n=2$. 
	\item[\rm (ii)]$V(f)$ is a \textbf{skew hyperplane} if $\deg(f)=1$, i.e., $f=r_0+r_1x_1+\cdots +r_nx_n$, with $r_i\in A$, $0\leq i\leq n$. When $n=2$, $V(f)$
 is a \textbf{skew line}.
	\end{enumerate}
\end{definition}

\begin{corollary}\label{corollary4.20}
Let $I$ be a left ideal of $A$. Then,
\begin{center}
$V(I)=V(f_1)\cap \cdots \cap V(f_r)$, where $I=Af_1+\cdots+Af_r$. 
\end{center}
Thus, if $f_i\in A-R$, for $1\leq i\leq r$, then $V(I)$ is a finite intersection of skew hypersurfaces.
\end{corollary}
\begin{proof}
For $I=0$, $r=1$ and $f_1=0$. Let $I\neq 0$, since $A$ is left noetherian, there exist $f_1,\dots, f_r\in A$ such that $I=Af_1+\cdots+Af_r$. From Theorem \ref{theorem4.4},
\begin{center}
$V(I)=V(Af_1+\cdots+Af_r)=V(Af_1)\cap \cdots \cap V(Af_r)=V(f_1)\cap \cdots \cap V(f_r)$.
\end{center}    
\end{proof}

\begin{remark}\label{remark4.21}
	(i) There exist skew $PBW$ extensions such that $V(A)\neq \emptyset$. In fact, let $A:=\sigma(\mathbb{Q})\langle x,y,z\rangle$ defined by
	\begin{center}
		$yx=xy-1$, $zx=xz$, $zy=yz$. 
	\end{center}
	Consider the left ideal $I:=A(x-1)+Ay+Az$ and observe that
	\begin{center}
		$1=-y(x-1)+(x-1)y+0z=-y(x-1)+(x-1)(y-0)+0(z-0)$,
	\end{center}
	i.e., $I=A$ and $(1,0,0)\in V(1)=V(A)$.
	
	(ii) According to Corollary \ref{corollary4.5}, if $R$ is finite, then $R^n$ is algebraic, and hence we do not need to assume this condition on $R^n$ in Definition \ref{definition4.14}. But if $R$ is infinite, we can not assert that there is $0\neq g\in A$ such that $R^n\subseteq V(g)$. Consider for example that $A:=\mathbb{F}[x_1,\dots,x_n]$ is the commutative ring of polynomials with coefficients in an infinite field $\mathbb{F}$, then $I(\mathbb{F}^n)=0$ (see \cite{Fulton}, Chapter 1) and contrary assume that there exists $0\neq g\in A$ such that $\mathbb{F}^n\subseteq V(g)$, hence $I(V(g))\subseteq I(\mathbb{F}^n)=0$, but $g\in I(V(g))$, a contradiction.
	 \end{remark}
 
 We conclude this subsection with a result that partially generalizes Proposition \ref{proposition2.30}.  Recall that $f\in A$ is \textbf{\textit{normal}} if $Af=fA$, i.e., $f$ is a two-sided polynomial. 
 
 \begin{proposition}\label{proposition4.22}
 	Assume that $A$ is quasi-commutative.
 	\begin{enumerate}
 		\item[\rm (i)]Let $f=cx^{\alpha}h\in A$, where $c\in R^*$ is normal in R, $x^{\alpha}\in {\rm Mon}(A)$ and $h\in Z(A)$. Then, $f$ is normal.
 	 \item[\rm (ii)]Let $f=c_1x^{\alpha_1}+\cdots
 	+c_tx^{\alpha_t}\in A$, with $c_i\in R-\{0\}$, $1\leq i\leq t$, and $x^{\alpha_1}\succ \cdots \succ
 	x^{\alpha_t}$. If $f$ is normal, then $c_i$ is normal in $R$, for every $1\leq i\leq t$.
 	\end{enumerate}   
 \end{proposition}
 \begin{proof}
 	$\succeq$ is the ${\rm deglex}$ order on ${\rm Mon}(A)$, but any other monomial order could be used (for other monomial orders see \cite{Lezama-sigmaPBW}, Chapter 13).
 	
 	(i) Since the product of normal elements is normal and clearly  $h$ is normal, then we have to show only that $c$ and $x^{\alpha}$ are normal elements of $A$. Let $a:=a_1x^{\beta_1}+\cdots
 	+a_sx^{\beta_s}\in A$,  with $a_j\in R-\{0\}$, $1\leq j\leq s$, and $x^{\beta_1}\succ \cdots \succ
 	x^{\beta_s}$. We have $ac=(a_1x^{\beta_1}+\cdots
 	+a_sx^{\beta_s})c=a_1x^{\beta_1}c+\cdots
 	+a_sx^{\beta_s}c=a_1\sigma^{\beta_1}(c)x^{\beta_1}+\cdots+a_s\sigma^{\beta_s}(c)x^{\beta_s}=a_1cc^{-1}\sigma^{\beta_1}(c)x^{\beta_1}+\cdots +a_scc^{-1}\sigma^{\beta_s}(c)x^{\beta_s}=c(a_1c^{-1}\sigma^{\beta_1}(c)x^{\beta_1}+\cdots +a_sc^{-1}\sigma^{\beta_s}(c)x^{\beta_s})$, thus $Ac\subseteq cA$. Since $A$ is bijective, then we can prove similarly that $cA\subseteq Ac$. Now, $ax^{\alpha}=(a_1x^{\beta_1}+\cdots
 	+a_sx^{\beta_s})x^{\alpha}=a_1x^{\beta_1}x^{\alpha}+\cdots+a_sx^{\beta_s}x^{\alpha}=a_1c_1'x^{\alpha}x^{\beta_1}+\cdots+a_sc_s'x^{\alpha}x^{\beta_s}=x^{\alpha}\sigma^{-\alpha}
 	(a_1c_1')x^{\beta_1}+\cdots+x^{\alpha}\sigma^{-\alpha}(a_sc_s')x^{\beta_s}$, for some $c_j'\in R^*$, $1\leq j\leq s$, thus $Ax^{\alpha}\subseteq x^{\alpha}A$. In a similar way we can prove that $x^{\alpha}A\subseteq Ax^{\alpha}$.
 	
 	(ii) Let $r\in R-\{0\}$, then $rf\in Af=fA$, so $rf=fg$, for some $g\in A$. Since $A$ is a domain, $\deg(rf)=\deg(f)=\deg(fg)=\deg(f)+\deg(g)$, so $g\in R-\{0\}$, but as $A$ is quasi-commutative, for every $1\leq i\leq t$, $rc_i=c_i\sigma^{\alpha_i}(g)$. This means that $Rc_i\subseteq c_iR$. Considering now $fr\in fA=Af$ we get that $fr=hf$, with $h\in R-\{0\}$, so for every $i$, $c_i\sigma^{\alpha_i}(r)=hc_i$, but since $A$ is bijective, $\sigma^{\alpha_i}(R)=R$, and hence $c_iR\subseteq Rc_i$. Thus, $c_iR= Rc_i$, i.e., $c_i$ is normal in $R$.  
 \end{proof} 
  For the habitual polynomials in one variable, the converse of the part (i) of the previous proposition is true, as the following corollary shows. This corollary also complements the part (i) of Remark \ref{remark4.21}.
 \begin{corollary}\label{proposition4.23}
 	Let $S$ be a left noetherian domain and $B:=S[x]$ be the habitual ring of polynomials.
 	\begin{enumerate}
 		\item[\rm (i)]Let $f\in B$, with ${\rm lc}(f)\in S^*$. $f$ is a normal polynomial if and only if $f=cx^th$, where $c\in S^*$ is normal, $t\geq 0$ and $h\in Z(B)$.
 		\item[\rm (ii)]Let $\mathbb{F}$ be an algebraically closed field and assume that $S$ is an $\mathbb{F}$-algebra with trivial center.
 		Let $I=Bf_1+\cdots+Bf_r$ be a left ideal of $B$, where $0\neq f_i$ is normal and ${\rm lc}(f_i)\in S^*$, for every $1\leq i\leq r$. If $V(I)=\emptyset$, then $I=B$.  
 	\end{enumerate}
 \end{corollary}
 \begin{proof}
 	Notice first that $B$ is a quasi-commutative bijective skew $PBW$ extension of $S$, so we can use all the previous results.
 	
 	(i) $\Rightarrow )$ Let $f:=f_0+f_1x\cdots +f_nx^n\neq 0$, with $f_n\neq 0$. We can assume that $f$ is monic. In fact, $f=f_n(f_n^{-1}f_0+f_n^{-1}f_1x+\cdots+x^n)=f_nf'$, with $f':=f_n^{-1}f_0+f_n^{-1}f_1x+\cdots+x^n$. We have $Bf=fB$, but from Proposition \ref{proposition4.22}, $f_n$ is a normal element of $S$, so  $Bf_n=f_nB$, and hence $Bf_nf'=f_nf'B=f_nBf'$, but since $B$ is a domain, then $f'B=Bf'$. Therefore, if we prove the claimed for monic polynomials, then $f'=c'x^th$, with $c'\in S^*$ normal, $t\geq 0$ and $h\in Z(B)$, so $f=cx^th$, with $c:=f_nc'\in S^*$ normal.  
 	
 	We will prove the claimed by induction on $n$. If $n=0$, then $c=f=1\in S^*$ normal, $t=0$ and $h=1$. If $n=1$, then $f=f_0+x$; if $f_0=0$, then we get the claimed with $c=1$, $t=1$ and $h=1$. Assume that $f_0\neq 0$; let $s\in S$, then $fs=bf$, with $b\in B$, this implies that $b:=b_0\in S$, and hence $f_0s+sx=b_0f_0+b_0x$, whence $s=b_0$ and $f_0s=sf_0$, i.e., $f_0\in Z(S)$. Therefore, $f\in Z(B)$ and we get the claimed with $c=1$, $t=0$ and $h=f$. This completes the proof of case $m=1$.
 	
 	Assume the claimed proved for non zero monic normal polynomials of degree $\leq n-1$. We will consider two possible cases. 
 	
 	\textit{Case 1}. $f_0\neq 0$. As before, let $s\in S$, then $fs=bf$, with $b\in B$, but this implies that $b:=b_0\in S$ and $f_is=sf_i$ for every $0\leq i\leq n-1$. Thus, $f\in Z(B)$ and we get the claimed with $c=1$, $t=0$ and $h=f$.
 	
 	\textit{Case 2}. $f_0=0$. Let $l$ be minimum such that $f_l\neq 0$. Then $f=x^{l}f'$, with $l\geq 1$ and $0\neq f'\in A$ monic. Observe that $x^l$ and $f'$ are normal: In fact, $x^lB=Bx^l$ and $Bf=Bx^lf'=x^lBf'=x^lf'B=fB$, but since $B$ is a domain, then $Bf'=f'B$. By induction, there exist $c'\in S^*$ normal, $l'\geq 0$ and $h\in Z(B)$ such that $f'=c'x^{l'}h$, hence $f=x^{l}c'x^{l'}h=c'x^{l+l'}h$, so we obtain the claimed with $c:=c'$ and $t:=l+l'$.
 	
 	$\Leftarrow )$ This follows from the previous proposition. 
 	
 	(ii) From Corollary \ref{corollary4.20} we have that $V(I)=V(f_1)\cap \cdots \cap V(f_r)$, but by (i), for every $1\leq i\leq r$, $f_i=c_ix^{t_i}h_i$, where $c_i\in S^*$ is normal, $t_i\geq 0$ and $h_i\in Z(B)=Z(S)[x]=\mathbb{F}[x]$. Since $c_i\in S^*$, then $V(f_i)=V(x^{t_i}h_i)$, hence
 	\begin{center} 
 		$V(I)\supseteq V(x^{t_1}h_1)\cap \cdots \cap V(x^{t_r}h_r)\supseteq V(d)$, where $d:={\rm gcrd}(x^{t_1}h_1,\dots,x^{t_r}h_r)$ in $\mathbb{F}[x]$. 
 	\end{center}
 	Since $V(I)=\emptyset$, then $V(d)=\emptyset$ with respect to $R^n$, whence, $V(d)=\emptyset$ with respect to $\mathbb{F}^n$, but $\mathbb{F}$ is algebraically closed, then $d\in \mathbb{F}^*$. We have $d=g_1x^{t_1}h_1+\cdots+g_rx^{t_r}h_r$, for some $g_1,\dots,g_r\in \mathbb{F}[x]$. Let $c:=c_1\cdots c_r\in S^*$, since $c_i$ is normal and $g_i\in Z(B)$, for every $1\leq i\leq r$, we get $cd=g_1'c_1x^{t_1}h_1+\cdots+g_r'c_rx^{t_r}h_r$, where every $g_i'\in B$. Thus, $cd=g_1'f_1+\cdots+g_r'f_r\in I\cap S^*$, i.e., $I=B$.  
 \end{proof}

\subsection{Hilbert's Nullstellensatz theorem for skew $PBW$ extensions}

The ring-theoretic version of the Hilbert's Nullstellensatz theorem for skew $PBW$ extensions has been considered in the beautiful paper \cite{reyes-jason}. The algebraic characterization of the theorem given by the authors of \cite{reyes-jason} does not use the notion of variety (see Theorem 3.1 in \cite{reyes-jason}). Applying the algebraic sets and the ideal of points introduced in Definition \ref{definition4.14} and Theorem \ref{theorem4.4}, we present next the classical version of this important theorem for quasi-commutative bijective skew $PBW$ extensions of algebraically closed fields. Our version covers the Nullstellensatz theorem of commutative algebraic geometry (see \cite{Fulton}, Chapter 1).

We start recalling some notions and results related to prime ideals of an arbitrary ring (see \cite{McConnell} and also \cite{Birkenmeier}, Definition 3).

\begin{definition}
	Let $S$ be a ring and $I,P$ be two-sided ideals of $S$, with $P\neq S$.
\begin{enumerate}
\item[\rm (i)]
$P$ is a \textbf{prime ideal} of $S$ if for any left ideals $L,J$ of $S$ the following condition holds:
$LJ\subseteq P$ if and only if $L\subseteq P$ or $J\subseteq P$.
\item[\rm (ii)]The \textbf{radical} of $I$, denoted $\sqrt{I}$, is the intersection of all prime ideals of $S$ containing $I$.
\item[\rm (iii)]An element $a\in S$ is \textbf{$I$-strongly nilpotent} if for any given sequence $\mathcal{S}:=\{a_i\}_{i\geq 1}$ of elements of $S$, with $a_1:=a$ and $a_{i+1}\in a_iSa_i$, there exists $m(\mathcal{S})\geq 1$ such that $a_{m(\mathcal{S})}\in I$. We say that $a$ is \textbf{$I$-nilpotent} if there exists $m\geq 1$ such that $a^m\in I$.
\item[\rm (iv)]$P$ is \textbf{completely prime} if the following condition holds for any $a,b\in S$: $ab\in P$ if and only if $a\in P$ or $b\in P$.  
\item[\rm (v)]$P$ is \textbf{completely semiprime} if the following condition holds for any $a\in S$: $a^2\in P$ if and only if $a\in P$. 
\end{enumerate}
\end{definition}
It is clear that if $a\in S$ is $I$-strongly nilpotent, then $a$ is $I$-nilpotent. If $a\in Z(S)$, then the converse is true. Observe that any element $a\in S$ is $S$-strongly nilpotent and $\sqrt{S}:=S$. If $P$ is completely prime, then $P$ is completely semiprime. By induction on $m$ it is easy to show that $P$ is completely semiprime if and only if the following condition holds: For any $a\in S$ and any $m\geq 1$,  $a^m\in P$ if and only if $a\in P$.

\begin{proposition}\label{proposition4.14}
	Let $S$ be a ring and $P$ be a proper two-sided ideal of $S$. $P$ is a prime ideal if and only if the following condition holds for any elements $a,b\in S$: $aSb\subseteq P$ if and only if $a\in P$ or $b\in P$.
\end{proposition}
\begin{proof}
$\Rightarrow)$: From $aSb\subseteq P$ we get that $SaSb\subseteq P$, whence, $Sa\subseteq P$ or $Sb\subseteq P$, i.e., $a\in P$ or $b\in P$. 

$\Leftarrow)$: Let $L,J$ be left ideals of $S$ such that $LJ\subseteq P$. Assume that $L\nsubseteq P$ and let $a\in L$ with $a\notin P$. Let $b\in J$, then $aSb\subseteq LJ\subseteq P$, whence $b\in P$. Thus, $J\subseteq P$. 
\end{proof}

\begin{proposition}\label{proposition4.15}
Let $S$ be a ring and $P$ be a proper two-sided ideal of $S$. $P$ is a prime ideal if and only if the following condition holds for any left ideals $L,J$ of $S$: If $P\subsetneq L$ and  $P\subsetneq J$, then $LJ\nsubseteq P$.
\end{proposition}
\begin{proof}
$\Rightarrow)$: Evident.

$\Leftarrow)$: Let $a,b\in S$ such that $aSb\subseteq P$ and suppose that $a\notin P$ and $b\notin P$. Then, $P\subsetneq P+Sa$ and $P\subsetneq P+Sb$, and from the hypothesis, $(P+Sa)(P+Sb)\nsubseteq P$, hence, there exist $p,p'\in P$ and $s,s'\in S$ such that $(p+sa)(p'+s'b)\notin P$, a contradiction.
\end{proof}

\begin{proposition}
	Let $S$ be a ring and $I$ be a two-sided ideal of $S$. Then,
\begin{center}
$\sqrt{I}=\{a\in S\mid a\ \text{is $I$-strongly nilpotent}\}$.
\end{center}
\end{proposition}
\begin{proof}
Let $a\in S$ such that $a\notin \sqrt{I}$, then there exists a prime ideal $P$ of $S$, containing $I$, such that $a\notin P$, hence, by Proposition \ref{proposition4.14},  $aSa\nsubseteq P$. This says that there exists $b\in S$ such that $aba\notin P$. Let $a_1:=a$ and $a_2:=aba$. Thus, $a_2Sa_2\nsubseteq P$ and hence there exists $c\in S$ such that $a_2ca_2\notin P$. Let $a_3:=a_2ca_2$. Continuing this way we get a sequence $\{a_i\}_{i\geq 1}$ of elements of $S$ such that $a_i\notin P$ for every $i\geq 1$, hence, $a_i\notin I$ for every $i\geq 1$. This means that $a$ is not $I$-strongly nilpotent.

Conversely, assume that $a\in S$ is not $I$-strongly nilpotent, then there exists a sequence $\mathcal{S}:=\{a_i\}_{i\geq 1}$ of elements of $S$, with $a_1:=a$ and $a_{i+1}\in a_iSa_i$, such that for every $i\geq 1$, $a_i\notin I$. By Zorn's lemma, there exists a two-sided ideal $P$ of $S$, containing $I$, maximal with respect to the condition $\mathcal{S}\cap P=\emptyset$ (observe that $I\supseteq I$ and $\mathcal{S}\cap I=\emptyset$). We will show that $P$ is a prime ideal of $S$. It is clear that $P\neq S$. Let $L,J$ be left ideals of $S$ such that $P\subsetneq L$ and $P\subsetneq J$. Then, $L\cap \mathcal{S}\neq \emptyset$ and $J\cap \mathcal{S}\neq \emptyset$, so there exists $a_i\in L$ and $a_j\in J$. Let $k:=\max\{i,j\}$, then $a_{k+1}\in LJ$, but $a_{k+1}\notin P$, i.e., $LJ\nsubseteq P$. Thus, from Proposition \ref{proposition4.15}, $P\supseteq I$ is a prime ideal such that $a\notin P$, so $a\notin \sqrt{I}$. 
\end{proof}

\begin{lemma}\label{lemma4.15}	Let $A:=\sigma(\mathbb{F})\langle x_1,\dots,x_n\rangle$ be a quasi-commutative bijective skew $PBW$ extension of $\mathbb{F}$, where $\mathbb{F}$ is a field. Then, for every $Z:=(z_1,\dots,z_n)\in \mathbb{F}^n$, $\langle Z\rangle$ is completely semiprime. 
\end{lemma}
\begin{proof}
We have to show first that $\langle Z\rangle\neq A$: Contrary, assume that $1\in A$, then 
\begin{center}
$1=p_1(x_{i_1}-z_{i_1})q_1+\cdots+p_m(x_{i_m}-z_{i_m})q_m$, with $p_j,q_j\in A$ and $i_j\in \{1,\dots,n\}$, $1\leq j\leq m$.
\end{center}
By Proposition \ref{algdivforPBW},
\begin{center} 
	$q_j=q_j'(x_{i_j}-z_{i_j})+h_j$, with $q_j',h\in A$, $h_j$ reduced w.r.t. $x_{i_j}-z_{i_j}$, $1\leq j\leq m$.
\end{center}
So, 
\begin{center}
	$1=p_1(x_{i_1}-z_{i_1})(q_1'(x_{i_1}-z_{i_1})+h_1)+\cdots+p_m(x_{i_m}-z_{i_m})(q_m'(x_{i_m}-z_{i_m})+h_m)=p_1(x_{i_1}-z_{i_1})q_1'(x_{i_1}-z_{i_1})+p_1(x_{i_1}-z_{i_1})h_1+
	\cdots+p_m(x_{i_m}-z_{i_m})q_m'(x_{i_m}-z_{i_m})+p_m(x_{i_m}-z_{i_m})h_m$.
\end{center}

Since every $h_j$ is reduced, then $h_j\in \mathbb{F}$, but $A$ is quasi-commutative, then $1\in I$, where $I$ is a left ideal of $A$ generated by elements of the form  $c_{i_j}x_{i_j}-z_{i_j}'$, where $c_{i_j},z_{i_j}'\in \mathbb{F}$ with $c_{i_j}\neq 0$, $1\leq j\leq m$ (actually, in some cases $c_{i_j}=1$ and $z_{i_j}'=z_{i_j}$, and in other cases, when $h_j\neq 0$, then $c_{i_j}=\sigma_{i_j}(h_j)$ and $z_{i_j}'=z_{i_j}h_j$, where $\sigma_{i_j}$ is as in Proposition \ref{sigmadefinition}). It is clear from Definition \ref{definition3.10} that the generators $c_{i_j}x_{i_j}-z_{i_j}'$ of $I$ conform a Gröbner basis of $I$. Since $1\in I$, from the part (iii) of Proposition \ref{153} we get that $x_{i_j}$ divides $1$ for some $j$, a contradiction. Hence, $\langle Z\rangle\neq A$.     

Now, let $f\in A$ such that $f^2\in \langle Z\rangle$. Applying again Proposition \ref{algdivforPBW}, there exist
polynomials $p_1,\dots ,p_t,h\in A$, with remainder $h$ reduced w.r.t.\ $F:=\{x_1-z_1,\dots,x_n-z_n\}$, such that
\begin{center}
$f=p_1(x_1-z_1)+\cdots +p_n(x_n-z_n)+h$.
\end{center}
Since $h$ is reduced, then $h\in \mathbb{F}$. If $h=0$, then $f\in \langle Z\rangle$ and the proof is over. Assume that $h\neq 0$, then
\begin{center} 
$h^2=(f-[p_1(x_1-z_1)+\cdots +p_n(x_n-z_n)])^2\in \langle Z\rangle$,
\end{center} 
hence $\langle Z\rangle=A$, a contradiction.
\end{proof} 

\begin{theorem}[\textbf{Hilbert's Nullstellensatz}]\label{Nullstellensatz}
Let $A:=\sigma(\mathbb{F})\langle x_1,\dots,x_n\rangle$ be a quasi-commutative bijective skew $PBW$ extension of $\mathbb{F}$, where $\mathbb{F}$ is an algebraically closed field. Assume that $Z(A)$ is a polynomial ring in $n$ variables with coefficients in $\mathbb{F}$.  
Let $I$ be a two-sided ideal of $A$. Then,  
	\begin{center}
		$\langle I_{Z(A)}(V_{Z(A)}(J))\rangle\subseteq \sqrt{I}\subseteq I(V(I))$,
	\end{center}
where $J:=l^{-1}(I)$, $l:Z(A)\to A$ is the inclusion of the center of $A$ in $A$, $V_{Z(A)}(J)$ is the vanishing set of $J$ with respect to $Z(A)$ and $I_{Z(A)}(V_{Z(A)}(J))$ is the ideal of points of $V_{Z(A)}(J)$ with respect to $Z(A)$. 
\end{theorem}
\begin{proof}
	
$\sqrt{I}\subseteq I(V(I))$: If $V(I)=\emptyset$, then $I(V(I))=A$ and hence there is nothing to prove. Assume that $V(I)\neq \emptyset$ and let $f\in \sqrt{I}$, then $f$ is $I$-strongly nilpotent, and hence, $I$-nilpotent, so there exists $m\geq 1$ such that $f^m\in I$. Let $Z:=(z_1,\dots,z_n)\in V(I)$, then $f^m\in \langle Z\rangle$. From Lemma \ref{lemma4.15}, $f\in \langle Z\rangle$, i.e., $f\in I(V(I))$. 

$\langle I_{Z(A)}(V_{Z(A)}(J))\rangle\subseteq \sqrt{I}$: Consider the inclusion $Z(A)\xrightarrow{l} A$ and let $J:=l^{-1}(I)$. Then, $J=l(J)=l(l^{-1}(I))\subseteq I$, and let $\langle J\rangle:=AJA$ be the two-sided ideal of $A$ generated by $J$. We have $\langle J\rangle\subseteq I$, so $\sqrt{\langle J\rangle}\subseteq \sqrt{I}$, but $\langle \sqrt{J}\rangle\subseteq \sqrt{\langle J\rangle}$, where $\sqrt{J}$ is the radical of $J$ in the ring $Z(A)$. In fact, let $w\in \sqrt{J}$, then there exists $m\geq 1$ such that $w^m\in J\subseteq \langle J\rangle$, but $w\in Z(A)$, then $w$ is $\langle J\rangle$-strongly nilpotent, i.e., $w\in \sqrt{\langle J\rangle}$. Thus, $\langle \sqrt{J}\rangle\subseteq \sqrt{I}$. Applying the classical Hilbert's Nullstellensatz for $Z(A)$ (here we use that $\mathbb{F}$ is algeraically closed) we have $\sqrt{J}=I_{Z(A)}(V_{Z(A)}(J))$, so we get that 
$\langle I_{Z(A)}(V_{Z(A)}(J))\rangle\subseteq \sqrt{I}$. 
\end{proof}

\begin{example}
Next we present some concrete examples of skew $PBW$ extensions that satisfy the hypothesis of Theorem \ref{Nullstellensatz}. $\mathbb{F}$ denotes an algebraically closed field.

(i) It is clear that if $A=\mathbb{F}[x_1,\dots,x_n]$ and $I$ is an ideal of $A$, then in Theorem \ref{Nullstellensatz} we have $\langle I_{Z(A)}(V_{Z(A)}(J))\rangle=I(V(I))$, and hence, $I(V(I))=\sqrt{I}$. 

(ii) If $q\neq 1$ is an arbitrary root of unity of degree $m\geq 2$, then the center of the quantum
plane $A:=\mathbb{F}_{q}[x,y]$ is the subalgebra generated by $x^m$ and $y^m$, i.e.,
$Z(\mathbb{F}_{q}[x,y])=\mathbb{F}[x^m,y^m]$ (see \cite{Shirikov} or also \cite{Lezama-sigmaPBW}, Proposition 3.3.14). Recall that the rule of multiplication in $A$ is given by $yx=qxy$.

(iii) The previous example can be generalized in the following way (see \cite{Zhangetal}, Lemma 4.1, or also \cite{Lezama-sigmaPBW}, Proposition 3.3.15): Let $q\in \mathbb{F}-\{0\}$ and $A:=\mathbb{F}_q[x_1,\dots,x_n]$ be the skew $PBW$ extension
defined by $x_jx_i=qx_ix_j$ for all $1\leq i<j\leq n$. If $n\geq 2$ and $q\neq 1$ is a root of
unity of degree $m\geq 2$, then
\begin{enumerate}
	\item[\rm (a)]If $q=-1$, then
	\begin{center}
		$Z(A)=\mathbb{F}[x_1^2,\dots,x_n^2]$ when $n$ is even.
		
\end{center}
	\item[\rm (b)]If $q\neq -1$, then
	\begin{center}
		$Z(A)=\mathbb{F}[x_1^m,\dots,x_n^m]$ when $n$ is even.
		
	\end{center}
\end{enumerate}

(iv) Consider that for every $1\leq i,j\leq n$, $q_{ij}\in \mathbb{F}-\{0\}$ is a nontrivial
root of unity of degree $d_{ij}<\infty$ and let $A:=\mathbb{F}_{\boldsymbol{\rm q}}[x_1,\dots,x_n]$ be the skew $PBW$ extension
	defined by $x_jx_i=q_{ij}x_ix_j$ for all $1\leq i<j\leq n$. Let $k_{ij}\in \mathbb{Z}$ such that $|k_{ij}|<d_{ij}$,
${\rm lcd}(k_{ij}, d_{ij})=1$ and $q_{ij}=\exp(2\pi\sqrt{-1}\frac{k_{ij}}{d_{ij}})$ $($choosing
$k_{ji}:=-k_{ij}$$)$. Let $L_i:={\rm lcm}\{d_{ij}|j=1,\dots,n\}$. Then 
$Z(\mathbb{F}_{\boldsymbol{\rm q}}[x_1,\dots,x_n])$ is a polynomial ring if and only if it is of the form
$\mathbb{F}[x_1^{L_1},\dots,x_n^{L_n}]$ (see \cite{Zhangetal3}, Theorem 0.3 or also \cite{Lezama-sigmaPBW}, Proposition 3.3.17).
\end{example}

\begin{remark}
In this last section we have studied only the algebraic sets of skew $PBW$ extensions and the results have been published separately in \cite{LezamaNullstellensatz}. The other properties of noncommutative coding theory considered in the second section for Ore extensions have not been investigated in this paper. Probably this represents a challenge for our readers.   
\end{remark}
%\begin{center}
%\textbf{Acknowledgements}
%\end{center}

%%%%%%%%%%%%%%%%%%%%%%%%%%%%%%%%%%%%%%%%%%%%%%%%%
%%%%%%%%%%%%%%%%%%%%%%%%%%%%%%%%%%%%%%%%%%%%%%%%%

\end{document}